\documentclass[leqno,11pt]{amsart}
\usepackage{amsmath,amscd,amsthm}
\usepackage{epsfig,graphics,ifsym,color}
\usepackage{amssymb,latexsym}
\usepackage{mathrsfs}
\usepackage{hyperref}
\usepackage[all,cmtip]{xy}

\newtheorem{thm}{\bf Theorem}[section]
\newtheorem{df}[thm]{\bf Definition}
\newtheorem{prop}[thm]{\bf Proposition}
\newtheorem{cor}[thm]{\bf Corollary}
\newtheorem{lem}[thm]{\bf Lemma}
\newtheorem{rem}[thm]{\bf Remark}
\newtheorem{ex}[thm]{\bf Example}

\newtheorem{nono-theorem}{Theorem}[]
\newtheorem*{thm*}{Theorem}

\newcommand{\A}{\mathscr{A}}

\newcommand{\B}{\mathbf{B}}

\newcommand{\cP}{\mathscr{P}}

\newcommand{\pf}{\noindent{\bfseries Proof. }}
\newcommand{\ov}{\overline}

\newcommand{\F}{\mathscr{F}}

\newcommand{\hf}{\frac{1}{2}}

\newcommand{\gl}{\mathfrak{gl}}
\newcommand{\Z}{\mathbb{Z}}
\newcommand{\C}{\mathbb{C}}
\newcommand{\h}{\mathfrak{h}}
\newcommand{\te}{\widetilde{e}}
\newcommand{\tf}{\widetilde{f}}
\newcommand{\g}{\mathfrak{g}}
\newcommand{\td}{\widetilde}

\newcommand{\mc}{\mathcal}
\newcommand{\mf}{\mathfrak}

\newcommand{\I}{\mathbb{I}}
\newcommand{\J}{\mathbb{J}}
\newcommand{\psp}{\psi}
\newcommand{\psm}{\psi^\ast}

\newcommand{\om}{\omega}

\numberwithin{equation}{section}

\begin{document}
\title[Crystal bases for quantum ortho-symplectic superalgebras]
{Super duality and Crystal bases for quantum ortho-symplectic superalgebras}
\author{JAE-HOON KWON}
\address{Department of Mathematics \\ Sungkyunkwan University \\ Suwon,  Republic of Korea}
\email{jaehoonkw@skku.edu}

\thanks{This work was  supported by Basic Science Research Program through the National Research Foundation of Korea (NRF)
funded by the Ministry of  Education, Science and Technology (No. 2011-0006735).}

\begin{abstract}
We introduce a semisimple tensor category $\mc{O}^{int}_q(m|n)$ of modules over an quantum ortho-symplectic superalgebra. It is a natural counterpart of the category  of finitely dominated integrable modules over the quantum classical (super) algebra of type $B_{m+n}$, $C_{m+n}$, $D_{m+n}$ or  $B(0,m+n)$ from a viewpoint of super duality. 
We classify the irreducible modules in $\mc{O}^{int}_q(m|n)$ and show that an irreducible module in $\mc{O}^{int}_q(m|n)$ has a unique crystal base in case of type $B$ and $C$. An explicit description of the crystal graph is given in terms of  a new combinatorial object called ortho-symplectic tableaux.
\end{abstract}

\maketitle
\setcounter{tocdepth}{1}
\tableofcontents

\section{Introduction}
\subsection{}
The  Kashiwara's crystal base \cite{Kas1} is a certain nice basis at $q=0$ of a module $M$ over the quantized enveloping algebra associated to a symmetrizable Kac-Moody algebra, which still contains rich combinatorial information on $M$, and it has been one of the most important and successful tools in representation theory of the quantum groups.

For the quantum superalgebra $U_q({\mf g})$ associated to a classical Lie superalgebra ${\mf g}$, the existence of a crystal base  is shown when $\mf{g}$ is a general linear Lie superalgebra $\mf{gl}_{m|n}$ \cite{BKK} and a queer Lie superalgebra $\mf{q}_n$ \cite{GJKKK-1,GJKKK-2} for a special but important class of finite dimensional modules appearing in a tensor power of the natural representation $V$ of $U_q(\mf{g})$, often called polynomial representations. We should remark that the crystal base theory for these two classical Lie superalgebras is far from being parallel to that of a symmetrizable Kac-Moody algebra due to the same substantial difficulties encountered when we consider the representations of classical Lie superalgebras compared to those of Lie algebras.  For example, a finite dimensional representation of $\mf{g}$ is not semisimple in general. Indeed, the above results for $\mf{g}=\mf{gl}_{m|n}$ and $\mf{q}_n$ are  based on the well-known semisimplicity of $V^{\otimes r}$, and closely related with the Schur-Weyl-Sergeev dualities \cite{Ser}. Also, there is a work on crystal bases of a family of infinite dimensional representations of the exceptional Lie superalgebra of type $D(2|1,\alpha)$ \cite{Zou}.

There is another important class in classical Lie superalgebras called ortho-symplectic Lie superalgebras $\mf{osp}_{k|2l}$. However, there is little known about the existence of  its crystal bases except for $\mf{osp}_{1|2l}$, which is of type $B(0,l)$, and  where we can apply the crystal base theory for a Kac-Moody superalgebra without odd isotropic simple root \cite{Je}. 

In this paper, we construct for the first time crystal bases of a large family of semisimple modules over a quantum ortho-symplectic superalgebra. We also prove the uniqueness of these crystal bases and give a combinatorial model for the associated crystals.  
\subsection{}
Let us explain our results in more details. Our first step is to find a nice semisimple category of modules over a quantum ortho-symplectic superalgebra. Since a tensor power of the natural representation of an ortho-symplectic Lie superalgebra is not semisimple in general, we take a completely different approach motivated by a recent work of Cheng, Lam and Wang on super duality \cite{CLW}.

Super duality is an equivalence between a parabolic BGG category $\mc{O}(m+\infty)$ of modules over the classical Lie algebras $\mf{g}_{m+\infty}$ and a category $\mc{O}(m|\infty)$ of modules over the basic classical Lie superalgebras $\mf{g}_{m|\infty}$ of infinite rank, where $\mf{g}=\mf{gl}, \mf{b}, \mf{b}^\bullet, \mf{c}, \mf{d}$. (From now on, we use $\mf{g}$ as a symbol representing the type of a Lie superalgebra.)
It was originally introduced in \cite{CW08,CWZ} as a conjecture in case of general linear Lie superalgebras and later proved by Cheng and Lam \cite{CL}. Then the duality for ortho-symplectic Lie superalgebras was established by Cheng, Lam and Wang \cite{CLW}.
One of its most remarkable and powerful features  is that super duality reveals a natural connection with the Kazhdan-Lusztig theory of Lie algebras.

We consider the semisimple subcategory  $\mc{O}^{int}(m|\infty)$ of $\mc{O}(m|\infty)$ equivalent to the subcategory $\mc{O}^{int}(m+\infty)$ of integrable modules in $\mc{O}(m+\infty)$ under super duality. It is known that $\mc{O}^{int}(m|\infty)$ is the category of polynomial modules when $\mf{g}_{m|\infty}$ is a general linear Lie superalgebra, that is,  $\mf{g}=\mf{gl}$, \cite{CW}. Motivated by this fact, we prove that when $\mf{g}_{m|\infty}$ is ortho-symplectic, that is, $\mf{g}=\mf{b}, \mf{b}^\bullet, \mf{c}, \mf{d}$, $\mc{O}^{int}(m|\infty)$ is a full subcategory of $\mc{O}(m|\infty)$ such that the weights of each object are polynomial with respect to a suitably chosen dual basis of the Cartan subalgebra (Theorem \ref{complete reducibility for integrable super}). We have a similar result for a category $\mc{O}^{int}(m|n)$ of modules over $\mf{g}_{m|n}$ of finite rank, where $\mc{O}^{int}(m|n)$ is obtained from $\mc{O}^{int}(m|\infty)$ by applying the truncation functor (Theorem \ref{complete reducibility for integrable super of finite rank}).
We should note that $\mc{O}^{int}(m|n)$ is not characterized only by locally nilpotent actions of positive simple root vectors since the odd isotropic root vectors  are always nilpotent on $\mf{g}_{m|n}$-modules.

Unlike the case of $\mf{g}=\mf{gl}$, the irreducible modules in $\mc{O}^{int}(m|n)$ are infinite dimensional when $\mf{g} =\mf{b}, \mf{b}^\bullet, \mf{c}, \mf{d}$ and $n>0$, which were called oscillator modules and studied via Howe duality in \cite{CKW}. But one may still regard $\mc{O}^{int}(m|n)$ as a natural counterpart of the category $\mc{O}^{int}(m+n)$ of finite dimensional modules over $\mf{g}_{m+n}$ of finite rank, since both of them are obtained from two equivalent categories $\mc{O}^{int}(m|\infty)$ and $\mc{O}^{int}(m+\infty)$ by truncation, respectively (see the diagram below, where $\mc{F}$ is the super duality functor and $\mf{tr}_n$ is a truncation functor).
\begin{equation*}
\begin{split}
&\ \  \xymatrixcolsep{3pc}\xymatrixrowsep{.3pc}\xymatrix{
\mc{O}(m+\infty)\  \ar@{->}[r]^-{\sim}_-{\mc{F}}
& \ \ \mc{O}(m|\infty) \\ \!\!\bigcup & \ \,\bigcup}   \\
&\xymatrixcolsep{3pc}\xymatrixrowsep{2pc}\xymatrix{
\mc{O}^{int}(m+\infty) \ar@{->}[d]^{\mf{tr}_n} \ar@{->}[r]^-{\sim}_-{\mc{F}}
& \mc{O}^{int}(m|\infty)  \ar@{->}[d]^{\mf{tr}_n} \\
\mc{O}^{int}(m+n)
& \mc{O}^{int}(m|n)}
\end{split}
\end{equation*}

Now we consider $q$-deformations of $\mf{g}_{m|n}$-modules in $\mc{O}^{int}(m|n)$  for $\mf{g} =\mf{b}, \mf{b}^\bullet, \mf{c}, \mf{d}$. Based on our characterization of $\mc{O}^{int}(m|n)$, we define a category  $\mc{O}^{int}_q(m|n)$ of $U_q(\mf{g}_{m|n})$-modules with the same conditions on weights (Definition \ref{category O^int_q(m|n)}). We can check by using the method of classical limit that  $\mc{O}^{int}_q(m|n)$ is a semisimple tensor category (Theorem \ref{complete reducibility}), and define the notion of a crystal base of a module in $\mc{O}^{int}_q(m|n)$ according to \cite{BKK}. 

We parametrize  the  irreducible highest weight modules in $\mc{O}^{int}(m|n)$ by $\cP(\mf{g})_{m|n}$, which is a set of pairs $(\lambda,\ell)$ of a partition and a positive integer with certain conditions, and let  $\Lambda_{m|n}(\lambda,\ell)$ be the corresponding highest weight. When $\mf{g}=\mf{b}, \mf{b}^\bullet, \mf{c}$, we show that an irreducible highest weight $U_q(\mf{g}_{m|n})$-module  $L_q(\mf{g}_{m|n},\Lambda_{m|n}(\lambda,\ell))$  with highest weight $\Lambda_{m|n}(\lambda,\ell)$ is an irreducible module in $\mc{O}^{int}_q(m|n)$ and it has a unique crystal base up to scalar multiplication. Moreover, we  realize the crystal of $L_q(\mf{g}_{m|n},\Lambda_{m|n}(\lambda,\ell))$ in terms of a new combinatorial object so-called ortho-symplectic tableaux of shape $(\lambda,\ell)$ (Theorem \ref{Existence of crystal base}). 
Since every irreducible module in $\mc{O}^{int}_q(m|n)$ is isomorphic to $L_q(\mf{g}_{m|n},\Lambda_{m|n}(\lambda,\ell))$ for some $(\lambda,\ell)\in \cP(\g)_{m|n}$, it follows that $\mc{O}^{int}_q(m|n)$ is equivalent to $\mc{O}^{int}(m|n)$.
Therefore, we obtain the following, which is the main result in this paper.

\begin{thm*}
For $\mf{g}=\mf{b}, \mf{b}^\bullet, \mf{c}$, $\mc{O}^{int}_q(m|n)$ is a semisimple tensor category, which is equivalent to $\mc{O}^{int}(m|n)$, and each irreducible module in  $\mc{O}^{int}_q(m|n)$ has a unique crystal base.
\end{thm*}

We remark that an ortho-symplectic tableau is defined over an arbitrary linearly ordered $\Z_2$-graded set $\mc{A}$ and its main advantage is compatibility with super duality functor $\mc{F}$ and truncation functor $\mf{tr}_n$. More precisely, the set of ortho-symplectic tableaux of shape $(\lambda,\ell)$ gives both irreducible characters in $\mc{O}^{int}_q(m+n)$ and $\mc{O}^{int}_q(m|n)$ with suitable choices of $\mc{A}$. The Schur positivity of the character of ortho-symplectic tableaux of shape $(\lambda,\ell)$ plays a crucial role in proving this compatibility. Also, when $\mc{A}$ is a finite set with even elements, we have an explicit bijection between ortho-symplectic tableaux and  Kashiwara-Nakashima tableaux \cite{KashNaka} of type $B$ and $C$.  

The above theorem is also true for ${\mf g}={\mf d}$. The major difference in the proof is that  more technical difficulty enters when we describe the crystal structure of an highest weight module and hence define the notion of ortho-symplectic tableaux of type $D$. We discuss the details in a subsequent paper \cite{K14}.

\subsection{}
The paper is organized as follows. In Section \ref{Lie superalgebra}, we recall the definition of ortho-symplectic Lie superalgebras $\mf{g}_{m|n}$ based on \cite{CW}. In Section \ref{Super duality and semisimple category}, we briefly review super duality and present a simple characterization of $\mc{O}^{int}(m|n)$. In Section \ref{Semisimple category for quantum ortho-symplectic}, we define $\mc{O}^{int}_q(m|n)$ and show that it is a semisimple tensor category.
In Section \ref{Crystal base of Fock spaces}, we review the notion of a crystal base for a quantum superalgebra \cite{BKK} and prove the existence of a crystal base of a $q$-deformed Fock space $\mathscr{V}_q$. Then  we introduce our main combinatorial object ${\bf T}_{m|n}(\lambda,\ell)$, the set of ortho-symplectic tableaux of shape $(\lambda,\ell)$,  for $\mf{g}=\mf{b}, \mf{b}^\bullet, \mf{c}$ in Section \ref{ortho-symplectic Tableaux of Type B and C}, and show that its character gives an irreducible character in $\mc{O}^{int}(m|n)$  in Section \ref{crystal structure for ortho-symplectic tableaux m+n}. Finally, in Section \ref{Crystal base for m|n}, we show that $L_q(\mf{g}_{m|n},\Lambda_{m|n}(\lambda,\ell))$ is an irreducible $U_q(\g_{m|n})$-module in $\mc{O}^{int}_q(m|n)$ and has a unique crystal base for $\mf{g}=\mf{b}, \mf{b}^\bullet, \mf{c}$ and $(\lambda,\ell)\in \cP(\mf{g})_{m|n}$, whose crystal is realized by ${\bf T}_{m|n}(\lambda,\ell)$.

\vskip 3mm
{\bf Acknowledgement} The author would like to thank S.-J. Cheng and W. Wang for valuable discussion on super duality and allowing him the draft of their book \cite{CW}, and the  referee for careful reading and pointing out some mistakes in the first version of the paper. This work was first announced at the 2nd Mini Symposium on Representation Theory in Jeju, Korea, December 2012. He thanks S.-J. Kang for his kind invitation and interest in this work.  

\section{Ortho-symplectic Lie superalgebra $\mf{g}_{m|n}$}\label{Lie superalgebra}
In this section, let us briefly recall some necessary background on Lie superalgebras (see \cite{CW,Kac77} for more details). Our exposition is based on \cite{CW} with a little modification. We assume that the base field is $\C$.

\subsection{General linear Lie superalgebras}
Throughout this paper, we fix a positive integer $m$ and let
\begin{equation*}
\begin{split}
&\td{\I}_m=\{\,\overline{k},-\overline{k}\,|\,1\leq k\leq m\,\}\cup \tfrac{1}{2}\Z,\\
&\I_m=\{\,\overline{k},-\overline{k} \,|\,1\leq k\leq m\,\}\cup  \Z,\\
&\ov{\I}_m=\{\,\overline{k},-\overline{k} \,|\,1\leq k\leq m\,\}\cup \left(\tfrac{1}{2}+\Z\right)\cup\{0\},\\
\end{split}
\end{equation*}
where $\td{\I}_m$ is a linearly ordered $\Z_2$-graded set with \begin{equation*}\label{the alphabet I}
\begin{split}
\cdots<-\tfrac{3}{2} <-1 <-\tfrac{1}{2}  <-\ov{1} <\cdots<-\ov{m} <&0<\ov{m}<\cdots<\ov{1}<\tfrac{1}{2}<1 < \tfrac{3}{2}<\ldots,\\
\left(\td{\I}_m\right)_0\supset \{\,\overline{k},-\overline{k}\,|\,1\leq k\leq m\,\}\cup \Z^\times, \ &\ \ \left(\td{\I}_m\right)_1=\tfrac{1}{2}+\Z,
\end{split}
\end{equation*}
(the parity of $0$ will be specified later) and the linear orderings and $\Z_2$-gradings on the other sets are induced from those on $\td{\I}_m$. For $a\in \td{\I}_m$, $|a|$ denotes the parity of $a$.
We put $\td{\I}_m^+=\{\,a\in \td{\I}_m^+\,|\,a>0\,\}$ and $\I^\times=\{\,a\in\I\,|\, a\neq 0\,\}$ for $\I\subset \td{\I}_m$.

For $\I\subset \td{\I}_m$, we denote by $V_\I$ the superspace with basis $\{\,v_a\,|\,a\in \I\,\}$, where the $\Z_2$-grading is induced from $\td{\I}_m$. Let $\gl(V_{\I})$ be the general linear Lie superalgebra of linear endomorphisms on $V_\I$ vanishing on $v_a$'s except for finitely many $a$'s. We identify $\gl(V_\I)$ with the space of matrices $(a_{ij})_{i,j\in \I}$ spanned  by
the elementary matrices $E_{i,j}$. We assume that $V_\I$ is a subspace of $V_{\td{\I}_m}$ and $\gl(V_\I)\subset \gl(V_{\td{\I}_m})$. Let
$\widehat{\gl}(V_{\I})$ be the central extension
of $\gl(V_{\I})$ by a one-dimensional center $\C K$ with respect to the
$2$-cocycle $\alpha(A,B)={\rm Str}([{J},A]B)$,
where ${\rm Str}$ is the supertrace with ${\rm Str}(a_{ij})=\sum_{i\in \I}(-1)^{|i|}a_{ii}$  and ${J}=\sum_{i\le 0}E_{i,i}$.

For $n\in \Z_{> 0}\cup\{\infty\}$, we put
\begin{equation*}
\begin{split}
\J_{m+n}&=\{\,a\in \I_{m}\,|\, \ov{m}\leq  a\leq n\,\}, \\   \J_{m|n}&=\{\,a\in \ov{\I}_{m}\,|\,\ov{m}\leq a\leq n-\tfrac{1}{2}\,\}.
\end{split}
\end{equation*}
We assume that $\J_{m+0}=\J_{m|0}=\{\,\ov{m},\ldots,\ov{1}\,\}$.

\subsection{Ortho-symplectic Lie superalgebras}\label{def of osp}
Suppose that the parity of $0\in\td{\I}_m$ is $1$. Define a skew-supersymmetric bilinear form $(\,\cdot\,|\,\cdot\,)$ on $V_{\td{\I}_m}$ by
\begin{equation}\label{skewsupersymmetric}
\begin{split}
&(v_{\pm a}|v_{\pm b})=0, \ \ (v_a|v_{-b})=-(-1)^{|a||b|}(v_{-b}|v_a)=\delta_{a b}, \\
&(v_{0}|v_{0})=1, \ \ \ \ \ (v_{0}|v_{\pm a})=0,
\end{split}
\end{equation}
for $a,b \in \td{\I}_m^+$.
For $\I\subset \td{\I}_m$, let $\frak{spo}(V_{\I})$ be the subalgebra  of $\gl(V_{\I})$ preserving the skew-supersymmetric bilinear form  on $V_\I$ induced from \eqref{skewsupersymmetric}. Then we define  $\mf{b}^\bullet_{m+\infty}$,  $\mf{b}^\bullet_{m|\infty}$, $\mf{c}_{m+\infty}$ and $\mf{c}_{m|\infty}$ to be the central extensions of $\frak{spo}(V_{\I})$ induced from  $\widehat{\gl}(V_{\I})$ when $\I$ is $\I_m$, $\ov{\I}_m$, $\I^\times_m$ and $\ov{\I}_m^\times$, respectively.

Next, suppose that the parity of $0\in\td{\I}_m$ is $0$. Define a supersymmetric bilinear form $(\,\cdot\,|\,\cdot\,)$ on $V_{\td{\I}_m}$ by
\begin{equation}\label{supersymmetric}
\begin{split}
&(v_{\pm a}|v_{\pm b})=0, \ \ (v_a|v_{-b})=(-1)^{|a||b|}(v_{-b}|v_a)=\delta_{a b}, \\
&(v_{0}|v_{0})=1, \ \ \ \ \ (v_{0}|v_{\pm a})=0,
\end{split}
\end{equation}
for $a,b \in \td{\I}_m^+$. For $\I\subset \td{\I}_m$, let $\frak{osp}(V_{\I})$ be the subalgebra  of $\gl(V_{\I})$ preserving the supersymmetric bilinear form  on $V_\I$ induced from \eqref{supersymmetric}. Then we define  $\mf{b}_{m+\infty}$,  $\mf{b}_{m|\infty}$, $\mf{d}_{m+\infty}$ and $\mf{d}_{m|\infty}$ to be the central extensions of $\frak{osp}(V_{\I})$ induced from  $\widehat{\gl}(V_{\I})$ when $\I$ is $\I_m$, $\ov{\I}_m$, $\I^\times_m$ and $\ov{\I}_m^\times$, respectively.

From now on, we assume that $\mf{g}$ is a symbol, which denotes one of $\mf{b}$, $\mf{b}^\bullet$, $\mf{c}$ and $\mf{d}$. Let $U(\mf{g}_{m+\infty})$ and $U(\mf{g}_{m|\infty})$ be the enveloping superalgebras associated to $\mf{g}_{m+\infty}$ and $\mf{g}_{m|\infty}$, respectively.

Let  $\h_{m+\infty}$ (resp. $\h_{m|\infty}$) be the Cartan subalgebra of $\mf{g}_{m+\infty}$ (resp. $\mf{g}_{m|\infty}$) spanned by  $K$ and $E_{a}:=E_{a,a}-E_{-a, -a}$ for $a\in \J_{m+\infty}$ (resp. $\J_{m|\infty}$), and let $\h_{m+\infty}^*$ (resp. $\h_{m|\infty}^*$) be the restricted dual of $\h_{m+\infty}$ (resp. $\h_{m|\infty}$) spanned by $\Lambda_{\ov{m}}$ and $\delta_a$ for $a\in \J_{m+\infty}$ (resp. $\J_{m|\infty}$), where
$\langle E_b,\delta_a\rangle=\delta_{ab}$, $\langle K,\delta_a\rangle=0$, $\langle E_a,\Lambda_{\ov{m}}\rangle=0$  for $a, b$ and $\langle K,\Lambda_{\ov{m}}\rangle=r$ with $r=1$ for $\mf{g}= \mf{c}$ and $r=\tfrac{1}{2}$ otherwise. Here $\langle\,\cdot\,,\,\cdot\,\rangle$ denotes the natural pairing on $\mf{h}_{m+\infty}\times \mf{h}^*_{m+\infty}$ or $\mf{h}_{m|\infty}\times \mf{h}^*_{m|\infty}$.

Let $I_{m+\infty}=\{\,\ov{m},\ldots,\ov{1},0\}\cup \Z_{> 0}$. Then  the set of simple roots $\Pi_{m+\infty}=\{\,\alpha_i\,|\,i\in I_{m+\infty}\,\}$,  the set of simple coroots $\Pi^\vee_{m+\infty}=\{\,\alpha_i^\vee\,|\,i\in I_{m+\infty}\,\}$ and the Dynkin diagram associated to  the Cartan matrix $(\langle \alpha_i^\vee,\alpha_j\rangle)_{i,j\in I_{m+\infty}}$  of $\mf{g}_{m+\infty}$ are listed below (the simple roots are with respect to a Borel subalgebra  spanned by the upper triangular matrices):\vskip 3mm

$\bullet$  $\mf{b}^\bullet_{m+\infty}$

\begin{equation*}
\alpha_i=
\begin{cases}
-\delta_{\ov{m}}, & \text{if $i=\ov{m}$},\\
\delta_{\ov{k+1}}-\delta_{\ov{k}}, & \text{if $i=\ov{k}\ (\neq \ov{m})$},\\
\delta_{\ov{1}}-\delta_{1}, & \text{if $i=0$},\\
\delta_i-\delta_{i+1}, & \text{if $i\in \Z_{>0}$},
\end{cases} \ \
\alpha_i^\vee=
\begin{cases}
-2E_{\ov{m}}+2K, & \text{if $i=\ov{m}$},\\
E_{\ov{k+1}}-E_{\ov{k}}, & \text{if $i=\ov{k}\ (\neq \ov{m})$},\\
E_{\ov{1}}-E_{1}, & \text{if $i=0$},\\
E_i-E_{i+1}, & \text{if $i\in \Z_{>0}$}.
\end{cases}
\end{equation*}
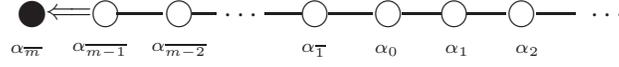
\begin{center}
\hskip -3cm \setlength{\unitlength}{0.16in}
\begin{picture}(24,4.3)
\put(5.6,2){\circle*{0.9}}
\put(8,2){\makebox(0,0)[c]{$\bigcirc$}}
\put(10.4,2){\makebox(0,0)[c]{$\bigcirc$}}
\put(14.85,2){\makebox(0,0)[c]{$\bigcirc$}}
\put(17.25,2){\makebox(0,0)[c]{$\bigcirc$}}
\put(19.4,2){\makebox(0,0)[c]{$\bigcirc$}}
\put(21.6,2){\makebox(0,0)[c]{$\bigcirc$}}
\put(8.35,2){\line(1,0){1.5}}
\put(10.82,2){\line(1,0){0.8}}
\put(13.2,2){\line(1,0){1.2}}
\put(15.28,2){\line(1,0){1.45}}
\put(17.7,2){\line(1,0){1.25}}
\put(19.81,2){\line(1,0){1.28}}
\put(22.1,2){\line(1,0){1.28}}
%
\put(6.8,2){\makebox(0,0)[c]{$\Longleftarrow$}}
\put(12.5,1.95){\makebox(0,0)[c]{$\cdots$}}
\put(24.5,1.95){\makebox(0,0)[c]{$\cdots$}}
\put(5.4,0.8){\makebox(0,0)[c]{\tiny $\alpha_{\ov{m}}$}}
\put(7.8,0.8){\makebox(0,0)[c]{\tiny $\alpha_{\ov{m-1}}$}}
\put(10.4,0.8){\makebox(0,0)[c]{\tiny $\alpha_{\ov{m-2}}$}}
\put(14.8,0.8){\makebox(0,0)[c]{\tiny $\alpha_{\ov{1}}$}}
\put(17.2,0.8){\makebox(0,0)[c]{\tiny $\alpha_{0}$}}
\put(19.5,0.8){\makebox(0,0)[c]{\tiny $\alpha_{1}$}}
\put(21.8,0.8){\makebox(0,0)[c]{\tiny $\alpha_{2}$}}
\end{picture}
\end{center}
\ \ \ \  \ \   (\ $
\begin{picture}(2,2)\setlength{\unitlength}{0.14in}
\put(0.3,0.3){\circle*{1}}
\end{picture}
$ \ \  denotes a non-isotropic odd simple root.)

$\bullet$ $\mf{c}_{m+\infty}$,

\begin{equation*}
\alpha_i=
\begin{cases}
-2\delta_{\ov{m}}, & \text{if $i=\ov{m}$},\\
\delta_{\ov{k+1}}-\delta_{\ov{k}}, & \text{if $i=\ov{k}\ (\neq \ov{m})$},\\
\delta_{\ov{1}}-\delta_{1}, & \text{if $i=0$},\\
\delta_i-\delta_{i+1}, & \text{if $i\in \Z_{>0}$},
\end{cases}\ \
\alpha_i^\vee=
\begin{cases}
-E_{\ov{m}}+K, & \text{if $i=\ov{m}$},\\
E_{\ov{k+1}}-E_{\ov{k}}, & \text{if $i=\ov{k}\ (\neq \ov{m})$},\\
E_{\ov{1}}-E_{1}, & \text{if $i=0$},\\
E_i-E_{i+1}, & \text{if $i\in \Z_{>0}$}.
\end{cases}
\end{equation*}
\begin{center}
\hskip -3cm \setlength{\unitlength}{0.16in}
\begin{picture}(24,4)
\put(5.6,2){\makebox(0,0)[c]{$\bigcirc$}}
\put(8,2){\makebox(0,0)[c]{$\bigcirc$}}
\put(10.4,2){\makebox(0,0)[c]{$\bigcirc$}}
\put(14.85,2){\makebox(0,0)[c]{$\bigcirc$}}
\put(17.25,2){\makebox(0,0)[c]{$\bigcirc$}}
\put(19.4,2){\makebox(0,0)[c]{$\bigcirc$}}
\put(21.6,2){\makebox(0,0)[c]{$\bigcirc$}}
\put(8.35,2){\line(1,0){1.5}}
\put(10.82,2){\line(1,0){0.8}}
\put(13.2,2){\line(1,0){1.2}}
\put(15.28,2){\line(1,0){1.45}}
\put(17.7,2){\line(1,0){1.25}}
\put(19.81,2){\line(1,0){1.28}}
\put(22.1,2){\line(1,0){1.28}}
%
\put(6.8,2){\makebox(0,0)[c]{$\Longrightarrow$}}
\put(12.5,1.95){\makebox(0,0)[c]{$\cdots$}}
\put(24.5,1.95){\makebox(0,0)[c]{$\cdots$}}
\put(5.4,0.8){\makebox(0,0)[c]{\tiny $\alpha_{\ov{m}}$}}
\put(7.8,0.8){\makebox(0,0)[c]{\tiny $\alpha_{\ov{m-1}}$}}
\put(10.4,0.8){\makebox(0,0)[c]{\tiny $\alpha_{\ov{m-2}}$}}
\put(14.8,0.8){\makebox(0,0)[c]{\tiny $\alpha_{\ov{1}}$}}
\put(17.2,0.8){\makebox(0,0)[c]{\tiny $\alpha_{0}$}}
\put(19.5,0.8){\makebox(0,0)[c]{\tiny $\alpha_{1}$}}
\put(21.8,0.8){\makebox(0,0)[c]{\tiny $\alpha_{2}$}}
\end{picture}
\end{center}

$\bullet$  $\mf{b}_{m+\infty}$

\begin{equation*}
\alpha_i=
\begin{cases}
-\delta_{\ov{m}}, & \text{if $i=\ov{m}$},\\
\delta_{\ov{k+1}}-\delta_{\ov{k}}, & \text{if $i=\ov{k}\ (\neq \ov{m})$},\\
\delta_{\ov{1}}-\delta_{1}, & \text{if $i=0$},\\
\delta_i-\delta_{i+1}, & \text{if $i\in \Z_{>0}$},
\end{cases} \ \
\alpha_i^\vee=
\begin{cases}
-2E_{\ov{m}}+2K, & \text{if $i=\ov{m}$},\\
E_{\ov{k+1}}-E_{\ov{k}}, & \text{if $i=\ov{k}\ (\neq \ov{m})$},\\
E_{\ov{1}}-E_{1}, & \text{if $i=0$},\\
E_i-E_{i+1}, & \text{if $i\in \Z_{>0}$}.
\end{cases}
\end{equation*}
\begin{center}
\hskip -3cm \setlength{\unitlength}{0.16in}
\begin{picture}(24,4)
\put(5.6,2){\makebox(0,0)[c]{$\bigcirc$}}
\put(8,2){\makebox(0,0)[c]{$\bigcirc$}}
\put(10.4,2){\makebox(0,0)[c]{$\bigcirc$}}
\put(14.85,2){\makebox(0,0)[c]{$\bigcirc$}}
\put(17.25,2){\makebox(0,0)[c]{$\bigcirc$}}
\put(19.4,2){\makebox(0,0)[c]{$\bigcirc$}}
\put(21.6,2){\makebox(0,0)[c]{$\bigcirc$}}
\put(8.35,2){\line(1,0){1.5}}
\put(10.82,2){\line(1,0){0.8}}
\put(13.2,2){\line(1,0){1.2}}
\put(15.28,2){\line(1,0){1.45}}
\put(17.7,2){\line(1,0){1.25}}
\put(19.81,2){\line(1,0){1.28}}
\put(22.1,2){\line(1,0){1.28}}
%
\put(6.8,2){\makebox(0,0)[c]{$\Longleftarrow$}}
\put(12.5,1.95){\makebox(0,0)[c]{$\cdots$}}
\put(24.5,1.95){\makebox(0,0)[c]{$\cdots$}}
\put(5.4,0.8){\makebox(0,0)[c]{\tiny $\alpha_{\ov{m}}$}}
\put(7.8,0.8){\makebox(0,0)[c]{\tiny $\alpha_{\ov{m-1}}$}}
\put(10.4,0.8){\makebox(0,0)[c]{\tiny $\alpha_{\ov{m-2}}$}}
\put(14.8,0.8){\makebox(0,0)[c]{\tiny $\alpha_{\ov{1}}$}}
\put(17.2,0.8){\makebox(0,0)[c]{\tiny $\alpha_{0}$}}
\put(19.5,0.8){\makebox(0,0)[c]{\tiny $\alpha_{1}$}}
\put(21.8,0.8){\makebox(0,0)[c]{\tiny $\alpha_{2}$}}
\end{picture}
\end{center}

$\bullet$  $\mf{d}_{m+\infty}$

\begin{equation*}
\alpha_i=
\begin{cases}
-\delta_{\ov{m}}-\delta_{\ov{m-1}}, & \text{if $i=\ov{m}$},\\
\delta_{\ov{k+1}}-\delta_{\ov{k}}, & \text{if $i=\ov{k}\ (\neq \ov{m})$},\\
\delta_{\ov{1}}-\delta_{1}, & \text{if $i=0$},\\
\delta_i-\delta_{i+1}, & \text{if $i\in \Z_{>0}$},
\end{cases} \ \,
\alpha_i^\vee=
\begin{cases}
-E_{\ov{m}}-E_{\ov{m-1}}+2K, & \text{if $i=\ov{m}$},\\
E_{\ov{k+1}}-E_{\ov{k}}, & \text{if $i=\ov{k}\ (\neq \ov{m})$},\\
E_{\ov{1}}-E_{1}, & \text{if $i=0$},\\
E_i-E_{i+1}, & \text{if $i\in \Z_{>0}$}.
\end{cases}
\end{equation*}
\begin{center}
\hskip -3cm \setlength{\unitlength}{0.16in} \medskip
\begin{picture}(24,5.8)
\put(6,0){\makebox(0,0)[c]{$\bigcirc$}}
\put(6,4){\makebox(0,0)[c]{$\bigcirc$}}
\put(8,2){\makebox(0,0)[c]{$\bigcirc$}}
\put(10.4,2){\makebox(0,0)[c]{$\bigcirc$}}
\put(14.85,2){\makebox(0,0)[c]{$\bigcirc$}}
\put(17.25,2){\makebox(0,0)[c]{$\bigcirc$}}
\put(19.4,2){\makebox(0,0)[c]{$\bigcirc$}}
\put(21.5,2){\makebox(0,0)[c]{$\bigcirc$}}
\put(6.35,0.3){\line(1,1){1.35}} \put(6.35,3.7){\line(1,-1){1.35}}
\put(8.4,2){\line(1,0){1.55}} \put(10.82,2){\line(1,0){0.8}}
\put(13.2,2){\line(1,0){1.2}} \put(15.28,2){\line(1,0){1.45}}
\put(17.7,2){\line(1,0){1.25}} \put(19.8,2){\line(1,0){1.25}}
\put(21.95,2){\line(1,0){1.4}} 
\put(12.5,1.95){\makebox(0,0)[c]{$\cdots$}}
\put(24.5,1.95){\makebox(0,0)[c]{$\cdots$}}
\put(6,5){\makebox(0,0)[c]{\tiny $\alpha_{\ov{m}}$}}
\put(6,-1.2){\makebox(0,0)[c]{\tiny $\alpha_{\ov{m-1}}$}}
\put(8.2,1){\makebox(0,0)[c]{\tiny $\alpha_{\ov{m-2}}$}}
\put(10.4,1){\makebox(0,0)[c]{\tiny $\alpha_{\ov{m-3}}$}}
\put(14.8,1){\makebox(0,0)[c]{\tiny $\alpha_{\ov{1}}$}}
\put(17.15,1){\makebox(0,0)[c]{\tiny $\alpha_0$}}
\put(19.5,1){\makebox(0,0)[c]{\tiny $\alpha_{1}$}}
\put(21.5,1){\makebox(0,0)[c]{\tiny $\alpha_{2}$}}
\end{picture}
\end{center}\vskip 5mm

Let $I_{m|\infty}=\{\,\ov{m},\ldots,\ov{1},0\}\cup \left(\hf+\Z_{\geq 0}\right)$. Then  the set of simple roots $\Pi_{m|\infty}=\{\,\beta_i\,|\,i\in I_{m|\infty}\,\}$,  the set of simple coroots $\Pi^\vee_{m|\infty}=\{\,\beta_i^\vee\,|\,i\in I_{m|\infty}\,\}$ and the Dynkin diagram associated to  the Cartan matrix $(\langle \beta_i^\vee,\beta_j\rangle)_{i,j\in I_{m|\infty}}$  of $\mf{g}_{m|\infty}$ are listed below (the simple roots are with respect to a Borel subalgebra  spanned by the upper triangular matrices):\vskip 3mm
\vskip 3mm

$\bullet$ $\mf{b}^\bullet_{m|\infty}$

\begin{equation*}
\beta_i=
\begin{cases}
-\delta_{\ov{m}}, & \text{if $i=\ov{m}$},\\
\delta_{\ov{k+1}}-\delta_{\ov{k}}, & \text{if $i=\ov{k}\ (\neq \ov{m})$},\\
\delta_{\ov{1}}-\delta_{\hf}, & \text{if $i=0$},\\
\delta_i-\delta_{i+1}, & \text{if $i\in \hf+\Z_{\geq 0}$},
\end{cases}\ \
\beta_i^\vee=
\begin{cases}
-2E_{\ov{m}}+2K, & \text{if $i=\ov{m}$},\\
E_{\ov{k+1}}-E_{\ov{k}}, & \text{if $i=\ov{k}\ (\neq \ov{m})$},\\
E_{\ov{1}}+E_{\hf}, & \text{if $i=0$},\\
E_i-E_{i+1}, & \text{if $i\in \hf+\Z_{\geq 0}$}.
\end{cases}
\end{equation*}
\begin{center}
\hskip -3cm \setlength{\unitlength}{0.16in}
\begin{picture}(24,4)
\put(5.6,2){\circle*{0.9}}
\put(8,2){\makebox(0,0)[c]{$\bigcirc$}}
\put(10.4,2){\makebox(0,0)[c]{$\bigcirc$}}
\put(14.85,2){\makebox(0,0)[c]{$\bigcirc$}}
\put(17.25,2){\makebox(0,0)[c]{$\bigotimes$}}
\put(19.4,2){\makebox(0,0)[c]{$\bigcirc$}}
\put(21.6,2){\makebox(0,0)[c]{$\bigcirc$}}
\put(8.35,2){\line(1,0){1.5}}
\put(10.82,2){\line(1,0){0.8}}
\put(13.2,2){\line(1,0){1.2}}
\put(15.28,2){\line(1,0){1.45}}
\put(17.7,2){\line(1,0){1.25}}
\put(19.81,2){\line(1,0){1.28}}
\put(22.1,2){\line(1,0){1.28}}
%
\put(6.8,2){\makebox(0,0)[c]{$\Longleftarrow$}}
\put(12.5,1.95){\makebox(0,0)[c]{$\cdots$}}
\put(24.5,1.95){\makebox(0,0)[c]{$\cdots$}}
\put(5.4,0.8){\makebox(0,0)[c]{\tiny $\beta_{\ov{m}}$}}
\put(7.8,0.8){\makebox(0,0)[c]{\tiny $\beta_{\ov{m-1}}$}}
\put(10.4,0.8){\makebox(0,0)[c]{\tiny $\beta_{\ov{m-2}}$}}
\put(14.8,0.8){\makebox(0,0)[c]{\tiny $\beta_{\ov{1}}$}}
\put(17.2,0.8){\makebox(0,0)[c]{\tiny $\beta_{0}$}}
\put(19.5,0.8){\makebox(0,0)[c]{\tiny $\beta_{\frac{1}{2}}$}}
\put(21.8,0.8){\makebox(0,0)[c]{\tiny $\beta_{\frac{3}{2}}$}}
\end{picture}
\end{center}
\ \ \ \  \ \ ( $\bigotimes$ denotes an isotropic odd simple root.)

$\bullet$ $\mf{c}_{m|\infty}$

\begin{equation*}
\beta_i=
\begin{cases}
-2\delta_{\ov{m}}, & \text{if $i=\ov{m}$},\\
\delta_{\ov{k+1}}-\delta_{\ov{k}}, & \text{if $i=\ov{k}\ (\neq \ov{m})$},\\
\delta_{\ov{1}}-\delta_{\hf}, & \text{if $i=0$},\\
\delta_i-\delta_{i+1}, & \text{if $i\in \hf+\Z_{\geq 0}$},
\end{cases}\ \
\beta_i^\vee=
\begin{cases}
-E_{\ov{m}}+K, & \text{if $i=\ov{m}$},\\
E_{\ov{k+1}}-E_{\ov{k}}, & \text{if $i=\ov{k}\ (\neq \ov{m})$},\\
E_{\ov{1}}+E_{\hf}, & \text{if $i=0$},\\
E_i-E_{i+1}, & \text{if $i\in \hf+\Z_{\geq 0}$}.
\end{cases}
\end{equation*}

\begin{center}
\hskip -3cm \setlength{\unitlength}{0.16in}
\begin{picture}(24,4)
\put(5.6,2){\makebox(0,0)[c]{$\bigcirc$}}
\put(8,2){\makebox(0,0)[c]{$\bigcirc$}}
\put(10.4,2){\makebox(0,0)[c]{$\bigcirc$}}
\put(14.85,2){\makebox(0,0)[c]{$\bigcirc$}}
\put(17.25,2){\makebox(0,0)[c]{$\bigotimes$}}
\put(19.4,2){\makebox(0,0)[c]{$\bigcirc$}}
\put(21.6,2){\makebox(0,0)[c]{$\bigcirc$}}
\put(8.35,2){\line(1,0){1.5}}
\put(10.82,2){\line(1,0){0.8}}
\put(13.2,2){\line(1,0){1.2}}
\put(15.28,2){\line(1,0){1.45}}
\put(17.7,2){\line(1,0){1.25}}
\put(19.81,2){\line(1,0){1.28}}
\put(22.1,2){\line(1,0){1.28}}
%
\put(6.8,2){\makebox(0,0)[c]{$\Longrightarrow$}}
\put(12.5,1.95){\makebox(0,0)[c]{$\cdots$}}
\put(24.5,1.95){\makebox(0,0)[c]{$\cdots$}}
\put(5.4,0.8){\makebox(0,0)[c]{\tiny $\beta_{\ov{m}}$}}
\put(7.8,0.8){\makebox(0,0)[c]{\tiny $\beta_{\ov{m-1}}$}}
\put(10.4,0.8){\makebox(0,0)[c]{\tiny $\beta_{\ov{m-2}}$}}
\put(14.8,0.8){\makebox(0,0)[c]{\tiny $\beta_{\ov{1}}$}}
\put(17.2,0.8){\makebox(0,0)[c]{\tiny $\beta_{0}$}}
\put(19.5,0.8){\makebox(0,0)[c]{\tiny $\beta_{\frac{1}{2}}$}}
\put(21.8,0.8){\makebox(0,0)[c]{\tiny $\beta_{\frac{3}{2}}$}}
\end{picture}
\end{center}

$\bullet$ $\mf{b}_{m|\infty}$

\begin{equation*}
\beta_i=
\begin{cases}
-\delta_{\ov{m}}, & \text{if $i=\ov{m}$},\\
\delta_{\ov{k+1}}-\delta_{\ov{k}}, & \text{if $i=\ov{k}\ (\neq \ov{m})$},\\
\delta_{\ov{1}}-\delta_{\hf}, & \text{if $i=0$},\\
\delta_i-\delta_{i+1}, & \text{if $i\in \hf+\Z_{\geq 0}$},
\end{cases}\ \
\beta_i^\vee=
\begin{cases}
-2E_{\ov{m}}+2K, & \text{if $i=\ov{m}$},\\
E_{\ov{k+1}}-E_{\ov{k}}, & \text{if $i=\ov{k}\ (\neq \ov{m})$},\\
E_{\ov{1}}+E_{\hf}, & \text{if $i=0$},\\
E_i-E_{i+1}, & \text{if $i\in \hf+\Z_{\geq 0}$}.
\end{cases}
\end{equation*}

\begin{center}
\hskip -3cm \setlength{\unitlength}{0.16in}
\begin{picture}(24,4)
\put(5.6,2){\makebox(0,0)[c]{$\bigcirc$}}
\put(8,2){\makebox(0,0)[c]{$\bigcirc$}}
\put(10.4,2){\makebox(0,0)[c]{$\bigcirc$}}
\put(14.85,2){\makebox(0,0)[c]{$\bigcirc$}}
\put(17.25,2){\makebox(0,0)[c]{$\bigotimes$}}
\put(19.4,2){\makebox(0,0)[c]{$\bigcirc$}}
\put(21.6,2){\makebox(0,0)[c]{$\bigcirc$}}
\put(8.35,2){\line(1,0){1.5}}
\put(10.82,2){\line(1,0){0.8}}
\put(13.2,2){\line(1,0){1.2}}
\put(15.28,2){\line(1,0){1.45}}
\put(17.7,2){\line(1,0){1.25}}
\put(19.81,2){\line(1,0){1.28}}
\put(22.1,2){\line(1,0){1.28}}
%
\put(6.8,2){\makebox(0,0)[c]{$\Longleftarrow$}}
\put(12.5,1.95){\makebox(0,0)[c]{$\cdots$}}
\put(24.5,1.95){\makebox(0,0)[c]{$\cdots$}}
\put(5.4,0.8){\makebox(0,0)[c]{\tiny $\beta_{\ov{m}}$}}
\put(7.8,0.8){\makebox(0,0)[c]{\tiny $\beta_{\ov{m-1}}$}}
\put(10.4,0.8){\makebox(0,0)[c]{\tiny $\beta_{\ov{m-2}}$}}
\put(14.8,0.8){\makebox(0,0)[c]{\tiny $\beta_{\ov{1}}$}}
\put(17.2,0.8){\makebox(0,0)[c]{\tiny $\beta_{0}$}}
\put(19.5,0.8){\makebox(0,0)[c]{\tiny $\beta_{\frac{1}{2}}$}}
\put(21.8,0.8){\makebox(0,0)[c]{\tiny $\beta_{\frac{3}{2}}$}}
\end{picture}
\end{center}

$\bullet$ $\mf{d}_{m|\infty}$

\begin{equation*}
\beta_i=
\begin{cases}
-\delta_{\ov{m}}-\delta_{\ov{m-1}}, & \text{if $i=\ov{m}$},\\
\delta_{\ov{k+1}}-\delta_{\ov{k}}, & \text{if $i=\ov{k}\ (\neq \ov{m})$},\\
\delta_{\ov{1}}-\delta_{\hf}, & \text{if $i=0$},\\
\delta_i-\delta_{i+1}, & \text{if $i\in \hf+\Z_{\geq 0}$},
\end{cases}\ \,
\beta_i^\vee=
\begin{cases}
-E_{\ov{m}}-E_{\ov{m-1}}+2K, & \text{if $i=\ov{m}$},\\
E_{\ov{k+1}}-E_{\ov{k}}, & \text{if $i=\ov{k}\ (\neq \ov{m})$},\\
E_{\ov{1}}+E_{\hf}, & \text{if $i=0$},\\
E_i-E_{i+1}, & \text{if $i\in \hf+\Z_{\geq 0}$}.
\end{cases}
\end{equation*}

\begin{center}
\hskip -3cm \setlength{\unitlength}{0.16in} \medskip
\begin{picture}(24,5.8)
\put(6,0){\makebox(0,0)[c]{$\bigcirc$}}
\put(6,4){\makebox(0,0)[c]{$\bigcirc$}}
\put(8,2){\makebox(0,0)[c]{$\bigcirc$}}
\put(10.4,2){\makebox(0,0)[c]{$\bigcirc$}}
\put(14.85,2){\makebox(0,0)[c]{$\bigcirc$}}
\put(17.25,2){\makebox(0,0)[c]{$\bigotimes$}}
\put(19.4,2){\makebox(0,0)[c]{$\bigcirc$}}
\put(21.5,2){\makebox(0,0)[c]{$\bigcirc$}}
\put(6.35,0.3){\line(1,1){1.35}} \put(6.35,3.7){\line(1,-1){1.35}}
\put(8.4,2){\line(1,0){1.55}} \put(10.82,2){\line(1,0){0.8}}
\put(13.2,2){\line(1,0){1.2}} \put(15.28,2){\line(1,0){1.45}}
\put(17.7,2){\line(1,0){1.25}} \put(19.8,2){\line(1,0){1.25}}
\put(21.95,2){\line(1,0){1.4}} 
\put(12.5,1.95){\makebox(0,0)[c]{$\cdots$}}
\put(24.5,1.95){\makebox(0,0)[c]{$\cdots$}}
\put(6,5){\makebox(0,0)[c]{\tiny $\beta_{\ov{m}}$}}
\put(6,-1.2){\makebox(0,0)[c]{\tiny $\beta_{\ov{m-1}}$}}
\put(8.2,1){\makebox(0,0)[c]{\tiny $\beta_{\ov{m-2}}$}}
\put(10.4,1){\makebox(0,0)[c]{\tiny $\beta_{\ov{m-3}}$}}
\put(14.9,1){\makebox(0,0)[c]{\tiny $\beta_{\ov{1}}$}}
\put(17.15,1){\makebox(0,0)[c]{\tiny $\beta_0$}}
\put(19.5,1){\makebox(0,0)[c]{\tiny $\beta_{\frac{1}{2}}$}}
\put(21.5,1){\makebox(0,0)[c]{\tiny $\beta_{\frac{3}{2}}$}}
\end{picture}\vskip 8mm
\end{center}
Note that $\alpha_i=\beta_i$ for $i=\ov{m},\ldots,\ov{1}$.

We  assume that   $\h_{m+\infty}^*$ and $\h_{m|\infty}^*$ have symmetric bilinear forms $(\,\cdot\,|\,\cdot\,)$  given by
\begin{equation*}\label{bilinear form spo}
\begin{split}
&(\lambda|\delta_a)=s\big\langle (-1)^{|a|}E_a -K, \lambda \big\rangle, \ \ \ (\Lambda_{\ov{m}}|\Lambda_{\ov{m}})=0,
\end{split}
\end{equation*}
for $a,b\in \J_{m+\infty}$ or $\J_{m|\infty}$ and $\lambda\in  \h_{m+\infty}^*$ or $\h_{m|\infty}^*$. Here we assume $s=2$ for $\mf{g}=\mf{b}, \mf{b}^\bullet$, and $s=1$ otherwise. We have
$(\delta_a|\delta_b)=s(-1)^{|a|}\delta_{ab}$ for $a,b$, and hence $(\alpha_i|\alpha_i)$, $(\beta_j|\beta_j)\in 2\Z$ for $i\in I_{m+\infty}$ and $j\in I_{m|\infty}$.
Let
\begin{equation}\label{symmetrizing Cartan matrix}
s_i=
\begin{cases}
\ \ \, 1 & \text{if $i=\ov{m}$ and $\mf{g}=\mf{b}^\bullet, \mf{b}, \mf{d}$},\\
\ \ \, 2 & \text{if $i=\ov{m}$ and $\mf{g}=\mf{c}$},\\
\ \ \,2 & \text{if $i\in\{\,\ov{m-1},\ldots,\ov{1},0\,\}\cup\Z_{>0}$ and $\mf{g}=\mf{b}^\bullet, \mf{b}$},\\
\ \ \,1 & \text{if $i\in\{\,\ov{m-1},\ldots,\ov{1},0\,\}\cup\Z_{>0}$ and $\mf{g}=\mf{c}, \mf{d}$},\\
- 2 & \text{if $i\in \hf+\Z_{\geq 0}$  and $\mf{g}=\mf{b}^\bullet, \mf{b}$},\\
- 1 & \text{if $i\in \hf+\Z_{\geq 0}$  and $\mf{g}=\mf{c}, \mf{d}$}.\\
\end{cases}
\end{equation}
Then $s_i \langle \alpha_i^\vee, \lambda \rangle = (\alpha_i | \lambda)$ for $i\in I_{m+\infty}$, $\lambda\in \mf{h}^*_{m+\infty}$, and $s_j \langle \beta_j^\vee, \mu \rangle = (\beta_j | \mu)$ for $j\in I_{m|\infty}$, $\mu\in \mf{h}^*_{m|\infty}$.\vskip 2mm

For $n\geq 0$, we put $I_{m+n}=\{\,i\in I_{m+\infty}\,|\,(\alpha_i|\delta_a)\neq 0\ \text{for some $a\in\J_{m+n}$}\,\}$ and $I_{m|n}=\{\,i\in I_{m|\infty}\,|\,(\beta_i|\delta_a)\neq 0\ \text{for some $a\in\J_{m|n}$}\,\}$.
Let  $\mf{g}_{m+n}$ (resp. $\mf{g}_{m|n}$) be the subalgebra of $\mf{g}_{m+\infty}$ (resp. $\mf{g}_{m|\infty}$) generated by root vectors $E_{\pm\gamma}$ for $\gamma\in \Pi_{m+n}:=\{\,\alpha_i\,|\,i\in I_{m+n}\,\}$  (resp. $\Pi_{m|n}:=\{\,\beta_i\,|\,i\in I_{m|n}\,\}$) and $K$.
The Cartan subalgebra $\mf{h}_{m+n}$ (resp. $\mf{h}_{m|n}$) of $\mf{g}_{m+n}$ (resp. $\mf{g}_{m|n}$) is spanned by $K$ and $E_a$ for $a\in \J_{m+n}$ (resp. $\J_{m|n}$).

\section{Super duality and a semisimple tensor category of $\mf{g}_{m|n}$-modules}\label{Super duality and semisimple category}
Throughout the paper, a module $M$ over a superalgebra $U$ is understood to be a supermodule, that is, $M=M_0\oplus M_1$ with $U_iM_j\subset M_{i+j}$ for $i,j\in \Z_2$.
If $U$ has a comultiplication $\Delta$, then we have a $U$-module structure on $M\otimes N$ via $\Delta$ for $U$-modules $M$ and $N$, where we have a superalgebra structure on $U\otimes U$ with multiplication $(u_1\otimes u_2)(v_1\otimes v_2)=(-1)^{|u_2||v_1|}(u_1v_1)\otimes (u_2v_2)$ ($|u|$ denotes the parity of a homogeneous element $u\in U$).

\subsection{Super duality}\label{super duality review}
Let us briefly recall the {\it super duality} for ortho-symplectic Lie superalgebras  \cite{CLW}. Let $\cP$ denote the set of partitions. For $\lambda=(\lambda_i)_{i\geq 1}\in\cP$, let $\lambda'=(\lambda'_i)_{i\geq 1}$ be the conjugate of $\lambda$.

Let ${\mf l}_{m+\infty}$  be the standard Levi subalgebra of $\mf{g}_{m+\infty}$  corresponding to $\{\,\alpha_i\,|\,\in J_{m+\infty}\,\}$ for some $J_{m+\infty}$ with $\Z_{>0}\subset J_{m+\infty}\subset I_{m+\infty}\setminus\{0\}$. 
Let
\begin{equation*}
\begin{split}
P^+_{m+\infty}&=  \left\{\,\Lambda=c\Lambda_{\ov{m}}+\sum_{a\in \J_{m+\infty}}\lambda_a\delta_a\,\Bigg\vert
\begin{array}{l}
\text{(1) $c\in\C$ and $\lambda^+:=(\lambda_1,\lambda_2,\ldots) \in \cP$, } \\
\text{(2)  $\langle \alpha_i^\vee, \Lambda\rangle\in\Z_{\geq 0}$ for $i\in J_{m+\infty}\setminus \Z_{>0}$}
\end{array} \right\}
\end{split}
\end{equation*}
be the set of $\mf{l}_{m+\infty}$-dominant integral weights in $\mf{h}^*_{m+\infty}$. For $\Lambda\in P^+_{m+\infty}$, let  $L({\mf l}_{m+\infty},\Lambda)$ be the irreducible ${\mf l}_{m+\infty}$-module with highest weight $\Lambda$, and  $L({\mf g}_{m+\infty},\Lambda)$  the irreducible quotient of $K({\mf g}_{m+\infty},\Lambda):={\rm Ind}^{{\mf g}_{m+\infty}}_{{\mf p}}L({\mf l}_{m+\infty},\Lambda)$, where ${\mf p}$ is the subalgebra spanned by $K$, upper triangular matrices and ${\mf l}_{m+\infty}$, and $L({\mf l}_{m+\infty},\Lambda)$ is extended to a ${\mf p}$-module in a trivial way.

Let $\mathcal{O}(m+\infty)$ be the category of ${\mf g}_{m+\infty}$-modules $M$ satisfying
\begin{itemize}
\item[(1)] $M=\bigoplus_{\gamma\in \mf{h}^*_{m+\infty}}M_\gamma$ and $\dim M_\gamma<\infty$ for $\gamma\in {\mf h}^\ast_{m+\infty}$,

\item[(2)] ${\rm wt}(M)\subset \bigcup_{i=1}^r\left(\Lambda_i  -\sum_{\Pi_{m+\infty}}\Z_{\geq 0}\alpha \right)$ for some $r\geq 1$ and $\Lambda_i\in P^+_{m+\infty}$,

\item[(3)] $M$ decomposes  into a direct sum of $L({\mf l}_{m+\infty},\Lambda)$'s for $\Lambda\in P^+_{m+\infty}$.
\end{itemize}
Here $M_\gamma=\{\,m\,|\,h\cdot m =\langle h, \gamma\rangle \, m \ \ (h\in \mf{h}_{m+\infty})\,\}$ and  ${\rm wt}(M)=\{\,\gamma\in {\mf h}^\ast_{m+\infty}\,|\,M_\gamma\neq 0\,\}$ called the set of weights of $M$.
Note that (3) can be replaced by the condition that $E_{-\alpha_i}$ is locally nilpotent on $M$ for $i\in J_{m+\infty}$, that is, $M$ is an integrable ${\mf l}_{m+\infty}$-module, where $E_{-\alpha_i}$ denotes as usual a non-zero root vector associated to $-\alpha_i$
(see \cite[Section 2.5]{Kac78} for the case of $\mf{osp}(1|2m)$-modules).\vskip 3mm

Next, let   ${\mf l}_{m|\infty}$   be the standard Levi subalgebra of   $\mf{g}_{m|\infty}$  corresponding to  $\{\,\beta_i\,|\,i\in J_{m|\infty}\,\}$, where $J_{m|\infty}=\left(J_{m+\infty}\cap \{\ov{m},\ldots,\ov{1}\}\right)\cup \left(\tfrac{1}{2} +\Z_{\geq 0}\right)$. 
Let
\begin{equation*}
\begin{split}
P^+_{m|\infty}&=\left\{\,\Lambda^\natural\ \Big|\ \Lambda\in P^+_{m+\infty}\,\right\},
\end{split}
\end{equation*}
where $\Lambda^\natural=c\Lambda_{\ov{m}}+ \sum_{a=\ov{m}}^{\ov{1}}\lambda_a\delta_a +\sum_{b\in \hf+\Z_{\geq 0}}(\lambda^+)'_{b+\hf}\delta_b$ for $\Lambda=c\Lambda_{\ov{m}}+\sum_{a\in \J_{m+\infty}}\lambda_a\delta_a$.
Then we define $L({\mf l}_{m|\infty},\Lambda^\natural)$, $K({\mf{g}}_{m|\infty},\Lambda^\natural)$, and $L({\mf{g}}_{m|\infty},\Lambda^\natural)$ for $\Lambda\in P^+_{m+\infty}$ in the same way as in $\mc{O}(m+\infty)$.

Let  $\mathcal{O}(m|\infty)$ be the category of ${\mf g}_{m|\infty}$-modules $M$ satisfying
\begin{itemize}
\item[(1)] $M=\bigoplus_{\gamma\in \mf{h}^*_{m|\infty}}M_\gamma$ and $\dim M_\gamma<\infty$ for $\gamma\in {\mf h}^\ast_{m|\infty}$,

\item[(2)] ${\rm wt}(M)\subset \bigcup_{i=1}^r\left(\Lambda_i  -\sum_{\Pi_{m|\infty}}\Z_{\geq 0}\beta \right)$ for some $r\geq 1$ and  $\Lambda_i\in P^+_{m|\infty}$,

\item[(3)]$M$ decomposes  into a direct sum of $L({\mf l}_{m|\infty},\Lambda)$'s for $\Lambda\in P^+_{m|\infty}$.
\end{itemize}

\begin{rem}\label{Parity of weight spaces}{\rm
For $\Lambda=c\Lambda_{\ov{m}}+\sum_{a\in \J_{m|\infty}}\lambda_a\delta_a\in {\rm wt}(M)$, we assume that the parity of  $\Lambda$ is $\sum_{a\geq \hf}\lambda_a$ $\pmod 2$ when $\mf{g}\neq \mf{b}^\bullet$, and $\sum_{a=\ov{m}}^{\ov{1}}\lambda_{a}$ $\pmod 2$ when $\mf{g}= \mf{b}^\bullet$, which we denote by $|\Lambda|$ (cf. \cite[Section 5.2]{CLW}). In particular, we have $|\beta_i|=0$ (resp. $1$) if and only if $\beta_i$ is even (resp. odd) for $i\in I_{m|\infty}$. We assume that the $\Z_2$-grading on $M\in\mc{O}(m|\infty)$ is induced from the parity of its weights.}
\end{rem}

Note that $\{\,L({\mf g}_{m+\infty},\Lambda)\,|\,\Lambda\in P^+_{m+\infty}\,\}$ and  $\{\,L({\mf g}_{m|\infty},\Lambda^\natural)\,|\,\Lambda\in P^+_{m+\infty}\,\}$ form complete lists of irreducibles (up to isomorphism) in $\mathcal{O}(m+\infty)$ and $\mathcal{O}(m|\infty)$, respectively.
By \cite[Theorems 4.6 and 5.4]{CLW}, we have the following, which is called {\it super duality}.

\begin{thm}\label{super duality}
There exists an  equivalence of categories $\mc{F} :  \mathcal{O}(m+\infty) \longrightarrow \mathcal{O}(m|\infty)$ such that $\mc{F}(L({\mf g}_{m+\infty},\Lambda)) \cong L({\mf{g}}_{m|\infty},\Lambda^\natural)$  for $\Lambda\in P^+_{m+\infty}$
\end{thm}

Let us give a brief description of $\mc{F}$ for the readers' convenience.
The super duality \cite{CLW} is indeed established by considering an intermediate category between $\mc{O}(m+\infty)$ and $\mc{O}(m|\infty)$, which plays a crucial role.

Let $\td{\mf{b}}^\bullet$, $\td{\mf{c}}$ and $\td{\mf{b}}$, $\td{\mf{d}}$ be the central extensions of $\mf{spo}(V_{\I})$ and $\mf{osp}(V_\I)$ with respect to the 2-cocylcle $\alpha$ for $\I=\td{\I}_m$, $\td{\I}_m^\times$, respectively.
Let $\td{\mf{h}}$ be the Cartan subalgebra of $\td{\mf{g}}$ with basis $\{\,K,\, E_{a} \ (a\in \td{\I}^+_m)\,\}$, and $\td{\mf{h}}^*$ the restricted dual spanned by $\{\,\Lambda_{\ov{m}},\, \delta_{a} \ (a\in \td{\I}^+_m)\,\}$.
The set of  simple roots of $\td{\mf g}$ is $\td{\Pi}=\{\,\gamma_i\,|\,i\in \td{I} \,\}$, where $\td{I}=\{\ov{m},\ldots,\ov{1},0\} \cup \hf\Z_{>0}$, and
$\gamma_{i}=\beta_{i}$ ($i=\ov{m},\ldots,\ov{1},0$), $\gamma_{j}=\delta_j-\delta_{j+\hf}$ ($r\in \hf\Z_{>0}$). By definition, we have $\mf{g}_{m+\infty}$, $\mf{g}_{m|\infty}\subset \td{\mf{g}}$.

Let $\td{\mf l}$  be the standard Levi subalgebra of $\td{\mf{g}}$  corresponding to  $\{\,\gamma_i\,|\,i\in \td{J}\,\}$ with $\td{J}= \left(J_{m+\infty}\cap \{\ov{m},\ldots,\ov{1}\}\right)\cup \hf \Z_{>0}$, and let $\td{P}^+=\{\,\Lambda^\theta\,|\,\Lambda\in P^+_{m+\infty}\,\}$, where $\Lambda^\theta=c\Lambda_{\ov{m}}+\sum_{a=\ov{m}}^{\ov{1}}\lambda_a\delta_a +\sum_{b\in \hf\Z_{>0}}\theta(\lambda^+)_b\delta_b,$ with $\theta(\lambda^+)_{i-\hf}=\max\{\lambda'_i-i+1,0 \}$ and $\theta(\lambda^+)_i=\max\{\lambda_i-i,0\}$ ($i\in\Z_{>0}$)  for $\Lambda=c\Lambda_{\ov{m}}+\sum_{a\in \J_{m+\infty}}\lambda_a\delta_a$. Define $L(\td{\mf l} ,\Lambda^\theta)$, $K(\td{\mf{g}},\Lambda^\theta)$, and $L(\td{\mf{g}},\Lambda^\theta)$ for $\Lambda\in   P^+_{m+\infty}$ as in $\mc{O}(m+\infty)$ and $\mc{O}(m|\infty)$. Let $\td{\mc{O}}$ be the category of $\td{\mf{g}}$-modules $M$ such that (1) $M$ is $\td{\mf h}$-semisimple with finite dimensional weight spaces, (2) ${\rm wt}(M)$ is dominated by a finite number of weights in $\td{P}^+$, (3) $M$ is a direct sum of $L(\td{\mf l} ,\Lambda^\theta)$'s  (see \cite[Section 3.2]{CLW} for more detail).

Now, for $M=\bigoplus_{\gamma \in \td{\mf{h}}^*}M_\gamma \in\td{\mc{O}}$, define
\begin{equation*}
T(M)=\bigoplus_{\gamma \in \mf{h}_{m+\infty}^*}M_\gamma, \ \ \ \
\ov{T}(M)=\bigoplus_{\gamma \in \mf{h}_{m|\infty}^*}M_\gamma.
\end{equation*}
Then we have equivalences of categories $T :  \td{\mc{O}} \longrightarrow \mc{O}(m+\infty)$ and  $\ov{T} :  \td{\mc{O}} \longrightarrow \mc{O}(m|\infty)$  \cite[Theorem 5.4]{CLW}
such that
\begin{equation}\label{T of K and L}
\begin{split}
&   \begin{cases} \ \, T(L(\td{\mf l},\Lambda^\theta))=L({\mf l}_{m+\infty},\Lambda),\\ T(X(\td{\mf{g}},\Lambda^\theta))=X(\mf{g}_{m+\infty},\Lambda),
\end{cases}\ \
 \begin{cases} \ \, \ov{T}(L(\td{\mf l},\Lambda^\theta))=L({\mf l}_{m|\infty},\Lambda^\natural),\\ \ov{T}(X(\td{\mf{g}},\Lambda^\theta))=X({\mf{g}}_{m|\infty},\Lambda^\natural), 
\end{cases}
\end{split}
\end{equation}
for $X=K, L$ (see also Remark \ref{remark on out setting}).
The functor $\mc{F}$ in Theorem \ref{super duality} is understood to be $\ov{T}\circ S$, where $S : \mc{O}(m+\infty) \longrightarrow \td{\mc{O}}$ is a functor such that $S\circ T$ (resp. $T\circ S$) is naturally isomorphic to the identity functor $1_{\td{\mc{O}}}$ (resp. $1_{\mc{O}(m+\infty)}$).
\begin{equation}\label{duality functors}
\xymatrixcolsep{3pc}\xymatrixrowsep{2pc}\xymatrix{
& \td{\mc{O}} \ar@{->}[dr]^-{\ov{T}}\ar@{->}[dl]_{T} &\\
\mathcal{O}(m+\infty) \ar@{->}[rr]^-{\mc{F}}
&  & \mathcal{O}(m|\infty)}
\end{equation}

\begin{rem}\label{remark on out setting}
{\rm  \mbox{}

(1) Our exposition is a special case of the results in \cite{CLW}, since we assume here that $m>0$.

(2) In \cite{CLW}, the authors use a central extension, say $\widehat{\mf{gl}}'(V_{\td{\I}_m})$, of ${\mf{gl}}(V_{\td{\I}_m})$  by a one-dimensional center $\C K$ with respect to the 2-cocycle $\alpha'(X,Y)={\rm Str}([J',X]Y)$, where $J'=E_{0,0}+\sum_{r\leq -1}E_{r,r}$ in order to describe in a unified way the truncation into modules over the ortho-symplectic Lie superalgebras of finite rank. On the other hand, our central extension $\widehat{\mf{gl}}(V_{\td{\I}_m})$ is given by using  $J=\sum_{r\leq 0}E_{r,r}$ to describe the fundamental weight $\Lambda_{\ov{m}}$ for $\mf{g}_{m+\infty}$ and $\mf{g}_{m|\infty}$. But, there is an isomorphism $\psi : \widehat{\mf{gl}}'(V_{\td{\I}_m}) \rightarrow \widehat{\mf{gl}}(V_{\td{\I}_m})$ given by $\psi(X)=X+{\rm Str}(J''X)K$ for $X\in \mf{gl}(V_{\td{\I}_m})$ and $\psi(K)=K$, where $J''=J-J'$ (cf. \cite[Section 2.4]{CLW}). So by  using $\psi$ one can translate the results in \cite{CLW} in terms of our setting without difficulty.

}
\end{rem}

\subsection{The category $\mathcal{O}^{int}(m|\infty)$}\label{subcategory of integrable modules}
For $\lambda=(\lambda_i)_{i\geq 1} \in \cP$ and $c\in \C$, put
\begin{equation*}
\begin{split}
\Lambda_{m+\infty}(\lambda,c)&=c\Lambda_{\ov{m}} + \lambda_1\delta_{\ov{m}}+\cdots+\lambda_m\delta_{\ov{1}}+\lambda_{m+1}\delta_{1}+\lambda_{m+2}\delta_{2}+\cdots,\\
\Lambda_{m|\infty}(\lambda,c)&=\Lambda_{m+\infty}(\lambda,c)^\natural.
\end{split}
\end{equation*}
Let $\cP(\mf{g})$ be given by
\begin{equation*}
\begin{split}
\cP(\mf{b}^\bullet)&=\{\,(\lambda,\ell)\in \cP\times\Z_{> 0}\,|\,\ell-2\lambda_1\in 2\Z_{\geq 0}\,\},\\
\cP(\mf{c})&=\{\,(\lambda,\ell)\in \cP\times\Z_{> 0}\,|\,\ell-\lambda_1\in \Z_{\geq 0}\,\},\\
\cP(\mf{b})&=\{\,(\lambda,\ell)\in \cP\times\Z_{> 0}\,|\,\ell-2\lambda_1\in \Z_{\geq 0}\,\},\\
\cP(\mf{d})&=\{\,(\lambda,\ell)\in \cP\times\Z_{> 0}\,|\,\ell-\lambda_1-\lambda_2\in \Z_{\geq 0}\,\}.
\end{split}
\end{equation*}

Let $\mathcal{O}^{int}(m+\infty)$ be a full subcategory of ${\mf g}_{m+\infty}$-modules $M$ in $\mathcal{O}(m+\infty)$ such that  a simple root vector $E_{-\alpha}$ is  locally nilpotent on $M$ for $\alpha\in\Pi_{m+\infty}$. It is well-known that $\mathcal{O}^{int}(m+\infty)$ is a semisimple tensor category, whose  irreducible objects are $L({\mf g}_{m+\infty},\Lambda_{m+\infty}(\lambda,\ell))$ for $(\lambda,\ell)\in \cP(\mf{g})$ \cite{Kac90}  (see \cite[Section 2.5]{Kac78} when $\mf{g}=\mf{b}^\bullet$).

Now, we want to characterize the $\mf{g}_{m|\infty}$-modules which correspond to integrable $\mf{g}_{m+\infty}$-modules in $\mc{O}^{int}(m+\infty)$ under the super duality functor $\mc{F}$. Since  $E_{-\beta_0}^2=0$ and hence $E_{-\beta_0}$ is always locally nilpotent on $M\in \mc{O}(m|\infty)$, we do not necessarily obtain such $\mf{g}_{m|\infty}$-modules by the condition that $E_{-\beta}$ is locally nilpotent on $M$ for all $\beta\in \Pi_{m|\infty}$.

\begin{df}\label{Oint(m|infty)}{\rm
Define $\mathcal{O}^{int}(m|\infty)$ to be the category  of ${\mf g}_{m|\infty}$-modules $M$ satisfying the following conditions:
\begin{itemize}
\item[($1$)] $M=\bigoplus_{\gamma\in \mf{h}^*_{m|\infty}}M_\gamma$ and $\dim M_\gamma<\infty$ for $\gamma\in {\mf h}^\ast_{m|\infty}$,
 
\item[($2$)] ${\rm wt}(M)\subset \bigcup_{i=1}^r\left(\ell_i\Lambda_{\ov{m}}+\sum_{a\in \J_{m|\infty}}\Z_{\geq 0}\delta_a \right)$ for some $r\geq 1$ and $\ell_i\in \Z_{\geq 0}$, 

\item[($3$)]  $E_{-\beta_{\ov{m}}}$ is locally nilpotent on $M$.
\end{itemize}}
\end{df}

Note that the conditions ($2$) and ($3$) imply that $E_{\pm\beta_{i}}$ is locally nilpotent on $M$ for $i\neq 0$, and hence $M$ is an integrable ${\mf l}_{m|\infty}$-module, where $J_{m|\infty}=I_{m|\infty}\setminus\{0\}$. 
Since any maximal weight of $M$ is of the form $\Lambda=\ell\Lambda_{\ov{m}}+\sum_{a\in \J_{m|\infty}}\lambda_a\delta_a \in P^+_{m|\infty}$ with $\ell>0$ and $\lambda_a\geq 0$, and $\Lambda\in \left(\ell\Lambda_{\ov{m}}-\sum_{\Pi_{m|\infty}}\Z_{\geq 0}\beta\right)\cup\left(\ell\Lambda_{\ov{m}}+\delta_{\ov{m}}-\sum_{\Pi_{m|\infty}}\Z_{\geq 0}\beta\right)$,  ${\rm wt}(M)$ is dominated by finitely many weights in $P^+_{m|\infty}$ by the condition ($2$). Then $M$ is a direct sum of $L(\mf{l}_{m|\infty},\Lambda)$'s for $\Lambda\in P^+_{m|\infty}$ again by \cite{Kac78,Kac90}. Therefore $M\in \mc{O}(m|\infty)$.
Equivalently, $\mathcal{O}^{int}(m|\infty)$ is a full subcategory of ${\mf g}_{m|\infty}$-modules $M$ in  $\mathcal{O}(m|\infty)$  with $\ov{m}\in J_{m|\infty}$  such that
\begin{equation}\label{polynomial weight}
{\rm wt}(M)\subset \Z_{\geq 0}\Lambda_{\ov{m}}+\sum_{a\in \J_{m|\infty}}\Z_{\geq 0}\delta_a.
\end{equation}

The goal of this section is to show that $\mathcal{O}^{int}(m|\infty)$ is a semisimple tensor category, which is equivalent to $\mathcal{O}^{int}(m+\infty)$ under $\mc{F}$.
For this, we need the following two lemmas.

Let $\mf{gl}_{m+\infty}$ be the standard Levi subalgebra of $\mf{g}_{m+\infty}$ corresponding to $\Pi_{m+\infty}\setminus\{\alpha_{\ov{m}}\}$. Let $\mf{h}^\circ_{m+\infty}$ be the Cartan subalgebra of $\mf{gl}_{m+\infty}$ with basis $\{\,E_a\,|\,a\in \J_{m+\infty}\,\}$, and
$(\mf{h}^\circ_{m+\infty})^*$ the restricted dual spanned by dual basis $\{\,\delta_a\,|\,a\in \J_{m+\infty}\,\}$.
Also, let $\mf{gl}_{m|\infty}$ be the standard Levi subalgebra of $\mf{g}_{m|\infty}$  corresponding to $\Pi_{m|\infty}\setminus\{\beta_{\ov{m}}\}$ with  $\mf{h}^\circ_{m|\infty}$ and  $(\mf{h}^\circ_{m|\infty})^*$ defined in a similar way.

For $\lambda\in\cP$, let $L(\mf{gl}_{m+\infty},\lambda)$ and $L(\mf{gl}_{m|\infty},\lambda^\natural)$ be the irreducible highest weight  modules over $\mf{gl}_{m+\infty}$ and $\mf{gl}_{m|\infty}$ with highest weights $\Lambda_{m+\infty}(\lambda,0)$ and $\Lambda_{m+\infty}(\lambda,0)^\natural$, respectively.

\begin{lem}\label{poly reps} Let $M\in \mathcal{O}^{int}(m|\infty)$. Then $M$ is a direct sum of $L(\mf{gl}_{m|\infty},\lambda^\natural)$'s for $\lambda\in\cP$.
\end{lem}
\pf By \eqref{polynomial weight} and \cite[Theorem 3.27]{CW} (cf. \cite[Section 3.2.2]{CK}), $M$ is a polynomial representation of $\mf{gl}_{m|\infty}$ and hence completely reducible.
That is, $M$ is a direct sum of $L(\mf{gl}_{m|\infty},\lambda^\natural)$'s for $\lambda\in\cP$.
\qed\vskip 2mm

Let $\mc{G}: \mathcal{O}(m|\infty) \longrightarrow \mathcal{O}(m+\infty)$ be an equivalence of categories such that $\mc{G}\circ\mc{F}$ (resp. $\mc{F}\circ\mc{G}$) is naturally isomorphic to ${\rm id}_{\mathcal{O}(m+\infty)}$ (resp. ${\rm id}_{\mathcal{O}(m|\infty)}$).

\begin{lem}\label{poly reps-2}{\rm
Let $M\in \mathcal{O}^{int}(m|\infty)$. Then $\mc{G}(M)$  is a direct sum of $L(\mf{gl}_{m+\infty},\lambda)$'s for $\lambda\in\cP$.
}
\end{lem}
\pf Let us briefly recall the super duality for general linear Lie superalgebras \cite{CL} (with respect to a maximal Levi subalgebra). Let $P^{+'}_{{m+\infty}}$ be the set of weights $\Lambda=\sum_{a\in \J_{m+\infty}}\lambda_a\delta_a\in \mf{h}^*_{m+\infty}$ such that  $\lambda^+=(\lambda_1,\lambda_2,\ldots) \in \cP$ and $\langle \alpha_i^\vee, \Lambda\rangle\in\Z_{\geq 0}$ for $i\in  \{\ov{m-1},\ldots,\ov{1}\}$.
Let $\mathcal{O}_{\mf{gl}_{m+\infty}}$ be the category of $\mf{gl}_{m+\infty}$-modules $M$ satisfying
\begin{itemize}
\item[(1)] $M=\bigoplus_{\gamma\in (\mf{h}^\circ_{m+\infty})^*}M_\gamma$ and
$\dim M_\gamma<\infty$ for $\gamma\in (\mf{h}^\circ_{m+\infty})^*$,

\item[(2)] ${\rm wt}(M)\subset \bigcup_{i=1}^r\left(\Lambda_i  -\sum_{\Pi_{m+\infty}\setminus\{\alpha_{\ov{m}\}}}\Z_{\geq 0}\alpha\right)$ for some $r\geq 1$ and  $\Lambda_i\in P^{+'}_{{m+\infty}}$,

\item[(3)] $M$ decomposes into a direct sum of $L(\mf{l}'_{m+\infty},\Lambda)$'s for $\Lambda\in P^{+'}_{{m+\infty}}$,
\end{itemize}
where ${\rm wt}(M)$ is the set of weights of $M$ with respect to $\mf{h}^\circ_{m+\infty}$, and $\mf{l}'_{m+\infty}$ is the standard Levi subalgebra corresponding to $\Pi_{m+\infty}\setminus\{\alpha_{\ov{m}},\alpha_0\}$.
Also we define  the category $\mathcal{O}_{\mf{gl}_{m|\infty}}$ of $\mf{gl}_{m|\infty}$-modules in a similar fashion.

Let $\td{\mf{gl}}$ be the standard Levi subalgebra of $\td{\mf{g}}$ corresponding to $\td{\Pi}\setminus\{\gamma_{\ov{m}}\}$, and let $\mc{O}_{\td{\mf{gl}}}$ be a  parabolic category of $\td{\mf{gl}}$-modules with respect to $\td{\Pi}\setminus\{\,\gamma_{\ov{m}},\gamma_0\,\}$ defined in the same way as $\mathcal{O}_{\mf{gl}_{m+\infty}}$ and $\mathcal{O}_{\mf{gl}_{m|\infty}}$ (see \cite{CL} for more detail). Let $\td{\mf{h}}^\circ$ be the Cartan subalgebra of $\td{\mf{gl}}$ with basis $\{\,E_a\,|\,a\in \td{\I}^+_m\,\}$ and $(\td{\mf{h}}^\circ)^*$ its restricted dual spanned by $\{\,\delta_a\,|\,a\in \td{\I}^+_m\,\}$.
By \cite[Theorem 5.1]{CL} (see also \cite[Theorem 6.38]{CW}), we have equivalences of categories
\begin{equation}\label{duality functors for gl}
\xymatrixcolsep{4pc}\xymatrixrowsep{4pc}\xymatrix{
\mathcal{O}_{\mf{gl}_{m+\infty}} \ar@{<-}[r]^-{T'}
& \mc{O}_{\td{\mf{gl}}} \ar@{->}[r]^-{\ov{T}'} & \mathcal{O}_{\mf{gl}_{m|\infty}}},
\end{equation}
where  for $M=\bigoplus_{\gamma \in (\td{\mf{h}}^\circ)^*}M_\gamma \in \mc{O}_{\td{\mf{gl}}}$,
\begin{equation*}
T'(M)=\bigoplus_{\gamma \in (\mf{h}_{m+\infty}^\circ)^*}M_\gamma, \ \ \ \
\ov{T}'(M)=\bigoplus_{\gamma \in (\mf{h}_{m|\infty}^\circ)^*}M_\gamma.
\end{equation*}

Now, let $\ov{M}\in \mathcal{O}(m|\infty)$ be given. Since $\ov{M}\cong \ov{T}(\td{M})$ for some $\td{M}\in \td{\mc{O}}$ by \eqref{duality functors}, let us identify $\ov{M}$ with $T(\td{M})$ and put $M=T(\td{M})$. Note that $M\cong \mc{G}(\ov{M})$ by Theorem \ref{super duality}.  By definition of $\td{\mc{O}}$, we have ${\rm wt}(\td{M})\subset \bigcup_{i=1}^r D(\Lambda_i)$ for some  $\Lambda_1,\ldots,\Lambda_r\in \td{P}^+$,
where $D(\Lambda_i) =\Lambda_i  -\sum_{\td{\Pi}}\Z_{\geq 0}\gamma$. We may assume that $D(\Lambda_i)\cap D(\Lambda_j)=\emptyset$ for $i\neq j$.

For $\Lambda\in \td{P}^+$ and $k\geq 1$, put $D(\Lambda)_k=\Lambda  - k\gamma_{\ov{m}}- \sum_{\td{\Pi}\setminus\{\gamma_{\ov{m}}\}}\Z_{\geq 0}\gamma$.
Then we have
$\td{M}=\bigoplus_{k\geq 0}\td{M}_k$,
where $\td{M}_k$ is the sum of  $\td{M}_\mu$ over  $\mu\in {\rm wt}(\td{M})$ such that $\mu\in  \bigcup_{i=1}^r D(\Lambda_i)_k$. We can check that $\td{M}_k$ is $\td{\mf{gl}}$-invariant and $\td{M}_k\in \mc{O}_{\td{\mf{gl}}}$ since its weights are finitely dominated as an $\td{\mf{gl}}$-module. Hence  we have
\begin{equation}\label{k grading}
\ov{M}=\bigoplus_{k\geq 0}\ov{M}_k, \ \ \ \ M=\bigoplus_{k\geq 0}M_k,
\end{equation}
where  $\ov{M}_k:=\ov{T}'(\td{M}_k)\in \mathcal{O}_{\mf{gl}_{m|\infty}}$, and  $M_k:={T}'(\td{M}_k)\in\mathcal{O}_{\mf{gl}_{m+\infty}}$. Also, we observe that $\ov{T}$ and $T$ when restricted on $\td{M}_k$ coincide with $\ov{T}'$ and $T'$, respectively.

Now, suppose that $\ov{M}\in \mc{O}^{int}(m|\infty)$. By Lemma \ref{poly reps}, $\ov{M}_k$ is a
direct sum of $L(\mf{gl}_{m|\infty},\lambda^\natural)$'s $(\lambda\in\cP)$. By \eqref{duality functors for gl}, $M_k$ is a direct sum of $L(\mf{gl}_{m+\infty},\lambda)$'s $(\lambda\in\cP)$ with the same multiplicity  as $\ov{M}_k$ for each $\lambda$. Finally, it follows  from \eqref{k grading} that  $M\cong \mc{G}(\ov{M})$ is a direct sum of $L(\mf{gl}_{m+\infty},\lambda)$'s $(\lambda\in\cP)$.\qed

\begin{thm}\label{complete reducibility for integrable super}
$\mathcal{O}^{int}(m|\infty)$ is a semisimple tensor category, which is equivalent to $\mathcal{O}^{int}(m+\infty)$ under $\mc{F}$, with  irreducible objects $L({\mf g}_{m|\infty},\Lambda_{m|\infty}(\lambda,\ell))$ for $(\lambda,\ell)\in \cP(\mf{g})$.
\end{thm}
\pf Let $(\lambda,\ell)\in \cP(\mf{g})$ be given. Since $\mc{F}(L({\mf g}_{m+\infty},\Lambda_{m+\infty}(\lambda,\ell)))$ is isomorphic to $L({\mf g}_{m|\infty},\Lambda_{m|\infty}(\lambda,\ell))$ and $L({\mf g}_{m+\infty},\Lambda_{m+\infty}(\lambda,\ell))$ is a direct sum of  $L(\mf{gl}_{m+\infty},\lambda)$'s $(\lambda\in\cP)$, $L({\mf g}_{m|\infty},\Lambda_{m|\infty}(\lambda,\ell))$ is a direct sum of $L(\mf{gl}_{m|\infty},\lambda^\natural)$'s $(\lambda\in\cP)$, which implies that $L({\mf g}_{m|\infty},\Lambda_{m|\infty}(\lambda,\ell))\in \mc{O}^{int}(m|\infty)$. 

Suppose that $M\in \mathcal{O}^{int}(m|\infty)$ is given. By the conditions ($2$) and  ($3$) in Definition \ref{Oint(m|infty)}, we may assume that $M\in \mc{O}(m|\infty)$ with $J_{m|\infty}=I_{m|\infty}\setminus\{0\}$. We first claim that $\mc{G}(M)\in \mathcal{O}^{int}(m+\infty)$.
By definition of $\mathcal{O}(m+\infty)$, $\mc{G}(M)$ is an integrable ${\mf l}_{m+\infty}$-module. Hence $E_{-\alpha_i}$ is locally nilpotent on $\mc{G}(M)$ for $i\neq 0$.
On the other hand, by Lemma \ref{poly reps-2}, $\mc{G}(M)$ is a direct sum of $L(\mf{gl}_{m+\infty},\lambda)$'s $(\lambda\in\cP)$. In particular, $E_{-\alpha_0}$ acts locally nilpotently on $\mc{G}(M)$. Therefore, $\mc{G}(M)\in \mc{O}^{int}(m+\infty)$.
Since $\mc{O}^{int}(m+\infty)$ is semisimple, we have
$\mc{G}(M)\cong \bigoplus_{(\lambda,\ell)\in\cP(\mf{g})}L({\mf{g}}_{m+\infty},\Lambda_{m+\infty}(\lambda,\ell))^{\oplus m_{(\lambda,\ell)}}$,
for some $m_{(\lambda,\ell)}\in\Z_{\geq 0}$. Then applying $\mc{F}$ to $\mc{G}(M)$, we get
\begin{equation*}
M\cong \bigoplus_{(\lambda,\ell)\in\cP(\mf{g})}L({\mf{g}}_{m|\infty},\Lambda_{m|\infty}(\lambda,\ell))^{\oplus m_{(\lambda,\ell)}}.
\end{equation*}
Hence $M$ is semisimple. Finally, given $M_1, M_2\in \mc{O}^{int}(m|\infty)$, it is not difficult to check that $M_1\otimes M_2\in \mc{O}^{int}(m|\infty)$. Therefore, $\mc{O}^{int}(m|\infty)$ is a semisimple tensor category. It is clear that the functor $\mc{F}$ when restricted to $\mc{O}^{int}(m+\infty)$ gives an equivalence between  $\mc{O}^{int}(m+\infty)$ and $\mc{O}^{int}(m|\infty)$.
\qed

\begin{cor}\label{highest weight module for integrable super}
Let $M$ be a highest weight $\mf{g}_{m|\infty}$-module in $\mathcal{O}^{int}(m|\infty)$. Then $M$ is  isomorphic to $L(\mf{g}_{m|\infty},\Lambda_{m|\infty}(\lambda,\ell))$ for some $(\lambda,\ell)\in \cP(\mf{g})$.
\end{cor}

\subsection{The category $\mc{O}^{int}(m|n)$}\label{category of integrable modules over finite rank}
Consider the representations of $\mf{g}_{m|n}$ of finite rank (see \cite[Sections 3.3 and 3.4]{CLW} for more detail).
 Let $\mf{l}_{m|n}=\mf{l}_{m|\infty}\cap\mf{g}_{m|n}$. Let $P^+_{m|n}$ be the set of $\Lambda\in P^+_{m|\infty}$ such that $\langle E_a,\Lambda\rangle =0$ for $a\not\in \J_{m|n}$, and $\cP(\mf{g})_{m|n}$ the set of $(\lambda,\ell)\in \cP(\mf{g})$ such that $\Lambda_{m|\infty}(\lambda,\ell)\in P^+_{m|n}$. (We assume that $\cP(\mf{g})_{m|\infty}=\cP(\mf{g})$.)  Let us write $\Lambda_{m|n}(\lambda,\ell)=\Lambda_{m|\infty}(\lambda,\ell)$ for $\Lambda_{m|\infty}(\lambda,\ell)\in P^+_{m|n}$.
As in Section \ref{super duality review}, one may define  $L(\mf{l}_{m|n},\Lambda)$, $K(\mf{g}_{m|n},\Lambda)$,  $L(\mf{g}_{m|n},\Lambda)$ for $\Lambda\in P^+_{m|n}$, and a parabolic category  $\mc{O}(m|n)$ of $\mf{g}_{m|n}$-modules.

Let $M\in \mc{O}(m|\infty)$ be given with $M=\bigoplus_{\gamma}M_\gamma$. We define the truncation functor $\mf{tr}_n : \mc{O}(m|\infty) \longrightarrow \mc{O}(m|n)$ by
\begin{equation*}
\mf{tr}_n(M) =\bigoplus_{\gamma}M_\gamma,
\end{equation*}
where the sum is over $\gamma$ with $\langle E_a,\gamma\rangle=0$ for $a\not\in \J_{m|n}$. For $\Lambda\in P^+_{m|\infty}$, we have $\mf{tr}_n(X(\mf{g}_{m|\infty},\Lambda))=X(\mf{g}_{m|n},\Lambda)$ if $\Lambda\in P^+_{m|n}$, and $0$ otherwise for $X=K, L$ by \cite[Lemma 3.2]{CLW}.
Also, it is easy to see that $\mf{tr}_n(L(\mf{l}_{m|\infty},\Lambda))=L(\mf{l}_{m|n},\Lambda)$ if  $\Lambda\in P^+_{m|n}$, and $0$ otherwise.

\begin{df}\label{Oint(m|n)}{\rm
Define $\mathcal{O}^{int}(m|n)$ to be the category  of ${\mf g}_{m|n}$-modules $M$ satisfying
\begin{itemize}
\item[($1$)] $M=\bigoplus_{\gamma\in \mf{h}_{m|n}^*}M_\gamma$ and $\dim M_\gamma<\infty$ for $\gamma\in  \mf{h}_{m|n}^*$,

\item[($2$)] ${\rm wt}(M)\subset \bigcup_{i=1}^r\left(\ell_i\Lambda_{\ov{m}}+\sum_{a\in \J_{m|n}}\Z_{\geq 0}\delta_a \right)$ for some $r\geq 1$ and $\ell_i\in \Z_{\geq 0}$,

\item[($3$)] $E_{-\beta_{\ov{m}}}$ is locally nilpotent on $M$.
\end{itemize}}
\end{df}
As in the case of $\mc{O}^{int}(m|\infty)$, each $M\in \mathcal{O}^{int}(m|n)$ is a direct sum of $L(\mf{l}_{m|n},\Lambda)$'s for $\Lambda\in P^+_{m|n}$, where $J_{m|n}=I_{m|n}\setminus\{0\}$, and $\mathcal{O}^{int}(m|n)$ is a full subcategory of $\mathcal{O}(m|n)$.

We will prove that $\mathcal{O}^{int}(m|n)$ is a semisimple tensor category. We remark that the super duality functor is not available to $\mathcal{O}^{int}(m|n)$. So we will prove this in a rather indirect way, but still using super duality. \vskip 3mm

\begin{lem}\label{lifting highest weight module} Let $V(\mf{g}_{m|n},\Lambda)$ be a highest weight $\mf{g}_{m|n}$-module in $\mc{O}(m|n)$ with highest weight $\Lambda\in P^+_{m|n}$. Then there exists a highest weight $\mf{g}_{m|\infty}$-module $V(\mf{g}_{m|\infty},\Lambda)$ with highest weight $\Lambda$ in $\mc{O}(m|\infty)$ such that $\mf{tr}_n(V(\mf{g}_{m|\infty},\Lambda))=V(\mf{g}_{m|n},\Lambda)$.
\end{lem}
\pf Note that $V(\mf{g}_{m|n},\Lambda)\cong K(\mf{g}_{m|n},\Lambda)/W_{n}$ for some $\mf{g}_{m|n}$-submodule $W_{n}$ of $K(\mf{g}_{m|n},\Lambda)$. By \cite[Lemma 2.1.10]{Kumar}, we have a filtration
$\{0\}= W_{n}^{(0)}\subset  W_{n}^{(1)}\subset  W_{n}^{(2)}\subset \cdots$
such that
\begin{itemize}
\item[(1)] $W_{n}=\bigcup_{i\geq 1}W_{n}^{(i)}$,

\item[(2)] $ W_{n}^{(i)}/ W_{n}^{(i-1)}$ is a highest weight $\mf{g}_{m|n}$-module with highest weight $\nu_i\in P^+_{m|n}$,

\item[(3)] if $\nu_i-\nu_j\in \sum_{\beta\in \Pi_{m|n}}\Z_{\geq 0}\beta$, then $i<j$,

\item[(4)] given $\gamma\in {\rm wt}(W_{n})$, $\left(W_{n}/W_{n}^{(i)} \right)_\gamma=0$ for $i\gg 0$.
\end{itemize}

We use induction to show that there exist $\mf{g}_{m|\infty}$-submodules $W^{(i)}$ of $ K(\mf{g}_{m|\infty},\Lambda)$ ($i\geq 0$) such that $\{0\}= W^{(0)}\subset  W^{(1)}\subset  W^{(2)}\subset \cdots$ and $\mf{tr}_n(W^{(i)})=W_n^{(i)}$ for $i\geq 1$.

Suppose that $i=1$. Let $v_n^{(1)}$ be a $\mf{g}_{m|n}$-highest weight vector of $W_n^{(1)}$.   Since we may regard $K(\mf{g}_{m|n},\Lambda)\subset K(\mf{g}_{m|\infty},\Lambda)$  with $\mf{tr}_n(K(\mf{g}_{m|\infty},\Lambda))=K(\mf{g}_{m|n},\Lambda)$ and $K(\mf{g}_{m|\infty},\Lambda)$ is an integrable $\mf{l}_{m|\infty}$-module, $v_n^{(1)}$ is also a $\mf{g}_{m|\infty}$-highest weight vector. We put $W^{(1)}=U(\mf{g}_{m|\infty}^{-})v_n^{(1)}$, where $\mf{g}_{m|\infty}^{-}$ is the subalgebra generated by $E_{-\beta}$ for $\beta\in \Pi_{m|\infty}$. By construction, it is clear that $\mf{tr}_n(W^{(1)})=W_n^{(1)}$

Suppose that there exist $W^{(1)}\subset \cdots \subset W^{(i-1)}$ such that $\mf{tr}_n(W^{(k)})=W_n^{(k)}$ for $k=1,\ldots,i-1$. Let $v_n^{(i)}$ be a $\mf{g}_{m|n}$-highest weight vector of $W_n^{(i)}/W_n^{(i-1)}$. Then by the same argument as above, $v_n^{(i)}$ is also a $\mf{g}_{m|\infty}$-highest weight vector. Let $W^{(i)}$ be the $\mf{g}_{m|\infty}$-submodule of $K(\mf{g}_{m|\infty},\Lambda)$ generated by $W^{(i-1)}$ and $v_n^{(i)}$. Note that $W^{(i)}/W^{(i-1)}=U(\mf{g}_{m|\infty}^-)v_n^{(i)}$ and hence $\mf{tr}_n(W^{(i)}/W^{(i-1)})=W_n^{(i)}/W_n^{(i-1)}$.
By the exactness of  $\mf{tr}_n$ and the induction hypothesis, we have $\mf{tr}_n(W^{(i)})=W_n^{(i)}$.
This completes the induction.

Now, if we put $W=\bigcup_{i\geq 1}W^{(i)}$, then $\mf{tr}_n(W)=W_n$. Since $\mf{tr}_n(K(\mf{g}_{m|\infty},\Lambda))=K(\mf{g}_{m|n},\Lambda)$ and $\mf{tr}_n$ is exact, it follows that $\mf{tr}_n(V(\mf{g}_{m|\infty},\Lambda))=V(\mf{g}_{m|n},\Lambda),$
where $V(\mf{g}_{m|\infty},\Lambda)=K(\mf{g}_{m|\infty},\Lambda)/W\in\mc{O}(m|\infty)$.
\qed

\begin{thm}\label{integrable highest weight module over finite rank}
Let $M$ be a highest weight $\mf{g}_{m|n}$-module in $\mathcal{O}^{int}(m|n)$. Then $M$ is isomorphic to $L({\mf g}_{m|n},\Lambda_{m|n}(\lambda,\ell))$ for some $(\lambda,\ell)\in \cP(\mf{g})_{m|n}$.
\end{thm}
\pf Let $v$ be a highest weight vector of $M$ with highest weight  $\gamma$. By the condition (2) in Definition \ref{Oint(m|n)} and \cite[Theorem 3.27]{CW}, $M$ is semisimple over $\mf{gl}_{m|n}:=\mf{gl}_{m|\infty}\cap\mf{g}_{m|n}$ and $v$ is also a highest weight vector of a polynomial module over $\mf{gl}_{m|n}$, which implies that $\gamma=\Lambda_{m|n}(\lambda,\ell)$ for some $\lambda\in\cP$ and $\ell\in\Z_{> 0}$. Furthermore, we have $(\lambda,\ell)\in\cP(\mf{g})_{m|n}$ by the condition (3) in Definition \ref{Oint(m|n)}.

Suppose that $M$ is not irreducible. Let $N$ be a proper maximal submodule of $M$. Choose a highest weight vector $v'$  of $N$ with highest weight $\eta$. By the same argument in the previous paragraph, we also have  $\eta=\Lambda_{m|n}(\mu,\ell')$ for some $(\mu,\ell')\in\cP(\mf{g})_{m|n}$. Since $\eta\in \gamma-\sum_{\beta\in\Pi_{m|n}}\Z_{\geq 0}\beta$, we have $\ell'=\ell$.

Next, by Lemma \ref{lifting highest weight module}, there exists a highest weight module $M_{m|\infty}\in\mc{O}(m|\infty)$ with highest weight $\Lambda_{m|\infty}(\lambda,\ell)=\Lambda_{m|n}(\lambda,\ell)\in P^+_{m|\infty}$ such that $\mf{tr}_n(M_{m|\infty})=M$. Here we may assume that $J_{m|\infty}=I_{m|\infty}\setminus\{0\}$. Then by \cite[Theorem 4.6 and Proposition 5.3]{CLW}, there exists a highest weight module $M_{m+\infty}\in\mc{O}(m+\infty)$ with highest weight $\Lambda_{m+\infty}(\lambda,\ell) \in P^+_{m+\infty}$ such that $\mc{F}(M_{m+\infty})=M_{m|\infty}$.

Note that $M$ is a semisimple $\mf{l}_{m|n}$-module and $v'$ generates a highest weight $\mf{l}_{m|n}$-submodule $L(\mf{l}_{m|n},\Lambda_{m|n}(\mu,\ell))$.
Since $M_{m|\infty}$ is a semisimple $\mf{l}_{m|\infty}$-module and
$\mf{tr}_n(L(\mf{l}_{m|\infty},\Lambda_{m|\infty}(\mu,\ell)))=L(\mf{l}_{m|n},\Lambda_{m|n}(\mu,\ell)),$
we conclude that the multiplicity of $L(\mf{l}_{m|\infty},\Lambda_{m|\infty}(\mu,\ell))$ in $M_{m|\infty}$ is non-zero. This implies that  the multiplicity of $L(\mf{l}_{m+\infty},\Lambda_{m+\infty}(\mu,\ell))$ in $M_{m+\infty}$ is also non-zero  by \eqref{T of K and L}. In particular, we have $\Lambda_{m+\infty}(\mu,\ell)\in \Lambda_{m+\infty}(\lambda,\ell)-\sum_{\alpha\in \Pi_{m+\infty}}\Z_{\geq 0}\alpha$.

Now, consider the eigenvalues of Casimir operator $\Omega$ on $M_{m+\infty}$. By \cite[Lemma 9.8]{Kac90}, we have
\begin{equation*}
(\Lambda_{m+\infty}(\lambda,\ell)+2\rho|\Lambda_{m+\infty}(\lambda,\ell))=(\Lambda_{m+\infty}(\mu,\ell)+2\rho|\Lambda_{m+\infty}(\mu,\ell)),
\end{equation*}
where $\rho$ is the Weyl vector for $\mf{g}_{m+\infty}$.
On the other hand, by \cite[Lemma 10.3]{Kac90},
\begin{equation*}
(\Lambda_{m+\infty}(\lambda,\ell)+2\rho|\Lambda_{m+\infty}(\lambda,\ell))>(\Lambda_{m+\infty}(\mu,\ell)+2\rho|\Lambda_{m+\infty}(\mu,\ell)),
\end{equation*}
which is a contradiction. Note that the arguments for the case of $\mf{g}=\mf{b}^\bullet$ is almost the same (see the proof of Theorem 1 in \cite[Section 2.1]{Kac78}).  Therefore, $M$ is an irreducible highest weight module with highest weight $\Lambda_{m|n}(\lambda,\ell)$.
\qed\vskip 2mm

For $M\in \mathcal{O}(m|n)$, let
$M^\vee=\bigoplus_{\gamma\in {\rm wt}(M)}M^*_\gamma$, where $M^*_\gamma={\rm Hom}_{\mathbb{C}}(M_\gamma,\C)$. 
Let $\tau$ be an anti-automorphism of ${\mf g}_{m|n}$ determined by $\tau(E_{\beta})= (-1)^{|\beta|} E_{-\beta}$, $\tau(E_{-\beta})=E_{\beta}$, and $\tau(h)=h$ for $\beta\in \Pi_{m|n}$ and $h\in {\mf h}_{m|n}$.
We define a $\mf{g}_{m|n}$-module structure on $M^\vee$ by
$\langle m,  u f\rangle  = (-1)^{|m||u|}\langle \tau(u) m,   f\rangle$,
for $f\in M^\vee$ and homogeneous elements  $u\in \mf{g}_{m|n}$, $m\in M$.  Then one can check that $M^\vee\in \mathcal{O}(m|n)$, and  $M \cong (M^\vee)^\vee$ where the isomorphism is given by mapping $m$ to $(-1)^{|\gamma|}m$ for $m\in M_\gamma$. Here we identify $(M^\vee)^\vee_\gamma$ with $M_\gamma$  as a vector space for each $\gamma$. We also have $L(\mf{g}_{m|n}, \Lambda)^\vee \cong L(\mf{g}_{m|n}, \Lambda)$ for $\Lambda\in P^+_{m|n}$.

\begin{thm}\label{complete reducibility for integrable super of finite rank}
$\mathcal{O}^{int}(m|n)$ is a semisimple tensor category, whose  irreducible objects are $L({\mf g}_{m|n},\Lambda_{m|n}(\lambda,\ell))$ for $(\lambda,\ell)\in \cP(\mf{g})_{m|n}$.
\end{thm}
\pf  Suppose that $M\in \mathcal{O}^{int}(m|n)$ is given. Choose a highest weight vector  $v$ with highest weight $\mu$, which is maximal (that is, $\mu+\beta\not\in {\rm wt}(M)$ for $\beta\in \Pi_{m|n}$), and let $N=U(\mf{g}_{m|n})v=U(\mf{g}_{m|n}^-)v$. By Theorem \ref{integrable highest weight module over finite rank}, we have $\mu=\Lambda_{m|n}(\lambda,\ell)$ for some  $(\lambda,\ell)\in \cP(\mf{g})_{m|n}$ and $N\cong L({\mf g}_{m|n},\Lambda_{m|n}(\lambda,\ell)).$

For $v^\ast \in M_\mu^*$ such that $\langle v,v^* \rangle\neq 0$, let $L=U(\mf{g}_{m|n})v^*=U(\mf{g}_{m|n}^-)v^* \subset M^\vee$. Since $M^\vee\in \mc{O}^{int}(m|n)$ and $\mu$ is maximal, $L$ is a highest weight module   in $\mc{O}^{int}(m|n)$ with highest weight $\Lambda_{m|n}(\lambda,\ell)$, and therefore,  $L \cong L({\mf g}_{m|n},\Lambda_{m|n}(\lambda,\ell))$ by Theorem \ref{integrable highest weight module over finite rank}.

By taking dual of the embedding $L\hookrightarrow M^\vee$ and then composing with $N\hookrightarrow M$, we have
\begin{equation}\label{left inverse over finite rank}
N \longrightarrow M\cong (M^\vee)^{\vee} \longrightarrow L^{\vee}.
\end{equation}
Since $v$ maps to a non-zero vector in $L^{\vee}$ and $L^{\vee}\cong N\cong L({\mf g}_{m|n},\Lambda_{m|n}(\lambda,\ell))$, \eqref{left inverse over finite rank} gives an isomorphism of $N$ onto itself by Schur's lemma, which implies that the short exact sequence
$0 \longrightarrow N \longrightarrow M \longrightarrow M/N \longrightarrow 0$
splits and $M\cong N\oplus M/N$. Therefore, $M$ is completely reducible.
Moreover, we see from Theorems \ref{complete reducibility for integrable super} and \ref{integrable highest weight module over finite rank} that $L({\mf g}_{m|n},\Lambda_{m|n}(\lambda,\ell))$ for $(\lambda,\ell)\in \cP(\mf{g})_{m|n}$ form a complete list of irreducibles in  $\mathcal{O}^{int}(m|n)$.

Finally, it is clear that $\mathcal{O}^{int}(m|n)$ is closed under tensor product. This completes the proof.
\qed


\section{Category $\mc{O}_q^{int}(m|n)$ over the quantum superalgebra $U_q({\mf g}_{m|n})$}\label{Semisimple category for quantum ortho-symplectic}
In this section, we consider the $q$-analogue of a module in $\mc{O}^{int}(m|n)$ over the quantized enveloping algebra $U_q(\mf{g}_{m|n})$, and prove its semisimplicity.

\subsection{The quantum superalgebra $U_q({\mf g}_{m|n})$}






From now on, we assume that $n\in\Z_{\geq 0}\cup\{\infty\}$.  Let $A=(a_{ij})=(\langle \beta_i^\vee,\beta_j\rangle)_{i,j\in I_{m|n}}$ the generalized Cartan matrix for $\mf{g}_{m|n}$. Let
\begin{equation*}
P_{m|n}=\bigoplus_{a\in \J_{m|n}}\Z \delta_a \oplus \Z \Lambda_{\ov{m}},\ \
P^\vee_{m|n}=\bigoplus_{a\in \J_{m|n}}\Z E_a \oplus \Z r K\ \
\end{equation*}
be the weight lattice and dual weight lattice, respectively, where $r=1$ for $\mf{g}=\mf{c}$ and $r=2$ otherwise.

Let $q$ be an indeterminate. Put $q_i = q^{{\ov s}_i}$  for $i\in I_{m|n}$, and
\begin{equation*}
[r]_i=\frac{q_i^r-q_i^{-r}}{q_i-q_i^{-1}}, \ \ [r]_i!=\prod_{k=1}^r[k]_i,
\end{equation*}
for $r\geq 0$, where ${\ov s}_i=-s_i$   (see Remarks \ref{crystal for gl type} and \ref{crystal for gl type-2} for the difference when we use $\ov{s}_i$ instead of $s_i$ \eqref{symmetrizing Cartan matrix}).

The  quantum superalgebra $U_q({\mf g}_{m|n})$ is the
associative superalgebra (or $\Z_2$-graded algebra) with $1$ over $\mathbb{Q}(q)$ generated by
$e_i$, $f_i$ $(\,i\in I_{m|n}\,)$ and $q^h$ $(\,h\in P^{\vee}_{m|n}\,)$, which are
subject to the following relations \cite{Ya}:
{\allowdisplaybreaks
\begin{align*}
\ \ \ & {\rm deg}(q^h)=0,\ \  {\rm deg}(e_i)={\rm deg}(f_i)=|\beta_i|,\\
& q^0=1, \quad q^{h +h'}=q^{h}q^{h'}, \\ &q^h e_i=q^{\langle h,\beta_i\rangle}
e_i q^h, \quad q^h f_i=q^{-\langle h,\beta_i\rangle} f_i q^h, \\
& e_i f_j-(-1)^{|\beta_i||\beta_j|}f_j e_i =\delta_{ij}\frac{t_i-t_i^{-1}}{q_i-q^{-1}_i}, \\
& z_i z_j - (-1)^{|\beta_i||\beta_j|} z_j z_i =0,     \ \ \ \ \ \ \ \ \ \ \ \ \ \ \ \ \ \ \ \ \ \, \text{if $a_{ij} =0$},
\\ & \sum_{r=0}^{1+|a_{ij}|}(-1)^r z_i^{(r)} z_j z_i^{(1+|a_{ij}|-r)}= 0,  \, \ \ \ \ \ \ \ \ \ \ \ \  \text{if $i\neq 0$ and $a_{ij} \neq 0$}, \\
& z_0 z_{\ov{1}} z_0 z_{\hf} + z_{\ov{1}} z_0
z_{\hf} z_{0} + z_{0} z_{\hf} z_{0} z_{\ov{1}}  + z_{\hf} z_{0} z_{\ov{1}} z_{0} -(q+q^{-1}) z_{0} z_{\ov{1}} z_{\hf} z_{0} =0,
\end{align*}
}
\hskip -1.5mm for $i,j\in I_{m|n}$, $h, h' \in P^{\vee}_{m|n}$ and $z=e,f$, where $t_i = q^{{\ov s}_i \beta_i^\vee}$ and $z_i^{(r)}=\frac{z_i^r}{[r]_i!}$ for $r\geq 0$.

Let $U^+_q$ (resp. $U^-_q$) be the subalgebra of $U_q(\mf{g}_{m|n})$ generated by $e_i$ (resp. $f_i$) for $i\in I_{m|n}$ and $U^0_q$ the subalgebra generated by $q^h$ for $h\in P^\vee_{m|n}$. We have triangular decomposition $U_q(\mf{g}_{m|n}) \cong U_q^-\otimes U^0_q\otimes U^+_q$ as a $\mathbb{Q}(q)$-space.
There is a Hopf superalgebra structure on $U_q(\mf{g}_{m|n})$, where the comultiplication $\Delta$ is given by
\begin{equation*}
\begin{split}
\Delta(q^h)&=q^h\otimes q^h, \\ \Delta(e_i)&=e_i\otimes t_i^{-1} + 1\otimes e_i, \\
\Delta(f_i)&=f_i\otimes 1+ t_i\otimes f_i.
\end{split}
\end{equation*}
We will also need the following subalgebras of $U_q(\mf{g}_{m|n})$:
\begin{align*}
U_q(\mf{gl}_{m|n})&=\left\langle\, e_i, f_i, q^{\pm E_a}\,\big|\,i\in I_{m|n}\setminus\{\ov{m}\}, \ a\in \J_{m|n}\, \right\rangle,\\
U_q(\mf{gl}_{m|0})& =\left\langle\, e_i, f_i, q^{\pm E_a}\,\big|\,i\in I_{m|0}\setminus\{\ov{m}\}, \ a \in (\J_{m|n})_0 \, \right\rangle, \\
U_q(\mf{gl}_{0|n})&=\left\langle\, e_i, f_i, q^{\pm E_a}\,\big|\,i\in I_{0|n}, a\in (\J_{m|n})_1\, \right\rangle,\\
\end{align*}
where  $I_{m|0}=\{\,\ov{m},\ldots,\ov{1}\,\}$ and $I_{0|n}=I_{m|n}\setminus \{\ov{m},\ldots,\ov{1},0\}$.

\subsection{Classical limit}\label{classical limit}
We can define the notion of a highest weight $U_q(\mf{g}_{m|n})$-module by the triangular decomposition of $U_q(\mf{g}_{m|n})$, and consider its classical limit  in the same way as in the case of symmetrizable Kac-Moody algebras. We leave the detailed verification to the reader (see \cite{Ja,Lu}).

For a $U_q(\mf{g}_{m|n})$-module $M$ and $\gamma\in P_{m|n}$, let $M_\gamma=\{\,m\,|\,q^h m =q^{\langle h,\mu \rangle}m\ (h\in P^\vee_{m|n})\,\}\subset M$ and  ${\rm wt}(M)=\{\,\gamma\in P_{m|n}\,|\,M_\gamma\neq 0\,\}$.
Suppose that $M$ is a $U_q(\mf{g}_{m|n})$-module generated by a highest weight vector $u$ of weight $\Lambda\in P_{m|n}$. Then  $M=\bigoplus_{\mu}M_\mu$, where the sum is over $\mu \in \Lambda-\sum_{\beta\in \Pi_{m|n}}\Z_{\geq 0}\beta$.

Let ${\bf A}=\mathbb{Q}[q,q^{-1}]$. Let $M_{\bf A}$ be the ${\bf A}$-span of $f_{i_1}\ldots f_{i_r}u$ for $r\geq 0$ and $i_1,\ldots,i_r\in I_{m|n}$, and
$M_{\mu,{\bf A}}=M_{\bf A}\cap M_\mu$.
Then $M_{{\bf A}}=\bigoplus_{\mu}M_{\mu,{\bf A}}$, and
${\rm rank}_{\bf A}M_{\mu,{\bf A}}=\dim_{\mathbb{Q}(q)}M_\mu$.
One can check that the ${\bf A}$-module $M_{\bf A}$ is invariant under $e_i$, $f_i$, $q^{h}$ and $\frac{q^h-q^{-h}}{q-q^{-1}}$ for $i\in I_{m|n}$ and $h\in P^\vee_{m|n}$. Set $\ov{M}=M_{\bf A}\otimes_{\bf A}\mathbb{C}$ and
$\ov{M}_{\mu}=M_{\mu, {\bf A}}\otimes_{\bf A}\mathbb{C}$.
Here $\mathbb{C}$ is understood to be an ${\bf A}$-module where $f(q)\cdot c =f(1)c$ for $f(q)\in {\bf A}$ and $c\in\C$. We have $\ov{M}=\bigoplus_{\mu}\ov{M}_\mu$ with $\dim_{\mathbb{C}}\ov{M}_\mu={\rm rank}_{\bf A}M_{\mu,{\bf A}}$.

Recall that the enveloping algebra  $U(\mf{g}_{m|n})$ is isomorphic to the associative superalgebra with $1$ over $\C$ generated by $x^\pm_i$ ($i\in I_{m|n}$) and $h\in P^{\vee}_{m|n}$ subject to the following relations \cite[Theorem 10.5.8]{Ya}:
{\allowdisplaybreaks
\begin{align*}
\ \ \ \ \ \ \ & {\rm deg}(h)=0,\ \  {\rm deg}(x^\pm_i)=|\beta_i|,\\
& [h,h']=0, \ \  [h,x^\pm_i]=\pm{\langle h,\beta_i\rangle}x^\pm_i, \ \
[x^+_i,x^-_j] =\delta_{ij}\beta_i^\vee , \\
& [z_i,z_j] = 0,   \  \text{if $a_{ij} =0$}, \ \
 \sum_{r=0}^{1+|a_{ij}|}(-1)^r z_i^{(r)} z_j z_i^{(1+|a_{ij}|-r)}= 0,\
\text{if $i\neq 0$ and $a_{ij} \neq 0$}, \\
& z_0 z_{\ov{1}} z_0 z_{\hf} + z_{\ov{1}} z_0
z_{\hf} z_{0} + z_{0} z_{\hf} z_{0} z_{\ov{1}}  + z_{\hf} z_{0} z_{\ov{1}} z_{0} -2 z_{0} z_{\ov{1}} z_{\hf} z_{0} =0,
\end{align*}
}
\hskip -1.5mm for $i,j\in I_{m|n}$ and $h, h' \in P^{\vee}_{m|n}$, where $z=x^\pm$ and $z^{(r)}=\frac{z^r}{r!}$ for $r\geq 0$. Here $[\ ,\ ]$ denotes the superbracket $[u,v]=uv-(-1)^{|u||v|}vu$ for homogeneous elements $u,v\in U(\g_{m|n})$.

Let $\ov{e_i}$, $\ov{f_i}$, and $\ov{h}$ be the $\mathbb{C}$-linear endomorphisms on $\ov{M}$ induced from $e_i$, $f_i$, and $\frac{q^h-q^{-h}}{q-q^{-1}}$ for $i\in I_{m|n}$ and $h\in P^\vee_{m|n}$, and $\ov{U}_{\ov{M}}$ the subalgebra of ${\rm End}_{\mathbb{C}}(\ov{M})$ generated by them.
Then there exists a $\mathbb{C}$-algebra homomorphism from $U(\mf{g}_{m|n})$ to $\ov{U}_{\ov{M}}$ sending $x^+_{i}$,  $x^-_{i}$, and $h$ to $\ov{e_i}$, $\ov{f_i}$, and $\ov{h}$, respectively.
Hence $\ov{M}$ is a $U(\mf{g}_{m|n})$-module with highest weight $\Lambda$, which is called the classical limit of $M$.

\subsection{The category $\mathcal{O}_q^{int}(m|n)$}
\begin{df}\label{category O^int_q(m|n)}{\rm
Let $\mathcal{O}_q^{int}(m|n)$ be  the category of $U_q(\mf{g}_{m|n})$-modules $M$ satisfying
\begin{itemize}
\item[(1)] $M=\bigoplus_{\gamma\in P_{m|n}}M_\gamma$ and $\dim M_\gamma <\infty$ for $\gamma\in P_{m|n}$,

\item[(2)] ${\rm wt}(M)\subset \bigcup_{i=1}^r\left(\ell_i\Lambda_{\ov{m}}+\sum_{a\in \J_{m|n}}\Z_{\geq 0}\delta_a \right)$ for some $r\geq 1$ and  $\ell_i\in \Z_{\geq 0}$,

\item[(3)] $f_{\ov{m}}$ is locally nilpotent on $M$.
\end{itemize}}
\end{df}
For $\Lambda\in P_{m|n}$, let $L_q(\mf{g}_{m|n},\Lambda)$ denote the irreducible highest weight $U_q(\mf{g}_{m|n})$-module with highest weight $\Lambda$.

\begin{thm}\label{highest weight module for integrable quantum super}
Let $M$ be a highest weight $U_q(\mf{g}_{m|n})$-module in $\mathcal{O}^{int}_q(m|n)$. Then $M$ is  isomorphic to $L_q(\mf{g}_{m|n},\Lambda_{m|n}(\lambda,\ell))$ for some $(\lambda,\ell)\in \cP(\mf{g})_{m|n}$, and its classical limit is isomorphic to $L(\mf{g}_{m|n},\Lambda_{m|n}(\lambda,\ell))$.
\end{thm}
\pf By Section \ref{classical limit}, the classical limit  $\ov{M}$ is a highest weight module in $\mathcal{O}^{int}(m|n)$. By Corollary \ref{highest weight module for integrable super} and Theorem \ref{integrable highest weight module over finite rank}, $\ov{M}$ is isomorphic to $L({\mf g}_{m|n},\Lambda_{m|n}(\lambda,\ell))$ for some $(\lambda,\ell)\in \cP(\mf{g})_{m|n}$. Since $\dim_{\mathbb{Q}(q)} M_\gamma =\dim_\C \ov{M}_\gamma$ for $\gamma\in P_{m|n}$, this forces $M$ to be an irreducible highest weight module with highest weight $\Lambda_{m|n}(\lambda,\ell)$, that is, $M\cong L_q(\mf{g}_{m|n},\Lambda_{m|n}(\lambda,\ell))$. \qed\vskip 2mm

\begin{thm}\label{complete reducibility}
$\mathcal{O}^{int}_q(m|n)$ is a semisimple tensor category.
\end{thm}
\pf Since it is clear that $\mathcal{O}^{int}_q(m|n)$ is closed under tensor product, it suffices to show that it is semisimple.
For $M\in \mathcal{O}^{int}_q(m|n)$, let
$M^\vee=\bigoplus_{\lambda\in P}M^*_\lambda$ with $M^*_\lambda={\rm Hom}_{\mathbb{Q}(q)}(M_\lambda,\mathbb{Q}(q))$. 
Let $\tau_q$ be an anti-automorphism of $U_q({\mf g}_{m|n})$ determined by $\tau_q(e_i)= (-1)^{|\beta_i|} f_i$, $\tau_q(f_i)=e_i$, and $\tau_q(q^h)=q^h$ for $i\in I_{m|n}$ and $h\in {\mf h}_{m|n}$.
We define a $U_q(\mf{g}_{m|n})$-module structure on $M^\vee$ by
$\langle m,  u f\rangle  = (-1)^{|m||u|}\langle \tau_q(u) m,   f\rangle$,
for $f\in M^*$ and homogeneous elements  $u\in U_q(\mf{g}_{m|n})$, $m\in M$.  Then one can check that $M^\vee\in \mathcal{O}^{int}_q(m|n) $, and $M \cong (M^\vee)^\vee$ (cf. Section \ref{category of integrable modules over finite rank}). Also, if  $L_q(\mf{g}_{m|n},\Lambda_{m|n}(\lambda,\ell))\in  \mathcal{O}^{int}_q(m|n)$, then we have  $L_q(\mf{g}_{m|n},\Lambda_{m|n}(\lambda,\ell))^\vee \cong L_q(\mf{g}_{m|n},\Lambda_{m|n}(\lambda,\ell))$ for $(\lambda,\ell)\in \cP(\mf{g})_{m|n}$.
Now we can apply the same arguments as in Theorem \ref{complete reducibility for integrable super of finite rank} to prove that $M$ is completely reducible. \qed

The main goal in the rest of this paper is to show that $L_q({\mf g}_{m|n},\Lambda_{m|n}(\lambda,\ell))$ for $(\lambda,\ell)\in \cP(\mf{g})_{m|n}$ are the irreducible objects in $\mathcal{O}^{int}_q(m|n)$ and they have crystal bases when $\g={\mf b}, {\mf b}^\bullet, {\mf c}$.

\section{Crystal base of a $q$-deformed Fock space}\label{Crystal base of Fock spaces}
In this section, we recall the crystal base theory for contragredient Lie superalgebras \cite[Section 2.3]{BKK}, and consider a crystal base of a $q$-deformed Fock space $\mathscr{V}_q$ over $U_q(\mf{g}_{m|n})$.

\subsection{Crystal bases}\label{crystal base}
Let $M$ be a $U_q(\mf{g}_{m|n})$-module in $\mathcal{O}_q^{int}(m|n)$.

First, suppose that $i\in I_{m|n}\setminus\{0\}$ is given, that is, $\beta_i$ is an even simple root.
For $u\in M$ of weight $\lambda$, we have a unique expression
\begin{equation*}
u=\sum_{k\geq 0, -\langle \beta^\vee_i,\lambda \rangle}f_i^{(k)}u_k,
\end{equation*}
where $e_iu_k=0$ for all $k\ge 0$. Then we define the Kashiwara operators $\te_i$ and $\tf_i$ as follows:\vskip 2mm

\begin{itemize}
\item[(1)] For $i\in I_{m|0}$,
\begin{equation*}
\begin{split}
\te_i u =\sum_{k}q_i^{l_k-2k+1}f_i^{(k-1)}u_k, \ \ \ \
\tf_i u =\sum_{k}q_i^{-l_k+2k+1}f_i^{(k+1)}u_k,
\end{split}
\end{equation*}
where  $l_k=\langle \beta^\vee_i,\lambda+k\beta_i \rangle$ for $k\geq 0$.

\item[(2)] For  $i \in I_{0|n}$,
\begin{equation*}
\begin{split}
\te_i u =\sum_{k}f_i^{(k-1)}u_k, \ \ \ \
\tf_i u =\sum_{k}f_i^{(k+1)}u_k.
\end{split}
\end{equation*}
\end{itemize}

\noindent Next, suppose that $i=0$, that is, $\beta_0$ is an odd isotropic simple root.  We define
\begin{equation*}
\begin{split}
\te_0 u =e_0u, \ \ \ \
\tf_0 u =q_0 f_0t_0^{-1}u .
\end{split}
\end{equation*}\vskip 2mm

Now, let $\mathbb{A}$ denote the subring of $\mathbb{Q}(q)$ consisting
of all rational functions which are regular at $q=0$.
Then a pair $(L,B)$ is called a {\it crystal base of $M$} if
\begin{itemize}
\item[(1)] $L$ is an $\mathbb{A}$-lattice of $M$, where  $L=\bigoplus_{\lambda\in P_{m|n}}L_{\lambda}$ with $L_{\lambda}=L\cap
M_{\lambda}$,
\item[(2)] $\tilde{e}_i L\subset L$ and $\tilde{f}_i L\subset L$
for $i\in I_{m|n}$,

\item[(3)] $B$ is a pseudo-basis of $L/qL$ (i.e.
$B=B^{\bullet}\cup(-B^{\bullet})$ for a $\mathbb{Q}$-basis
$B^{\bullet}$ of $L/qL$),

\item[(4)] $B=\bigsqcup_{\lambda\in P_{m|n}}B_{\lambda}$ with
$B_{\lambda}=B\cap(L/qL)_{\lambda}$,

\item[(5)] $\tilde{e}_iB \subset B\sqcup \{0\}$,
$\tilde{f}_i B\subset B\sqcup \{0\}$ for $i\in I_{m|n}$,

\item[(6)] for $b,b'\in B$ and $k\in I_{m|n}$,
$\tilde{f}_i b = b'$ if and only if $b=\tilde{e}_i b'$.
\end{itemize}
The set $B / \{\pm 1 \}$ has an $I_{m|n}$-colored oriented graph
structure, where  $b\stackrel{i}{\rightarrow}
b'$ if and only if $\tf_i b=b'$ for $i\in I_{m|n}$ and  $\,b, b' \in B / \{\pm
1\}$. We call $B / \{\pm 1 \}$ the {\it crystal} of $M$.
For $b\in B$ and $i\in I_{m|n}$, we set
$\varepsilon_i(b)= \max\{\,r\in\mathbb{Z}_{\geq
0}\,|\,\tilde{e}_k^r b\neq 0\,\}$ and $\varphi_i(b)=
\max\{\,r\in\mathbb{Z}_{\geq 0}\,|\,\tilde{f}_k^r b\neq 0\,\}$. We denote the weight of $b$ by ${\rm wt}(b)$.

\begin{rem}\label{crystal for gl type}{\rm
A crystal base $(L,B)$ of a $U_q(\mf{g}_{m|n})$-module is also a crystal base as a $U_q(\mf{gl}_{m|n})$-module in the sense of \cite{BKK}, which in particular implies that it is a upper crystal base of $M$ as a $U_q(\mf{gl}_{m|0})$-module, and a lower crystal base of $M$ as a $U_q(\mf{gl}_{0|n})$-module (see \cite[Lemma 2.5]{BKK}).}
\end{rem}

Let $M_i$ $(\,i=1,2\,)$ be a
$U_q(\mf{g}_{m|n})$-module in $\mathcal{O}^{int}_q(m|n)$
with a crystal base $(L_i,B_i)$. Then  $(L_1\otimes L_2,B_1\otimes B_2)$ is a crystal
base of $M_1\otimes M_2$ \cite[Proposition 2.8]{BKK}. The actions of $\te_i$ and $\tf_i$ on $B_1\otimes B_2$ are give as follows.

For $i\in I_{0|n}$, we have {\allowdisplaybreaks
\begin{equation}\label{lower tensor product rule}
\begin{split}
&\tilde{e}_i(b_1\otimes b_2)= \begin{cases}
(\tilde{e}_i b_1) \otimes b_2, & \text{if $\varphi_i(b_1)\geq\varepsilon_i(b_2)$}, \\
b_1 \otimes (\tilde{e}_i b_2), & \text{if $\varphi_i(b_1)<\varepsilon_i(b_2)$},\\
\end{cases}
\\
&\tilde{f}_i(b_1\otimes b_2)=
\begin{cases}
(\tilde{f}_i b_1) \otimes b_2, & \text{if $\varphi_i(b_1)>\varepsilon_i(b_2)$}, \\
 b_1 \otimes (\tilde{f}_i b_2), & \text{if $\varphi_i(b_1)\leq\varepsilon_i(b_2)$}.
\end{cases}
\end{split}
\end{equation}}
For $i\in I_{m|0}$, we have
{\allowdisplaybreaks
\begin{equation}\label{upper tensor product rule}
\begin{split}
&\tilde{e}_i(b_1\otimes b_2)= \begin{cases}
 b_1 \otimes (\tilde{e}_ib_2), & \text{if $\varphi_i(b_2)\geq\varepsilon_i(b_1)$}, \\
(\tilde{e}_i b_1) \otimes b_2, & \text{if $\varphi_i(b_2)<\varepsilon_i(b_1)$},\\
\end{cases}
\\
&\tilde{f}_i(b_1\otimes b_2)=
\begin{cases}
b_1 \otimes (\tilde{f}_i b_2), & \text{if $\varphi_i(b_2)>\varepsilon_i(b_1)$}, \\
(\tilde{f}_i b_1) \otimes  b_2, & \text{if $\varphi_i(b_2)\leq\varepsilon_i(b_1)$}.
\end{cases}
\end{split}
\end{equation}}
For $i=0$, we have{\allowdisplaybreaks
\begin{equation}\label{tensor product rule for 0}
\begin{split}
\tilde{e}_0(b_1\otimes b_2)=&
\begin{cases}
\pm b_1\otimes (\tilde{e}_0 b_2), & \text{if }\langle \beta^\vee_0,{\rm
wt}(b_2)\rangle>0, \\ (\tilde{e}_0 b_1)\otimes  b_2, & \text{if
}\langle \beta^\vee_0,{\rm wt}(b_2)\rangle=0,
\end{cases}
\\
\tilde{f}_0(b_1\otimes b_2)=&
\begin{cases}
\pm b_1\otimes (\tilde{f}_0 b_2), & \text{if }\langle \beta^\vee_0,{\rm
wt}(b_2)\rangle>0, \\ (\tilde{f}_0 b_1)\otimes  b_2, & \text{if
}\langle \beta^\vee_0,{\rm wt}(b_2)\rangle=0,
\end{cases}
\end{split}
\end{equation}}
where $\pm$ is determined by the parity of $b_1$.

\subsection{$q$-deformed Clifford-Weyl algebra}
Let $\A_q$ be the associative $\mathbb{Q}(q)$-algebra with $1$ generated by $\psp_{a}$, $\psm_a$, $\om_a$, and $\om^{-1}_a$ for $a\in \pm\,\J_{m|n}$
subject to the following relations:
{\allowdisplaybreaks
\begin{align*}
&\om_a \om_b= \om_b\om_a, \ \ \om_a\om_a^{-1}=1, \\
&\om_a \psp_b \om_a^{-1} =
q^{(-1)^{|a|}\delta_{ab}}\psp_b ,\ \ \ \ \ \om_a \psm_b \om_a^{-1} =
q^{-(-1)^{|a|}\delta_{ab}}\psm_b  \\
&\psp_a\psp_b+(-1)^{|a||b|}\psp_b\psp_a=0,\ \ \ \ \psm_a\psm_b+(-1)^{|a||b|}\psm_b\psm_a=0 , \\
&\psp_a\psm_b+(-1)^{|a||b|}\psm_b\psp_a=0 \ \ \ \ (a\neq b), \\
&\psp_a\psm_a =[q\om_a],\ \ \
\psm_a\psp_a =
(-1)^{1+|a|}[\om_a].
\end{align*}}
Here $[q^k\om_a^{\pm 1}]=\frac{q^k\om^{\pm 1}_a-q^{-k}\om_a^{\mp 1}}{q-q^{-1}}$ for $k\in\Z$ and $a\in \pm\,\J_{m|n}$.
Note that the subalgebra generated by  $\psp_{a}$, $\psm_a$ and $\om^{\pm 1}_a$ for $a\in \pm(\J_{m|n})_0$ (resp. $\pm(\J_{m|n})_1$) is a $q$-deformed Clifford algebra (resp. $q$-deformed Weyl algebra) introduced by Hayashi \cite{Ha}. Let $\A_q^{-}$ (resp. $\A_q^{+}$) be the subalgebra generated by $\psp_a$, $\psm_a$ for $a\in -\J_{m|n}$ (resp. $a\in \J_{m|n}$).\vskip 2mm

Let $\F_q$ be the $\A_q$-module generated by $|0\rangle$ satisfying
\begin{equation*}
\psp_{-a} |0\rangle=\psm_b|0\rangle =0,\ \  \omega_{-a}|0\rangle =|0\rangle,\ \  \omega_b|0\rangle =q^{-1}|0\rangle,
\end{equation*}
for $a, b\in \J_{m|n}$. Let  $\F_q^{-}$ (resp. $\F_q^{+}$) be the $\A_q^{-}$-submodule (resp. $\A_q^{+}$-submodule) of $\F_q$ generated by $|0\rangle$.

Let $\bf{B}$ be the set of sequences ${\bf m}=(m_a)$ of non-negative integers indexed by $\pm\J_{m|n}$ such that $m_a\leq 1$ for $|a|=1$. For ${\bf m}=(m_a)\in \bf{B}$,  let
\begin{equation*}
\begin{split}
&\psi_{{\bf m}}=\overrightarrow{\prod_{a\in -\J_{m|n}}}{\psm_a}^{(m_a)}  \ \overrightarrow{\prod_{b\in\J_{m|n}}}{\psp_b}^{(m_b)},
\end{split}
\end{equation*}
where the product is taken in the order of $<$ on $\ov{\I}_{m}$ and $$\psp_a^{(r)}=\frac{(\psp_a)^{r}}{[r]!}, \ \ \ {\psm_a}^{(r)}=\frac{(\psm_a)^{r}}{[r]!},$$ with $[r]!=[r]\ldots[1]$ and $[k]=\frac{q^k-q^{-k}}{q-q^{-1}}$ for  $k, r\geq 0$.
By similar arguments as in \cite[Proposition 2.1]{Ha}, we can check that $\F_q$ is an irreducible $\A_q$-module with a $\mathbb{Q}(q)$-linear basis $\{\,\psi_{{\bf m}}|0\rangle\,|\,{\bf m}\in \B\,\}$.

Let $\B^{-}$ (resp. $\B^{+}$) be the set of ${\bf m}=(m_a)\in \B$ such that $m_a=0$ for $a\in \J_{m|n}$ (resp. $a\in -\J_{m|n}$).
Then $\{\,\psi_{{\bf m}}|0\rangle\,|\,{\bf m}\in \B^{-}\,\}$ and $\{\,\psi_{{\bf m}}|0\rangle\,|\,{\bf m}\in \B^{+}\,\}$ are $\mathbb{Q}(q)$-bases of $\F_q^{-}$ and $\F_q^{+}$, respectively.\vskip 2mm

Let us describe an action of  $U_q(\gl_{m|n})\subset U_q(\mf{g}_{m|n})$  on $\F_q$.
For $i\in I_{m|n}\setminus\{\ov{m}\}$, put
{\allowdisplaybreaks
\begin{align*}
{\bf t}^{+}_i = &
\begin{cases}
\om_{\ov{k+1}}^{-1}\om_{\ov{k}},  & \text{if $i = \ov{k}\in I_{m|0}$}, \\
\om_{\ov{1}}^{-1}\om_{\hf},  & \text{if $i = 0$}, \\
\om_{ i}^{-1}\om_{i+1},  & \text{if $i \in I_{0|n}$}, \\
\end{cases}\ \ \ \
{\bf t}^{-}_i =
\begin{cases}
\om_{-\ov{k}}^{-1}\om_{-\ov{k+1}},  & \text{if $i = \ov{k}\in I_{m|0}$}, \\
\om_{-\hf}^{-1}\om_{-\ov{1}},  & \text{if $i = 0$}, \\
\om_{-i-1}^{-1}\om_{-i},  & \text{if $i \in I_{0|n}$}, \\
\end{cases}\\
{\bf e}^{+}_i = &
\begin{cases}
\psp_{\ov{k+1}}\psm_{\ov{k}},  & \text{if $i = \ov{k}\in I_{m|0}$}, \\
\psp_{\ov{1}}\psm_{\hf},  & \text{if $i = 0$}, \\
\psp_{i}\psm_{i+1},  & \text{if $i \in  I_{0|n}$},
\end{cases}
\ \ \ \
{\bf e}^{-}_i =
\begin{cases}
\psp_{-\ov{k}}\psm_{-\ov{k+1}},  & \text{if $i = \ov{k}\in I_{m|0}$}, \\
\psp_{-\hf}\psm_{-\ov{1}},  & \text{if $i = 0$}, \\
\psp_{-i-1}\psm_{-i},  & \text{if $i \in  I_{0|n}$}, \\
\end{cases}
\\
{\bf f}^{+}_i = &
\begin{cases}
\psp_{\ov{k}}\psm_{\ov{k+1}},  & \text{if $i = \ov{k}\in I_{m|0}$}, \\
-\psp_{\hf}\psm_{\ov{1}},  & \text{if $i = 0$}, \\
\psp_{i+1}\psm_{ i},  & \text{if $i \in  I_{0|n}$},
\end{cases}  \ \ \ \
{\bf f}^{-}_i =
\begin{cases}
\psp_{-\ov{k+1}}\psm_{-\ov{k}},  & \text{if $i = \ov{k}\in I_{m|0}$}, \\
\psp_{-\ov{1}}\psm_{-\hf},  & \text{if $i = 0$}, \\
\psp_{-i}\psm_{-i-1},  & \text{if $i \in  I_{0|n}$}.
\end{cases}
\end{align*}}

We may regard ${\bf t}^{\pm}_i$, ${\bf e}^{\pm}_i$ and ${\bf f}^{\pm}_i $ as $\mathbb{Q}(q)$-linear operators on $\F_q$ under left multiplication.

\begin{lem}\label{rho<>0}
 $\F_{q^r}^{\pm}$ has a $U_q(\mf{gl}_{m|n})$-module structure  $\rho^{\pm} : U_q(\mf{gl}_{m|n}) \longrightarrow {\rm End}_{\mathbb{Q}(q)}(\F_{q^r}^{\pm})$ such that
\begin{align*}
\rho^+(q^{\pm E_a})=
\begin{cases}
q^{\pm 1}\om_a^{\pm \frac{1}{r}}, & \text{if $|a|=0$},\\
q^{\mp 1}\om_a^{\mp\frac{1}{r}}, & \text{if $|a|=1$},\\
\end{cases}&\   \ \
\rho^-(q^{\pm E_a})=
\begin{cases}
\om_{-a}^{\mp\frac{1}{r}}, & \text{if $|a|=0$},\\
\om_{-a}^{\pm \frac{1}{r}}, & \text{if $|a|=1$},\\
\end{cases} \\
\rho^{\pm}(e_i)= {\bf e}^{\pm}_i, & \ \ \ \rho^{\pm}(f_i)={\bf f}^{\pm}_i,
\end{align*}
for $a\in\J_{m|n}$ and $i\in I_{m|n}\setminus\{\ov{m}\}$. Here $r=2$ when $\mf{g}_{m|n}=\mf{b}_{m|n}$, and $r=1$ otherwise. We understand $\omega_a^{\pm\frac{1}{2}}$ as an operator given by $\omega_a^{\pm\frac{1}{2}}v=q^k v$ for $v\in \mathscr{F}_{q^2}^\pm$ with $\omega_a^{\pm 1}v=q^{2k} v$.
\end{lem}
\pf Put $\rho^{\pm}(t_i)={\bf t}^{\pm}_i$  for $i\in I_{m|n}\setminus\{\ov{m}\}$. We can check that for $i,j\in I_{m|n}\setminus\{\ov{m}\}$ and $h\in \bigoplus_{a\in \J_{m|n}}\Z E_a$
\begin{equation*}
\begin{split}
&{\bf e}^{\pm}_i{\bf f}^{\pm}_j -(-1)^{|\beta_i||\beta_j|} {\bf f}^{\pm}_j{\bf e}^{\pm}_i=\delta_{ij}\frac{{\bf t}^{\pm}_i-({\bf t}^{\pm}_i)^{-1}}{q_i-q_i^{-1}},\\
&{\bf e}^{\pm}_i{\bf e}^{\pm}_j -(-1)^{|\beta_i||\beta_j|} {\bf e}^{\pm}_j{\bf e}^{\pm}_i=
{\bf f}^{\pm}_i{\bf f}^{\pm}_j -(-1)^{|\beta_i||\beta_j|} {\bf f}^{\pm}_j{\bf f}^{\pm}_i=0\ \ \text{if $a_{ij}=0$},\\
&\rho^{\pm}(q^h){\bf e}^{\pm}_i \rho^{\pm}(q^{-h})= q^{\langle h,\beta_i \rangle}{\bf e}^{\pm}_i, \ \ \ \rho^{\pm}(q^h){\bf f}^{\pm}_i \rho^{\pm}(q^{-h})= q^{-\langle h,\beta_i \rangle}{\bf f}^{\pm}_i  
\end{split}
\end{equation*}
 (cf. \cite[Lemma 3.1]{Ha}).
Then we see that $\F_{q^r}^{\pm}$ is a $U_q(\mf{gl}_{m|0})\oplus U_q(\mf{gl}_{0|n})$-module by \cite[Proposition B.1]{KMPY} since $e_i, f_i$ are locally nilpotent on $\F_{q^r}^{\pm}$ for $i\in I_{m|n}\setminus\{\ov{m},0\}$.
Also, it is straightforward to check that ${\bf e}^{\pm}_i$ and ${\bf f}^{\pm}_i$ ($i=\ov{1}, 0, 1$) satisfy the other relevant relations in $U_q(\mf{gl}_{m|n})$. This implies that $\F_{q^r}^{\pm}$ is a $U_q(\mf{gl}_{m|n})$-module.\qed

\subsection{Crystal bases of $q$-deformed Fock spaces}
Let us put
\begin{equation}\label{Fock space}
\begin{split}
\mathscr{V}_q=
\begin{cases}
\F_q, & \text{if $\mf{g}=\mf{c}$},\\
\F^+_{q^2}, & \text{if $\mf{g}=\mf{b}$},\\
\F^+_q, & \text{if $\mf{g}=\mf{d}$},\\
\F^+_{q^2}\otimes\F^+_{q^2}, & \text{if $\mf{g}=\mf{b}^\bullet$}.\\
\end{cases}
\end{split}
\end{equation}

\begin{prop}\label{rho}
$\mathscr{V}_q$ has a $U_q(\mf{g}_{m|n})$-module structure  $\rho : U_q(\mf{g}_{m|n}) \longrightarrow {\rm End}_{\mathbb{Q}(q)}(\mathscr{V}_q)$ as follows: for $a\in\J_{m|n}$ and $i\in I_{m|n}$,
\begin{itemize}
\item[(1)] if $\mf{g}=\mf{c}$, then

{\allowdisplaybreaks
\begin{align*}
&\rho(q^{\pm E_a})= \rho^+(q^{\pm E_a})\rho^-(q^{\pm E_a}),
\ \ \ \ \rho(q^K)=q, \\
&\rho(e_i)   =
\begin{cases}
\psp_{-\ov{m}}\psm_{\ov{m}}, & \text{if $i = \ov{m}$}, \\
{\bf e}^{-}_i ({\bf t}^{+}_i)^{-1}+{\bf e}^{+}_i, & \text{if $i \neq \ov{m}$}, \\
\end{cases}
\ \ \ \
\rho(f_i)  =
\begin{cases}
\psp_{\ov{m}}\psm_{-\ov{m}}, & \text{if $i = \ov{m}$}, \\
{\bf f}^{-}_i +{\bf t}^{-}_i{\bf f}^{+}_i,  & \text{if $i \neq \ov{m}$}, \\
\end{cases}
\end{align*}}

\item[(2)] if $\mf{g}=\mf{b}$, then
{\allowdisplaybreaks
\begin{align*}
&\rho(q^{\pm E_a})= \rho^+(q^{\pm E_a}),
\ \ \ \ \rho(q^{2K})=q, \\
&\rho(e_i)   =
\begin{cases}
\psm_{\ov{m}}, & \text{if $i = \ov{m}$}, \\
{\bf e}^{+}_i, & \text{if $i \neq \ov{m}$}, \\
\end{cases}
\ \ \ \
\rho(f_i)  =
\begin{cases}
\psp_{\ov{m}}, & \text{if $i = \ov{m}$}, \\
{\bf f}^{+}_i,  & \text{if $i \neq \ov{m}$}, \\
\end{cases}
\end{align*}}

\item[(3)]  if $\mf{g}=\mf{d}$, then
{\allowdisplaybreaks
\begin{align*}
&\rho(q^{\pm E_a})= \rho^+(q^{\pm E_a}),
\ \ \ \ \rho(q^{2K})=q, \\
&\rho(e_i)   =
\begin{cases}
\psm_{\ov{m}}\psm_{\ov{m-1}}, & \text{if $i = \ov{m}$}, \\
{\bf e}^{+}_i, & \text{if $i \neq \ov{m}$}, \\
\end{cases}
\ \ \ \
\rho(f_i)  =
\begin{cases}
-\psp_{\ov{m}}\psp_{\ov{m-1}}, & \text{if $i = \ov{m}$}, \\
{\bf f}^{+}_i,  & \text{if $i \neq \ov{m}$}, \\
\end{cases}
\end{align*}}

\item[(4)]  if $\mf{g}=\mf{b}^\bullet$, then
{\allowdisplaybreaks
\begin{align*}
&\rho(q^{\pm E_a})= \rho^{+}(q^{\pm E_a})\otimes \rho^{+}(q^{\pm E_a}),
\ \ \ \ \rho(q^{2K})=q\otimes q, \\
&\rho(e_i)   =
\begin{cases}
\psm_{\ov{m}}\otimes \rho(t_{\ov{m}}^{-1})+1\otimes \psm_{\ov{m}}, & \text{if $i = \ov{m}$}, \\
{\bf e}^{+}_i ({\bf t}^{+}_i)^{-1}+{\bf e}^{+}_i, & \text{if $i \neq \ov{m}$}, \\
\end{cases}
\\
&\rho(f_i)  =
\begin{cases}
\left(\psp_{\ov{m}}\otimes 1 + \rho(t_{\ov{m}})\otimes \psp_{\ov{m}}\right)\sigma, & \text{if $i = \ov{m}$}, \\
{\bf f}^{+}_i +{\bf t}^{+}_i{\bf f}^{+}_i,  & \text{if $i \neq \ov{m}$}.
\end{cases}
\end{align*}}
Here $\rho$ is as in $(2)$ and $\sigma$ is a $\mathbb{Q}(q)$-linear operator on $\F^+_{q^2}\otimes\F^+_{q^2}$ defined by $$\sigma (\psi_{{\bf m}}|0\rangle\otimes \psi_{{\bf m}'}|0\rangle)=(-1)^{m_{\ov{m}}+m'_{\ov{m}}}\psi_{{\bf m}}|0\rangle\otimes \psi_{{\bf m}'}|0\rangle,$$ for ${\bf m}=(m_a)$ and ${\bf m}'=(m'_a)\in \B^+$.
\end{itemize}
\end{prop}
\pf The proof is almost the same as in Lemma \ref{rho<>0}. So, we leave the details to the readers.\qed

\begin{cor}\label{F as gl(m|n)-module}
We have $\F_q\cong \F_q^{-}\otimes \F_q^{+}$ as a module over $U_q(\mf{gl}_{m|n})\subset U_q({\mf c}_{m|n})$.
\end{cor}
\pf It follows immediately from comparing the actions of $U_q(\mf{gl}_{m|n})$ in Lemma \ref{rho<>0} and Proposition \ref{rho}.
\qed

\begin{cor}\label{tensor power of Fock space}
$\mathscr{V}_q$ is a $U_q(\g_{m|n})$-module in $\mc{O}^{int}_q(m|n)$. In particular,
$\mathscr{V}_q^{\otimes \ell}$ is completely reducible for $\ell\geq 1$.
\end{cor}
\pf By Proposition \ref{rho}, $\mathscr{V}_q\in \mc{O}^{int}_q(m|n)$ with ${\rm wt}(\mathscr{V}_q)\subset \ell\Lambda_{\ov{m}}+\sum_{a\in \J_{m|n}}\Z_{\geq 0}\delta_a$ with $\ell=1,2$. 
For example, when $\g={\mf c}$ we have ${\rm wt}(\psi_{{\bf m}}|0\rangle)=\Lambda_{\ov{m}}+\sum_{a\in \J_{m|n}}(m_a+m_{-a})\delta_a$ for ${\bf m}\in {\bf B}$.
By Theorem \ref{highest weight module for integrable quantum super}, $\mathscr{V}_q^{\otimes \ell}\in \mc{O}^{int}_q(m|n)$ for $\ell\geq 1$, and hence it is completely reducible.  \qed\vskip 2mm


Put
\begin{equation*}\label{Crystal base (L,B) of Fock space}
\begin{split}
\mathscr{L}&=\sum_{\bf{m}\in \B}\mathbb{A} \psp_{\bf{m}}|0\rangle,\ \ \ \ \
\mathscr{B}=\left\{\, \pm\psp_{\bf{m}}|0\rangle \!\!\!\! \pmod{q\mathscr{L}}\,|\,\bf{m}\in \B\, \right\},\\
\mathscr{L}^{\pm}&=\sum_{\bf{m}\in \B^\pm}\mathbb{A} \psp_{\bf{m}}|0\rangle, \ \ \
\mathscr{B}^{\pm}=\left\{\, \pm\psp_{\bf{m}}|0\rangle \!\!\!\! \pmod{q\mathscr{L}}\,|\,\bf{m}\in \B^\pm\, \right\}.
\end{split}
\end{equation*}

\begin{thm}\label{crystal base of a Fock space} The following is a crystal base of $\mathscr{V}_q$.
\begin{equation*}
\begin{cases}
(\mathscr{L},\mathscr{B}), & \text{if $\mf{g}=\mf{c}$},\\
(\mathscr{L}^+,\mathscr{B}^+), & \text{if $\mf{g}=\mf{b}, \mf{d}$},\\
(\mathscr{L}^+\otimes \mathscr{L}^+,\mathscr{B}^+\otimes \mathscr{B}^+), & \text{if $\mf{g}=\mf{b}^\bullet$}.\end{cases}
\end{equation*}
\end{thm}
\pf Suppose that $\mf{g}=\mf{c}$. It is clear by definition that $\mathscr{L}$ is an $\mathbb{A}$-lattice of $\F_q$.
We first claim that $(\mathscr{L}^\pm,\mathscr{B}^\pm)$ is a crystal base of $\F_q^\pm$ as a $U_q(\gl_{m|n})$-module, and hence $(\mathscr{L},\mathscr{B})$ is a crystal base of $\F_q$ as a $U_q(\gl_{m|n})$-module by Corollary \ref{F as gl(m|n)-module}.
Consider $(\mathscr{L}^+,\mathscr{B}^+)$. The proof for $(\mathscr{L}^-,\mathscr{B}^-)$ is the same.

Let $a\in \J_{m|n}$ with $|a|=1$ be given. Then for $m\geq 1$, $\psp_a^{(m)}$ is a highest weight vector with respect to $\langle e_a, f_a, t_a^{\pm 1} \rangle\cong U_q(\mf{sl}_2)$. We have
\begin{equation*}
\begin{split}
\tf_a \psp_a^{(m)}|0\rangle &=f_a \psp_a^{(m)}|0\rangle= {\bf f}^+_a \psp_a^{(m)}|0\rangle=\psp_{a+1}\psm_{a}\psp_a^{(m)}|0\rangle = \frac{1}{[m]}\psp_{a+1}[\om_a]\psp_a^{(m-1)}|0\rangle\\
&=\frac{1}{[m]}\psp_{a+1}\psp_a^{(m-1)}[q^{-m+1}\om_a]|0\rangle=\psp_a^{(m-1)} \psp_{a+1}|0\rangle.
\end{split}
\end{equation*}
Similarly, we have $$\tf_a^k\psp_a^{(m)}|0\rangle = f_a^{(k)} \psp_a^{(m)}|0\rangle = \psp_a^{(m-k)}\psp_{a+1}^{(k)}|0\rangle,$$ for $1\leq k\leq m$. This implies that  $\mathscr{L}^+$  is invariant under $\te_a$ and $\tf_a$ for $a\in I_{0|n}$.

On the other hand, $\mathscr{L}^+$  is invariant under $\te_a$ and $\tf_a$ for $a\in I_{m|0}$ since the weight of $\psp_{\bf m}|0\rangle$ with respect to $\langle e_a, f_a, t_a^{\pm 1} \rangle\cong U_q(\mf{sl}_2)$ is minuscule for ${\bf m}\in \B^+$ and $a=\ov{m-1},\ldots,\ov{1}$.
Also, we observe that for $m\geq 1$
\begin{equation*}
\begin{split}
e_0 \psp_{\hf}^{(m)}|0\rangle &= {\bf e}^+_0 \psp_\hf^{(m)}|0\rangle=\psp_{\ov{1}}\psm_{\hf}\psp_\hf^{(m)}|0\rangle = \psp_{\ov{1}}\psp_\hf^{(m-1)}|0\rangle,
\end{split}
\end{equation*}
which implies that $\mathscr{L}^+$ is invariant under $\te_0$ and $\tf_0$. Therefore, $\mathscr{L}^+$ is invariant under $\te_i$ and $\tf_i$ for $i\in I_{m|n}\setminus\{\ov{m}\}$. It is straightforward to check the other conditions for  $(\mathscr{L}^+,\mathscr{B}^+)$ to be  a crystal base of $\F_q^+$ as a $U_q(\gl_{m|n})$-module. This proves our claim.

By Corollary \ref{F as gl(m|n)-module}, $(\mathscr{L},\mathscr{B})$ is isomorphic to $(\mathscr{L}^-\otimes\mathscr{L}^+, \mathscr{B}^-\otimes\mathscr{B}^+)$ as a crystal base of the $U_q(\gl_{m|n})$-module $\F_q$, where $\psi_{{\bf m}^-}|0\rangle \otimes \psi_{{\bf m}^+}|0\rangle$ is mapped to $\psi_{{\bf m}^-}\psi_{{\bf m}^+}|0\rangle$ for ${\bf m}^{\pm}\in \B^{\pm}$.

Finally, for ${\bf m}\in \B$, the weight of $\psp_{\bf m}|0\rangle$ with respect to $\langle e_{\ov{m}}, f_{\ov{m}}, t_{\ov{m}}^{\pm 1} \rangle\cong U_q(\mf{sl}_2)$ is minuscule, and hence $\mathscr{L}$  is invariant under $\te_{\ov{m}}$ and $\tf_{\ov{m}}$. For example,  for ${\bf m}$ with $m_a=0$ when $a=\pm\ov{m}$
\begin{equation*}
\tf_{\ov{m}}\psi_{\bf m}|0\rangle =f_{\ov{m}}\psi_{\bf m}|0\rangle=\pm \psi_{\bf m'}|0\rangle,
\end{equation*}
where ${\bf m}'=(m'_a)\in \B$ is such that $m'_a=1$ for $a=\pm\ov{m}$  and $m'_a=m_a$ for $a\neq \pm\ov{m}$.
Therefore, we conclude that $(\mathscr{L},\mathscr{B})$ is a crystal base of $\F_q$.

We omit the proof for $\mf{g}=\mf{b}, \mf{d}$ since it is similar to the case of $\mf{g}=\mf{c}$.
Since the action of $U_q(\mf{b}_{m|n})$ on $\F^+_{q^2}\otimes\F^+_{q^2}$ is the same as $U_q(\mf{b}^\bullet_{m|n})$ except $f_{\ov{m}}$ by scalar multiplication, $(\mathscr{L}^+\otimes \mathscr{L}^+,\mathscr{B}^+\otimes \mathscr{B}^+)$ is a also crystal base of $\mathscr{V}_q$ when $\mf{g}=\mf{b}^\bullet$.
\qed

\begin{cor}
For $\ell\geq 1$, $\mathscr{V}^{\otimes \ell}_q$ has a crystal base.
\end{cor}


\section{Ortho-symplectic tableaux of type $B$ and $C$}\label{ortho-symplectic Tableaux of Type B and C}
In this section, we introduce our main combinatorial object called ortho-symplectic tableaux, which play a crucial role in the next sections.

\subsection{Semistandard tableaux}
We assume that  $\mc{A}$ is a linearly
ordered  set with a $\mathbb{Z}_2$-grading $\mc{A}=\mc{A}_0\sqcup\mc{A}_1$. When $\mc{A}=\mathbb{N}$, we assume that $\mathbb{N}_0=\mathbb{N}$ with a usual linear ordering. For a skew  Young diagram  $\lambda/\mu$, a tableau $T$ obtained by
filling $\lambda/\mu$ with entries in $\mc{A}$ is called
$\mc{A}$-semistandard  if (1) the entries in each row (resp. column) are
weakly increasing from left to right (resp. from top to bottom), (2)
the entries in $\mc{A}_0$ (resp. $\mc{A}_1$) are strictly increasing in each
column (resp. row).  We say that the shape of
$T$ is $\lambda/\mu$, and write ${\rm sh}(T)=\lambda/\mu$. The weight of $T$ is the sequence $(m_a)_{a\in\mc{A}}$, where $m_a$ is the number of occurrences of $a$ in $T$.
We denote by ${ SST}_{\mc{A}}(\lambda/\mu)$ the set of all
$\mc{A}$-semistandard tableaux of shape $\lambda/\mu$ (cf. \cite{Ful}). 

Let $x_{\mc{A}}=\{\,x_a\,|\,a\in\mc{A}\,\}$ be the set of formal commuting variables indexed by $\mc{A}$. For $\lambda\in\cP$, let $s_{\lambda}(x_{\mc{A}})$ be the super Schur function corresponding to $\lambda$, which is the weight generating function of $SST_\mc{A}(\lambda)$, that is, $s_{\lambda}(x_{\mc{A}})=\sum_{T}x_\mc{A}^T$, where $x_\mc{A}^T=\prod_{a}x_a^{m_a}$ and $(m_a)_{a\in \mc{A}}$ is the weight of $T$.

For $T\in SST_\mc{A}(\lambda)$ and
$a\in \mc{A}$, we denote by $a \rightarrow T$ the tableau obtained by applying the usual Schensted column insertion of $a$ into $T$ (cf. \cite{BR,Ful}). For a finite word $w=w_1\ldots w_r$ with letters in $\mc{A}$, we define $(w \rightarrow T)=(
w_r\rightarrow(\cdots(w_1\rightarrow T)))$. For an $\mc{A}$-semistandard tableau $S$ of skew Young diagram shape, let $w(S)$ be the word obtained by reading the entries of $S$ column by column from right to left, where in each column we read the entries from top to bottom. We define $(S\rightarrow T)=(w(S)\rightarrow T)$. We denote by $w^{\rm rev}(S)$ the reverse word of $w(S)$.

\subsection{Combinatorial $R$-matrix}
For a single-columned tableau $S$, let $S(i)$ ($i\geq 1$) denote the $i$-th entry from the bottom, and ${\rm ht}(S)$ the height of $S$.

For $k,l\in\Z_{>0}$ with $k\geq l$, let us describe a bijection \cite[Example 5.9]{KK}
\begin{equation}\label{R matrix}
R : SST_{\mc{A}}(1^k)\times SST_{\mc{A}}(1^l) \longrightarrow SST_{\mc{A}}(1^l)\times SST_{\mc{A}}(1^k),
\end{equation}
which preserves the Knuth equivalence (cf. \cite[Section 4.4]{BKK}).  It coincides with a combinatorial $R$-matrix when $\mc{A}=\mc{A}_0$ \cite{NY}.

Let $(S,T)\in SST_{\mc{A}}(1^k)\times SST_{\mc{A}}(1^l)$ be given. We choose the entries $S(i_1),\ldots, S(i_l)$ in $S$ inductively as follows: (1) If $T(1)\in \mc{A}_0$ (resp. $\mc{A}_1$), then let $S(i_1)$ be the largest one which is no greater than (resp. less than) $T(1)$. If there is no such entry, then put $i_1=1$. (2) For $2\leq u\leq l$, if $T(u)\in \mc{A}_0$ (resp. $\mc{A}_1$), then let $S(i_u)$ be the largest one in $\{\,S(1),\ldots, S(k)\,\}\setminus\{\,S(i_{1}),\ldots, S(i_{u-1})\,\}$ which is no greater than (resp. less than) $T(u)$. If there is no such entry, then choose $S(i_u)$ to be the largest one in $\{\,S(1),\ldots, S(k)\,\}\setminus\{\,S(i_{1}),\ldots, S(i_{u-1})\,\}$ at the lowest position.
Then we have $$R(S,T)=(T^\sharp,S^\sharp),$$ where $T^\sharp$ is the tableau of shape $(1^l)$ having $\{\,S(i_1),\ldots,S(i_l)\,\}$ as its entries, and $S^\sharp$ is the tableau of shape $(1^k)$ obtained by adding additional $k-l$ entries $\{\,S(1),\ldots, S(k)\,\}\setminus\{\,S(i_{1}),\ldots, S(i_{l})\,\}$ to $T$ and rearranging the entries with respect to the ordering on $\mc{A}$.\vskip 2mm

\begin{rem}{\rm In particular, if the pair $(S,T)$ forms an $\mc{A}$-semistandard tableau $U$ of shape $(k,l)'$, then $(T^\sharp \rightarrow S^{\sharp})=U$.

}
\end{rem}

\subsection{Signature of recording tableaux}\label{signature rule of a tableau}
Let us first recall a combinatorial algorithm often called 
signature rule.
Let $\sigma=(\sigma_{1},\sigma_{2},\ldots)$ be a sequence (not
necessarily of finite length) with $\sigma_{i}\in \{\,+\,,\,-\, , \ \cdot\ \}$ such that
$\sigma_i=\, +$ or $\cdot\,$ for $i\gg 0$.  We replace a pair
$(\sigma_{j},\sigma_{j'})=(+,-)$ in $\sigma$, where $j<j'$ and $\sigma_{j''}=\,\cdot\,$
for $j<j''<j'$, with $(\,\cdot\,,\,\cdot\,)$, and repeat this process
as far as possible until we get a sequence with no $-$ placed to the
right of $+$.  We denote the resulting reduced sequence by $\td{\sigma}$.

Let $w=w_1\ldots w_r$ be a finite word with letters in $\mathbb{N}$. Fix $k\in \mathbb{N}$. We associate a sequence $\sigma=(\sigma_1,\ldots,\sigma_r)$, where $\sigma_i=+$ (resp. $-$) if $w_i=k$ (resp. $k+1$) and $\sigma_i=\cdot$\ otherwise, for $1\leq i\leq r$.
We say that the {\it $k$-signature of $w$} is $(a,b)$, where $a$ (resp. $b$) is the number of $-$'s (resp. $+$'s) in $\td{\sigma}$.
For an $\mathbb{N}$-semistandard tableau $T$, we define the {\it $k$-signature of $T$} to be that of $w(T)$

Let $S_1$ and $S_2$ be two single-columned $\mc{A}$-semistandard tableaux.  Let $w=w(S_2)w(S_1)$ with $w(S_2)=w_1\ldots w_r$ and $w(S_1)=w_{r+1}\ldots w_s$. Given $U\in SST_{\mc{A}}(\mu)$ ($\mu\in \cP$), let $$P(S_1,S_2;U)=(S_1\rightarrow (S_2\rightarrow U)).$$  Suppose that $\lambda={\rm sh}(P(S_1,S_2;U))$.  Define  $$Q(S_1, S_2;U)_{[k]}$$ to be a tableau of shape $\lambda'/\mu'$  with entries in $\{k,k+1\}$ such that ${\rm sh}(w_1\ldots w_i \rightarrow U)'/{\rm sh}(w_1\ldots w_{i-1} \rightarrow U)'$ is filled with $k$  (resp. $k+1$) for $1\leq i\leq r$ (resp. $r+1\leq i\leq s$). We have $Q(S_1, S_2;U)_{[k]}\in SST_{\{k,k+1\}}(\lambda'/\mu')$ by definition. The correspondence $(S_1,S_2)\mapsto (P(S_1,S_2;U),Q(S_1, S_2;U)_{[k]})$ is reversible  by reverse bumping and hence one-to-one.

We define the {\it signature of $(S_1,S_2)$} to be the  $k$-signature of $Q(S_1, S_2;U)_{[k]}$.
Indeed, if $(a,b)$ is the $k$-signature of $Q(S_1, S_2;U)_{[k]}$, then we can check by standard arguments (cf. \cite[Chapter 1]{Ful}) that it does not depend on the choice of $U$. In particular, if  $U=\emptyset$, then $\lambda=(r+a,r-b)'$ for some $a,b$ with $a+(r-b)=s$. This means that when $S_1$ is inserted into $S_2$, $r-b$ entries in $S_2$ are bumped out and $a$ entries in $S_1$ are placed below the bottom of $S_2$. \vskip 2mm

For $a,b,c\in\Z_{\geq 0}$, let $$\lambda(a,b,c)=(2^{b+c},1^a)/(1^b)=(a+b+c,b+c)'/(b)',$$ which is a skew Young diagram with two columns of heights $a+c$ and $b+c$. For example,
$$\lambda(1,3,2)\ =\ {\def\lr#1{\multicolumn{1}{|@{\hspace{.75ex}}c@{\hspace{.75ex}}|}{\raisebox{-.04ex}{$#1$}}}\raisebox{-.6ex}
{$\begin{array}{cc}
\cline{2-2}
\cdot &\lr{}\\ 
\cline{2-2}
\cdot &\lr{}\\ 
\cline{2-2}
\cdot&\lr{}\\ 
\cline{1-1}\cline{2-2}
\lr{}&\lr{}\\ 
\cline{1-1}\cline{2-2}
\lr{\ \ }&\lr{\ \ }\\ 
\cline{1-1}\cline{2-2} 
\lr{} \\ 
\cline{1-1}
\end{array}$}}$$
In other words, $a$ (resp. $b$) denotes the length of the lower (resp. upper) tail of $\lambda(a,b,c)$, and $c$ denotes the height of its rectangular body.\vskip 2mm

Let $S$ be a tableau of shape $\lambda(a,b,c)$, whose column is $\mc{A}$-semistandard. Let $S_i$ denote the $i$-th column of $S$ from the left for $i=1,2$. 

\begin{lem}\label{signature condition}{\rm Under the above hypothesis, $S$ is $\mc{A}$-semistandard of shape $\lambda(a,b,c)$ if and only if the signature of $(S_1, S_2)$ is $(a-p,b-p)$ for some $0\leq p\leq \min \{a,b\}$.
}
\end{lem}
\pf Suppose that $S$ is $\mc{A}$-semistandard. Let $w_1\ldots w_{b+c}$ (resp. $w_{b+c+1}\ldots w_{a+b+2c}$) be the subword of $w(S)$ corresponding to the entries in $S_2$ (resp. $S_1$). Given $U\in SST_{\mc{A}}(\mu)$, let $i_s$ ($1\leq s\leq b+c$) be the row index (enumerated from the top) of $1$ in $Q(S_1, S_2;U)_{[1]}$ corresponding to $w_s$, and $j_t$ ($1\leq t\leq a+c$) the row index of $2$ in $Q(S_1, S_2;U)_{[1]}$ corresponding to $w_{t+b+c}$. Note that $i_1\geq \ldots \geq i_{b+c}$ and $j_1\geq \ldots \geq j_{a+c}$.

Since $S$ is $\mc{A}$-semistandard, we have $i_{b+u}< j_{u}$ for $1\leq u\leq c$, which implies that a cancellation of $(+,-)$ or $(1,2)$-pair occurs in $w=w(Q(S_1, S_2;U)_{[1]})$ at least as many as $c$ times. Hence the $1$-signature of $w$ is $(a-p,b-p)$ for some $0\leq p\leq \min \{a,b\}$.

Conversely, suppose that $S$ is not $\mc{A}$-semistandard. We can choose $p\geq 1$ such that there exists $\td{S}\in SST_\mc{A}(\lambda(a+p,b+p,c-p))$ whose columns are the same as those of $S$. We assume that $p$ is minimal among such ones. Then by the minimality of $p$, we can check that $j_1\leq i_{b+p}$ and hence the $1$-signature of $w$ is $(a+p,b+p)$, which is a contradiction.
\qed\vskip 2mm

\subsection{Combinatorial $R$-matrix and recording tableaux}
Given $k\in \mathbb{N}$, let us define an operator  $r_k$ on a finite word $w$ with letters in $\mathbb{N}$.  Let $w$ be given with the $k$-signature $(a,b)$, and  let $\sigma$ and  $\td{\sigma}$ denote the associated sequences with $+, - ,\,\cdot$  defined in Section \ref{signature rule of a tableau}.

If $a\leq b$, then we define $r_k w$ to be the word obtained by replacing  $k$'s corresponding to the leftmost $(b-a)$ $+$'s in $\td{\sigma}$ with $k+1$'s. Similarly, if $a\geq b$ then we define $r_k w$ to be the word obtained by replacing $k+1$ corresponding to the rightmost $(a-b)$ $-$'s in $\td{\sigma}$ with $k$. Then the $k$-signature of $r_kw$ is $(b,a)$ in both cases. Note that
$r_k$ is the Weyl group action on a tensor product of the crystal of the natural representation of $U_q(\mf{sl}_2)$, that is, $k \stackrel{k}{\rightarrow} k+1$, with respect to the tensor product rule \eqref{lower tensor product rule}.

Also, we define $\varrho_k w$ to be the word obtained by replacing $k+1$'s in $w$, which do not come from a cancelled pair $(+,-)$ in $\sigma$, with $k$'s. Note that the $k$-signature of $\varrho_k w$ is $(0, a+b)$, and $\varrho_k w=r_k w$ when $b=0$.

For an $\mathbb{N}$-semistandard tableau $T$, we define $r_k(T)$ (resp. $\varrho_k(T)$) to be the tableau obtained from $T$ by applying $r_k$ (resp. $\varrho_k$) to $w(T)$, that is, $w(r_k(T))=r_k(w(T))$ and  $w(\varrho_k(T))=\varrho_k(w(T))$.  Note that $r_k(T)$ and $\varrho_k(T)$ are well defined $\mathbb{N}$-semistandard tableaux. \vskip 2mm

Now, let $S\in SST_{\mc{A}}(\lambda(a,b,c))$ be given, and let $S_i$ be the $i$-th column  of $S$  from the left ($i=1,2$). Suppose first that $b=0$. Let
$$(T_1 ,T_2)=R(S_1, S_2)\in SST_{\mc{A}}(1^c)\times SST_{\mc{A}}(1^{a+c}),$$
where $R$ is the combinatorial $R$-matrix in  \eqref{R matrix}. The signature of $(S_1, S_2)$ is $(a,0)$ by Lemma \ref{signature condition}, and the signature of $(T_1 ,T_2)$ is $(0,a)$.
Let $U\in SST_{\mc{A}}(\lambda)$ be given for $\lambda\in \cP$. By considering the bumping paths in the insertion of $(S_1\rightarrow (S_2\rightarrow U))$ (cf. \cite[Chapter 1]{Ful}) and the definition of $R$, we have the following.

\begin{lem}\label{R matrix and signature} Under the above hypothesis,
\begin{equation*}
r_k Q(S_1, S_2;U)_{[k]}= Q(T_1 ,T_2;U)_{[k]}.
\end{equation*}
\end{lem}

Next, we suppose that $b\geq 0$ and the signature of $(S_1, S_2)$ is $(a,b)$. Let $\td{S}=(S_1\rightarrow S_2)\in SST_{\mc{A}}(\lambda(a+b,0,c))$ with $\td{S}_i$ the $i$-th column of $\td{S}$ from the left ($i=1,2$). Let
$$(T_1 ,T_2)=R(\td{S}_1, \td{S}_2)\in SST_{\mc{A}}(1^c)\times SST_{\mc{A}}(1^{a+b+c}).$$

\begin{lem}\label{R matrix and signature-2} Under the above hypothesis,
\begin{equation*}
\varrho_k Q(S_1, S_2;U)_{[k]}= Q(T_1 ,T_2;U)_{[k]}.
\end{equation*}
\end{lem}
\pf Let $S^{\rm low}_2$ be the subtableau of $S_2$ consisting of $S_2(i)$ for $1\leq i\leq c$, and let $S^{\rm up}_2$ be its complement in $S_2$. Then $(S_1,S^{\rm low}_2)$ forms an $\mc{A}$-semistandard tableau of shape $\lambda(a,0,c)$. Let $(V_1,V_2^{\rm low})=R(S_1,S^{\rm low}_2)\in SST_{\mc{A}}(1^c)\times  SST_{\mc{A}}(1^{a+c})$.

Put $U^*=(S^{\rm up}_2\rightarrow U)$ and $\mu={\rm sh}(U^*)$.
Let $Q_1$ be the tableau of shape $\mu'/\lambda'$ filled with $k$ and let $Q_2=Q(S_1,S_2^{\rm low};U^*)_{[k]}$, which is of shape $\nu'/\mu'$ with $\nu={\rm sh}(S_1\rightarrow(S_2^{\rm low}\rightarrow U^*))$.
We see by definition that the tableau $Q$ obtained by glueing $Q_1$ and $Q_2$ is equal to $Q(S_1, S_2;U)_{[k]}$.

By Lemma \ref{R matrix and signature}, we have $r_kQ(S_1,S_2^{\rm low};U^*)_{[k]}=Q(V_1,V_2^{\rm low};U^*)_{[k]}$. Also we see that each $k+1$ in $Q_2$, which is replaced by $k$ under $r_k$, is always to the east of $k$'s in $Q_1$, since the $k$-signature of $Q$ is $(a,b)$. This implies that the tableau obtained by glueing $Q_1$ and $r_kQ_2$ is equal to $\varrho_kQ$. Hence, if $V_2$ is the single-columned tableau of height $a+b+c$ obtained by placing $V^{\rm low}_2$ below $S_2^{\rm up}$, then $V_2$ is $\mc{A}$-semistandard, and $Q(V_1,V_2;U)_{[k]}=\varrho_kQ$.

Since $(V_1\rightarrow V_2)$ and $\td{S}$ are Knuth equivalent, we have  $(V_1\rightarrow V_2)=\td{S}$. Furthermore,  $V_1=T_1$ and $V_2=T_2$, since $(T_1\rightarrow T_2)=\td{S}$ and there exist unique $(Y_1,Y_2)\in SST_{\mc{A}}(1^c)\times  SST_{\mc{A}}(1^{a+b+c})$ such that $(Y_1\rightarrow Y_2)=\td{S}\in SST_{\mc{A}}(\lambda(a+b,0,c))$. Therefore, we have $\varrho_k Q(S_1, S_2;U)_{[k]}=Q(V_1,V_2;U)_{[k]}= Q(T_1 ,T_2;U)_{[k]}$.
\qed

\subsection{Ortho-symplectic tableaux of type $B$ and $C$}
Let us assume that  $\mf{g}=\mf{b}, \mf{b}^\bullet, \mf{c}$. 
From now on, for a two-columned $\mc{A}$-semistandard tableau $T$, we denote by $T^{\tt L}$ and $T^{\tt R}$ the left and right columns of $T$, respectively, and often identify $T$ with $(T^{\tt L}, T^{\tt R})$.

\begin{df}\label{Definition of T_A(a)}{\rm
For $a \in \Z_{\geq  0}$, we define ${\bf T}^{\mf{g}}_\mc{A}(a)$  to be  the set of $\mc{A}$-semistandard tableaux $T$ of shape $\lambda(a,b,c)$ such that
\begin{itemize}
\item[(1)] $(b,c)\in \{0\}\times  \Z_{\geq 0}$ if $\mf{g}=\mf{c}$,

\item[(2)] $(b,c)\in \Z_{\geq 0}\times  \Z_{\geq 0}$ and the signature of $(T^{\tt L}, T^{\tt R})$ is $(a,b)$ if $\mf{g}=\mf{b},\mf{b}^\bullet$.
\end{itemize}
We define ${\bf T}^{\rm sp}_\mc{A}$ to be the set of single-columned $\mc{A}$-semistandard tableaux.
}
\end{df}

Let $T\in {\bf T}^{\mf{g}}_\mc{A}(a)$  be given with ${\rm sh}(T)=\lambda(a,b,c)$ for some $a,b,c\geq 0$. Let us define
\begin{equation*}
({}^{\tt L}T, {}^{\tt R}T)=R(T^{\tt L}\rightarrow T^{\tt R}).
\end{equation*}
By Definition \ref{Definition of T_A(a)}, ${\rm sh}(T^{\tt L}\rightarrow T^{\tt R})=\lambda(a+b,0,c)$, which is two-columned, so that we may apply $R$ to the pair of its left and right columns and $({}^{\tt L}T, {}^{\tt R}T)\in SST_\mc{A}(1^{c})\times SST_\mc{A}(1^{a+b+c})$. Note that for $\mf{g}=\mf{c}$,
we have $({}^{\tt L}T, {}^{\tt R}T)=R(T^{\tt L},T^{\tt R})$ since $(T^{\tt L}\rightarrow T^{\tt R})=T$.
For $T\in {\bf T}^{\rm sp}_\mc{A} $, we assume that ${T}^{\tt R}={}^{\tt R}T=T$, and $T^{\tt L}={}^{\tt L}T=\emptyset$.

\begin{ex}\label{Example 1}{\rm Suppose that $\mc{A}=\J_{4|\infty}$.
$${\bf T}^{\mf{b}}_\mc{A}(2) \ni S=(S^{\tt L}, S^{\tt R}) \ = \ 
{\def\lr#1{\multicolumn{1}{|@{\hspace{.75ex}}c@{\hspace{.75ex}}|}{\raisebox{-.04ex}{$#1$}}}\raisebox{-.6ex}
{$\begin{array}{cc}
\cline{1-1}\cline{2-2}
\lr{\overline{4}}&\lr{\overline{3}}\\ 
\cline{1-1}\cline{2-2}
\lr{\overline{3}}&\lr{\overline{1}}\\ 
\cline{1-1}\cline{2-2}
\lr{\overline{1}}&\lr{\tfrac{3}{2}}\\ 
\cline{1-1}\cline{2-2} 
\lr{\tfrac{1}{2}} \\ 
\cline{1-1}
\lr{\tfrac{1}{2}} \\ 
\cline{1-1}
\end{array}$}}\ \ \ \ \ \ 
({}^{\tt L}S, {}^{\tt R}S)\  = \ 
{\def\lr#1{\multicolumn{1}{|@{\hspace{.75ex}}c@{\hspace{.75ex}}|}{\raisebox{-.04ex}{$#1$}}}\raisebox{-.6ex}
{$\begin{array}{cc}
\cline{1-1}\cline{2-2}
\lr{\overline{3}}&\lr{\overline{4}}\\ 
\cline{1-1}\cline{2-2}
\lr{\overline{1}}&\lr{\overline{3}}\\ 
\cline{1-1}\cline{2-2}
\lr{\tfrac{1}{2}}&\lr{\overline{1}}\\ 
\cline{1-1}\cline{2-2} 
&\lr{\tfrac{1}{2}} \\ 
\cline{2-2}
&\lr{\tfrac{3}{2}} \\ 
\cline{2-2}
\end{array}$}}
$$
$$
{\bf T}^{\mf{b}}_\mc{A}(3) \ni T=(T^{\tt L}, T^{\tt R}) \ = \ 
{\def\lr#1{\multicolumn{1}{|@{\hspace{.75ex}}c@{\hspace{.75ex}}|}{\raisebox{-.04ex}{$#1$}}}\raisebox{-.6ex}
{$\begin{array}{cc}
\cline{2-2}
&\lr{\overline{4}}\\ 
\cline{2-2}
&\lr{\overline{3}}\\ 
\cline{1-1}\cline{2-2}
\lr{\overline{2}}&\lr{\overline{2}}\\ 
\cline{1-1}\cline{2-2}
\lr{\overline{1}}&\lr{\overline{1}}\\ 
\cline{1-1}\cline{2-2}
\lr{\tfrac{1}{2}}&\lr{\tfrac{5}{2}}\\ 
\cline{1-1}\cline{2-2} 
\lr{\tfrac{3}{2}} \\ 
\cline{1-1}
\lr{\tfrac{3}{2}} \\ 
\cline{1-1}
\lr{\tfrac{5}{2}} \\ 
\cline{1-1}
\end{array}$}}\ \ \ \
(T^{\tt L}\rightarrow T^{\tt R})=\
{\def\lr#1{\multicolumn{1}{|@{\hspace{.75ex}}c@{\hspace{.75ex}}|}{\raisebox{-.04ex}{$#1$}}}\raisebox{-.6ex}
{$\begin{array}{cc}
\cline{1-1}\cline{2-2}
\lr{\overline{4}}&\lr{\overline{2}}\\ 
\cline{1-1}\cline{2-2}
\lr{\overline{3}}&\lr{\overline{1}}\\ 
\cline{1-1}\cline{2-2}
\lr{\overline{2}}&\lr{\tfrac{5}{2}}\\ 
\cline{1-1}\cline{2-2} 
\lr{\overline{1}}\\ 
\cline{1-1}
\lr{\tfrac{1}{2}}\\ 
\cline{1-1} 
\lr{\tfrac{3}{2}} \\ 
\cline{1-1}
\lr{\tfrac{3}{2}} \\ 
\cline{1-1}
\lr{\tfrac{5}{2}} \\ 
\cline{1-1}
\end{array}$}}\ \ \ 
({}^{\tt L}T, {}^{\tt R}T)\  = \ 
{\def\lr#1{\multicolumn{1}{|@{\hspace{.75ex}}c@{\hspace{.75ex}}|}{\raisebox{-.04ex}{$#1$}}}\raisebox{-.6ex}
{$\begin{array}{cc}
\cline{2-2}
&\lr{\overline{4}}\\ 
\cline{2-2}
&\lr{\overline{3}}\\ 
\cline{1-1}\cline{2-2}
\lr{\overline{2}}&\lr{\overline{2}}\\ 
\cline{1-1}\cline{2-2}
\lr{\overline{1}}&\lr{\overline{1}}\\ 
\cline{1-1}\cline{2-2}
\lr{\tfrac{3}{2}}&\lr{\tfrac{1}{2}}\\ 
\cline{1-1}\cline{2-2} 
&\lr{\tfrac{3}{2}} \\ 
\cline{2-2}
&\lr{\tfrac{5}{2}} \\ 
\cline{2-2}
&\lr{\tfrac{5}{2}} \\ 
\cline{2-2}
\end{array}$}}
$$

}
\end{ex}

\begin{df}\label{admissible}{\rm
Suppose that $S \in {\bf T}^{\mf{g}}_\mc{A}(p) \cup {{\bf T}^{\rm sp}_\mc{A}}$, and $T \in {\bf T}^{\mf{g}}_\mc{A}(q)$   are given for $p\leq q$. We assume that $p=0$ when $T\in {{\bf T}^{\rm sp}_\mc{A}}$. We say that the pair $(S,T)$ is {\it admissible} and write $S\prec T$ if it satisfies the following conditions:

\begin{itemize}
\item[(1)] ${\rm ht}(S^{\tt R})+p\leq {\rm ht}(T^{\tt L})$,

\item[(2)]  $^{\tt R}S(i)\leq {}T^{\tt L}(i)$ for $i\geq 1$,

\item[(3)] $S^{\tt R}(i+q-p)\leq{}^{\tt L}T(i)$ for $i\geq 1$,
\end{itemize}
where the equality holds in (2) and (3) only if the entries are even or in $\mc{A}_0$.}
\end{df}

\begin{ex}{\rm Let $S$ and $T$ be as in Example \ref{Example 1}. Then $(S,T)$ is admissible or  $S\prec T$ since
$$
({}^{\tt R}S, T^{\tt L}) =\ 
{\def\lr#1{\multicolumn{1}{|@{\hspace{.75ex}}c@{\hspace{.75ex}}|}{\raisebox{-.04ex}{$#1$}}}\raisebox{-.6ex}
{$\begin{array}{cc}
\cline{2-2}
&\lr{\overline{2}}\\ 
\cline{1-1}\cline{2-2}
\lr{\overline{4}}&\lr{\overline{1}}\\ 
\cline{1-1}\cline{2-2}
\lr{\overline{3}}&\lr{\tfrac{1}{2}}\\ 
\cline{1-1}\cline{2-2}
\lr{\overline{1}}&\lr{\tfrac{3}{2}}\\ 
\cline{1-1}\cline{2-2} 
\lr{\tfrac{1}{2}}&\lr{\tfrac{3}{2}} \\ 
\cline{1-1}\cline{2-2}
\lr{\tfrac{3}{2}}&\lr{\tfrac{5}{2}} \\ 
\cline{1-1}\cline{2-2}
\end{array}$}}\ \ \ \ \ \ \ \ \ \ 
({S}^{\tt R} ,{}^{\tt L}T) =\ 
{\def\lr#1{\multicolumn{1}{|@{\hspace{.75ex}}c@{\hspace{.75ex}}|}{\raisebox{-.04ex}{$#1$}}}\raisebox{-.6ex}
{$\begin{array}{cc}
\cline{2-2}
&\lr{\overline{2}}\\ 
\cline{1-1}\cline{2-2}
\lr{\overline{3}}&\lr{\overline{1}}\\ 
\cline{1-1}\cline{2-2}
\lr{\overline{1}}&\lr{\tfrac{3}{2}}\\ 
\cline{1-1}\cline{2-2}
\lr{\tfrac{3}{2}}& \cdot\\ 
\cline{1-1}
\cdot & \cdot \\
\cdot & \cdot 
\end{array}$}}
$$
}
\end{ex}

\begin{rem}{\rm The conditions (2) and (3) in Definition \ref{admissible} are equivalent to saying that $(^{\tt R}S,T^{\tt L})$ and $(S^{\tt R},{}^{\tt L}T)$ form    $\mc{A}$-semistandard tableaux of shape $\lambda(a,b,c)$ and $\lambda(a^*,b^*,c^*)$, respectively, where
\begin{equation*}
\begin{cases}
a= 0, \\
b= {\rm ht}(T^{\tt L})-{\rm ht}(S^{\tt L}),\\
c= {\rm ht}(S^{\tt L}),
\end{cases} \ \ \
\begin{cases}
a^*= q-p, \\
b^*={\rm ht}(T^{\tt L})-{\rm ht}(S^{\tt R})-p,\\
c^*= {\rm ht}(S^{\tt R})+p -q.
\end{cases}
\end{equation*}

}
\end{rem}





Now we are in a position to introduce the notion of ortho-symplectic tableaux of type ${\mf g}$, which is our main combinatorial object.

\begin{df}{\rm
Let $(\lambda,\ell)\in\cP({\mf g})$ be given. Let
\begin{equation}\label{Length of tuples}
L=
\begin{cases}
\ell, & \text{if $\mf{g}=\mf{c}$},\\
\ell/2, & \text{if $\mf{g}=\mf{b}^\bullet$ or $\mf{g}=\mf{b}$ with $\ell-2\lambda_1$ even},\\
(\ell+1)/2, & \text{if $\mf{g}=\mf{b}$ with $\ell-2\lambda_1$ odd}.
\end{cases}
\end{equation}
We define ${\bf T}^{\mf{g}}_\mc{A}(\lambda,\ell)$ to be the set of ${\bf T}=(T_L,\ldots,T_1)$  in
\begin{equation*}
\begin{cases}
{\bf T}^{\rm sp}_\mc{A}\times{\bf T}^{\mf{b}}_\mc{A}(\lambda'_{L-1})\times\cdots \times {\bf T}^{\mf{b}}_\mc{A}(\lambda'_1), & \text{if $\mf{g}=\mf{b}$ with $\ell-2\lambda_1$ odd},\\
{\bf T}^{\mf{g}}_\mc{A}(\lambda'_{L})\times\cdots \times {\bf T}^{\mf{g}}_\mc{A}(\lambda'_1), & \text{otherwise},\\
\end{cases}
\end{equation*}
such that $T_{k+1}\prec T_{k}$  for $1\leq k\leq L-1$. We call ${\bf T}\in {\bf T}^{\mf{g}}_\mc{A}(\lambda,\ell)$ an {\it ortho-symplectic tableau of type ${\mf g}$ and shape $(\lambda,\ell)$}.
}
\end{df}

Note that $\cP(\mf{b}^\bullet)\varsubsetneq \cP(\mf{b})$ and ${\bf T}^{\mf{b}^\bullet}_\mc{A}(\lambda,\ell)={\bf T}^{\mf{b}}_\mc{A}(\lambda,\ell)$ as a set for $(\lambda,\ell)\in \cP(\mf{b}^\bullet)$.

\subsection{Schur positivity}
Let $z$ be an indeterminate. 
For $(\lambda,\ell)\in\cP({\mf g})$, put $$S^{\mf g}_{(\lambda,\ell)}(x_{\mc{A}})=z^\ell\sum_{{\bf T}\in {\bf T}^{\mf g}_\mc{A}(\lambda,\ell)}\prod_{k=1}^L x_\mc{A}^{T_k},$$
which is the weight generating function of ${\bf T}^{\mf{g}}_\mc{A}(\lambda,\ell)$.
We will show that $S^{\mf g}_{(\lambda,\ell)}(x_{\mc{A}})$ can be written as a non-negative integral (possibly infinite) linear combination  of $s_\mu(x_\mc{A})$. For this, we introduce the following.

\begin{df}\label{ortho-symplectic Kostka}{\rm
Let $(\lambda,\ell)\in\cP(\mf{g})$ be given with $L$ as in \eqref{Length of tuples}.  For $\mu\in\cP$, we
define ${\bf K}^{\mf{g}}_{\mu\,(\lambda,\ell)}$ to be the set of $Q\in SST_{\{1,\ldots,2L\}}(\mu)$ with weight $(m_1,\ldots, m_{2L})$ satisfying the following conditions:\vskip 2mm


$\bullet$ $\mf{g}=\mf{c}$ 
\begin{itemize}
\item[(1)] $m_{2k}-m_{2k-1}=\lambda'_k$ for $1\leq k\leq L$,

\item[(2)] $m_{2k}\geq m_{2k+2}$ for $1\leq k\leq L-1$,

\item[(3)] the $(2k-1)$-signature of $Q$ is $(\lambda'_k,0)$ for $1\leq k\leq L$,

\item[(4)] the $2k$-signature of $r_{2k+1}Q$ is $(0,m_{2k}-m_{2k+2})$ for $1\leq k\leq L-1$,

\item[(5)] the $2k$-signature of $r_{2k-1}Q$ is $(\lambda'_k-\lambda'_{k+1}-p,m_{2k}-m_{2k+2}-p)$ with $p\geq 0$ for $1\leq k\leq L-1$.
\end{itemize}\vskip 2mm

$\bullet$ $\mf{g}=\mf{b}, \mf{b}^\bullet$
\begin{itemize}
\item[(1)] $m_{2L}=0$ if $\ell-2\lambda_1$ is odd,

\item[(2)] $m_{2k-1}\geq m_{2k}-\lambda'_k\geq 0$ for $1\leq k\leq L$,

\item[(3)] $m_{2k}\geq m_{2k+1}+\lambda'_{k+1}$ for $1\leq k\leq L-1$,


\item[(4)] the $(2k-1)$-signature of $Q$ is $(\lambda'_k,m_{2k-1}-m_{2k}+\lambda'_k)$ for $1\leq k\leq L$,

\item[(5)] the $2k$-signature of $\varrho_{2k+1}Q$ is $(0,m_{2k}-m_{2k+1}-\lambda'_{k+1})$ for $1\leq k\leq L-1$,

\item[(6)] the $2k$-signature of $\varrho_{2k-1}Q$ is $(\lambda'_k-\lambda'_{k+1}-p,m_{2k}-m_{2k+1}-\lambda'_{k+1}-p)$ with $p\geq 0$ for $1\leq k\leq L-1$.

\end{itemize}\vskip 2mm

}
\end{df}

Then we have the following.

\begin{thm}\label{Schur positivity} For $(\lambda,\ell)\in\cP({\mf g})$, we have a weight preserving bijection
\begin{equation*}
\psi_{(\lambda,\ell)} : {\bf T}^{\mf g}_\mc{A}(\lambda,\ell) \longrightarrow \bigsqcup_{\mu\in\cP}  SST_{\mc{A}}(\mu)\times  {\bf K}^{\mf g}_{\mu\,(\lambda,\ell)}.
\end{equation*}
\end{thm}
\pf Let us first prove the case when $\mf{g}=\mf{c}$.

Let ${\bf T}=(T_\ell,\ldots,T_1)\in {\bf T}^{\mf c}_\mc{A}(\lambda'_\ell)\times\cdots \times {\bf T}^{\mf c}_\mc{A}(\lambda'_1)
$ be given. 
Put $P_1=T^{\tt R}_1$ and $P_2=(T^{\tt L}_1\rightarrow P_1)$, and  define  inductively
$$P_{2k-1}=(T^{\tt R}_k\rightarrow P_{2k-2}),\ \ \ P_{2k}=(T^{\tt L}_k\rightarrow P_{2k-1}),$$ for $2\leq k\leq \ell$.  Let $P=P_{2\ell}$. Suppose that $\mu={\rm sh}(P)$. Define $Q$ to be a tableau in $SST_{\{1,\ldots, 2\ell\}}(\mu')$ such that ${\rm sh}(P_{k})'/{\rm sh}(P_{k-1})'$ is filled with $k$ for $1\leq k\leq 2\ell$, where we assume that $P_0=\emptyset$. Note that the subtableau of $Q$ with entries $\{\,2k-1,2k\,\}$ (resp. $\{\,2k,2k+1\,\}$) is $Q(T^{\tt L}_{k},T^{\tt R}_k;P_{2k-2})_{[2k-1]}$ (resp.  $Q(T^{\tt R}_{k+1},T^{\tt L}_k;P_{2k-1})_{[2k]}$) for $1\leq k\leq \ell$ (resp. $1\leq k\leq \ell-1$).

Let $(m_1,\ldots, m_{2\ell})$ be the weight of $Q$.
Then $m_k$'s satisfy the conditions (1) and (2) for ${\bf K}^{\mf{c}}_{\mu\,(\lambda,\ell)}$ in Definition \ref{ortho-symplectic Kostka} since $m_{2k-1}={\rm ht}(T^{\tt R}_k)$ and $m_{2k}={\rm ht}(T^{\tt L}_k)$ for $1\leq k\leq \ell$.
Also, by Lemma \ref{signature condition}, $Q$ satisfies the condition (3) since each $T_k$ is $\mc{A}$-semistandard of shape $\lambda(m_{2k}-m_{2k-1},0,m_{2k-1})$, where $m_{2k}-m_{2k-1}=\lambda'_k$.

For $1\leq k\leq \ell$, let us define a sequence of tableaux $P^{(k)}_i$  ($1\leq i\leq 2\ell$)  in the same way as $P_i$'s except $$P^{(k)}_{2k-1}=({}^{\tt R}T_k\rightarrow P_{2k-2}^{(k)}),\ \ \ P^{(k)}_{2k}=({^{\tt L}}T_k\rightarrow P_{2k-1}^{(k)}).$$ Since $(T^{\tt L}_k\rightarrow T^{\tt R}_k)=({}^{\tt L}T_k\rightarrow {}^{\tt R}T_k)=T_k$, we have $P^{(k)}_i=P_i$ for $i\neq 2k-1, 2k$. We define $Q^{(k)}$ in a similar way using $P^{(k)}_i$'s.  By Lemma \ref{R matrix and signature}, we have $Q^{(k)}=r_{2k-1}Q$. Then by Lemma \ref{signature condition}, the conditions (2) and (3) in Definition \ref{admissible} imply the conditions (4) and (5)  in Definition \ref{ortho-symplectic Kostka}, respectively. Hence, we have $Q\in {\bf K}^{\mf c}_{\mu\,(\lambda,\ell)}$.

Now, we define $\psi_{(\lambda,\ell)}({\bf T})=(P,Q)$.  Since the construction of $(P,Q)$ is reversible, it is not difficult to see that $\psi_{(\lambda,\ell)}$ is a bijection to the set of pairs  $(P,Q)$ with $P\in SST_\mc{A}(\mu)$ and $Q\in {\bf K}^{\mf c}_{\mu\,(\lambda,\ell)}$ for $\mu\in\cP$.

The proof for the cases when $\mf{g}={\mf b}, {\mf b}^\bullet$ is almost identical, where we replace $r_k$ with $\varrho_k$ and use Lemma \ref{R matrix and signature-2} instead of Lemma \ref{R matrix and signature}. We leave the details to the reader. \qed

\begin{cor}\label{Schur positivity-2}
 For $(\lambda,\ell)\in\cP({\mf g})$, we have
\begin{equation*}
S^{\mf g}_{(\lambda,\ell)}(x_{\mc{A}})=z^{\ell}\sum_{\mu\in\cP}K^{\mf{g}}_{\mu\,(\lambda,\ell)}s_\mu(x_\mc{A}),
\end{equation*}
where $K^{\mf{g}}_{\mu\,(\lambda,\ell)}$ is the number of tableaux in ${\bf K}^{\mf{g}}_{\mu\,(\lambda,\ell)}$.
\end{cor}

\section{Character of a highest weight  module in $\mc{O}^{int}(m|n)$}\label{crystal structure for ortho-symplectic tableaux m+n}
In this section, we show that the weight generating function of ortho-symplectic tableaux of shape $(\lambda,\ell)$ is the character of $L(\mf{g}_{m|n},\Lambda_{m|n}(\lambda,\ell))$ for $(\lambda,\ell)\in \cP(\mf{g})_{m|n}$, when $\mf{g}=\mf{b}, \mf{b}^\bullet, \mf{c}$.

\subsection{Combinatorial character formula for irreducibles in $\mc{O}^{int}(m|n)$}
We assume that  $\mf{g}=\mf{b}, \mf{b}^\bullet, \mf{c}$. Let 
$P_{m+n}=\bigoplus_{a\in \J_{m+n}}\Z \delta_a \oplus \Z \Lambda_{\ov{m}}$, and let
$\cP(\mf{g})_{m+n}$ be the set of $(\lambda,\ell)\in \cP(\mf{g})$ such that $\Lambda_{m+\infty}(\lambda,\ell)\in P_{m+n}$. Write $\Lambda_{m+n}(\lambda,\ell)=\Lambda_{m+\infty}(\lambda,\ell)$ for $\Lambda_{m+\infty}(\lambda,\ell)\in P_{m+n}$. Let $U_q({\mf g}_{m+n})$ be the quantized enveloping algebra associated to ${\mf g}_{m+n}$, and $L_q({\mf g}_{m+n},\Lambda)$ its irreducible highest weight module with highest weight $\Lambda\in P_{m+n}$.

 For simplicity, we put
\begin{equation*}
\begin{split}
{\bf T}_{m+n}(a)={\bf T}^{\mf{g}}_{\J_{m+n}}(a),\ \
{\bf T}_{m+n}(\lambda,\ell)={\bf T}^{\mf{g}}_{\J_{m+n}}(\lambda,\ell),\ \
{\bf T}_{m+n}^{\rm sp}={\bf T}^{\rm sp}_{\J_{m+n}},
\end{split}
\end{equation*}
for $a\in\Z_{\geq 0}$ and $(\lambda,\ell)\in \cP(\mf{g})_{m+n}$.
We first define an (abstract)  ${\mf g}_{m+n}$-crystal structure on ${\bf T}_{m+n}(\lambda,\ell)$ for $(\lambda,\ell)\in\cP(\mf{g})_{m+n}$, and show that it is isomorphic to the crystal of $L_q(\mf{g}_{m+n},\Lambda_{m+n}(\lambda,\ell))$.

We refer the reader to \cite{HK,Kas94} and references therein for general exposition on (abstract) crystals associated to a symmetrizable Kac-Moody algebra. Note that $\mf{b}^\bullet_{m+n}$ is a Kac-Moody superalgebra, and one can consider abstract $\mf{b}^\bullet_{m+n}$-crystals in the same way (cf. \cite{Je}). To avoid confusion with those for $U_q(\mf{g}_{m|n})$, we denote the Kashiwara operators on ${\mf g}_{m+n}$-crystals  by $\td{\mathsf e}_i$ and $\td{\mathsf f}_i$ for $i\in I_{m+n}$, and assume that the tensor product rule in this case follows \eqref{lower tensor product rule}.

Recall that $\J_{m+n}$ can be identified with the crystal of the natural representation of $U_q(\gl_{m+n})$ with respect to $\td{\mathsf e}_i$ and $\td{\mathsf f}_i$ for $i\in I_{m+n}\setminus\{\ov{m}\}$, where
\begin{equation*}
\begin{split}
&\ov{m}\ \stackrel{^{\ov{m-1}}}{\longrightarrow}\ \ov{m-1}\ \stackrel{^{\ov{m-2}}}{\longrightarrow}
\cdots\stackrel{^{\ov{1}}}{\longrightarrow}\ \ov{1}\ \stackrel{^{0}}{\longrightarrow}\ 1 \
\stackrel{^1}{\longrightarrow}\ 2\stackrel{^2}{\longrightarrow}\cdots\ \ \ \ \ \ \ \ \ \ \ (n=\infty), \\
&\ov{m}\ \stackrel{^{\ov{m-1}}}{\longrightarrow}\ \ov{m-1}\ \stackrel{^{\ov{m-2}}}{\longrightarrow}
\cdots\stackrel{^{\ov{1}}}{\longrightarrow}\ \ov{1}\ \stackrel{^{0}}{\longrightarrow}\ 1 \
\stackrel{^1}{\longrightarrow}\ 2\stackrel{^2}{\longrightarrow}\cdots \stackrel{^{n-1}}{\longrightarrow} n \, \ \ \  (n<\infty),
\end{split}
\end{equation*}
with ${\rm wt}(a)=\delta_a$ for $a\in \J_{m+n}$. Then the set of finite words with letters in $\J_{m+n}$ has a $\gl_{m+n}$-crystal structure, where each non-empty word $w=w_1\cdots w_r$ is identified with $w_1\otimes\cdots\otimes w_r\in (\J_{m+n})^{\otimes r}$. Also, by applying $\td{\mathsf e}_i$ and $\td{\mathsf f}_i$ to the word of an $\J_{m+n}$-semistandard tableau, we have an (abstract) $\gl_{m+n}$-crystal structure on $SST_{\J_{m+n}}(\lambda/\mu)$ for a skew Young diagram $\lambda/\mu$ \cite{KashNaka}. For $\lambda\in\cP$, we denote by $H_{\lambda}$ the highest weight element in $SST_{\J_{m+n}}(\lambda)$.
\vskip 2mm

Let ${\mc B}$ denote either ${\bf T}_{m+n}^{\rm sp}$ or ${\bf T}_{m+n}(a)$ for $0\leq a\leq m+n$. Let us define a ${\mf g}_{m+n}$-crystal structure on ${\mc B}$. 

For $T\in {\mc B}$, let
\begin{equation*}
{\rm wt}(T)=
\begin{cases}
r\Lambda_{\ov{m}}+\sum_{s\in \J_{m+n}}m_s\delta_s, & \text{if  $\mc{B}={\bf T}_{m+n}(a)$ for $0\leq a\leq m+n$},\\
\Lambda_{\ov{m}}+\sum_{s\in \J_{m+n}}m_s\delta_s, & \text{if  $\mc{B}={\bf T}_{m+n}^{\rm sp}$},\\
\end{cases}
\end{equation*}
where $r=1$ for $\mf{g}=\mf{c}$ and $r=2$ otherwise, and $(m_s)_{s\in \J_{m+n}}$ is the weight of $T$. 
Note that ${\mc B}$ is a set of $\J_{m+n}$-semistandard tableaux, and the signature of $T\in {\bf T}_{m+n}(a)$ is invariant under $\td{\mathsf e}_i$ and $\td{\mathsf f}_i$ for $i\in I_{m|n}\setminus\{0\}$ such that $\td{\mathsf e}_i T\neq {\bf 0}$ or $\td{\mathsf f}_i T\neq {\bf 0}$, since ${\rm sh}(T^{\tt L}\rightarrow T^{\tt R})$ is invariant under  $\td{\mathsf e}_{i}$, $\td{\mathsf f}_{i}$.
Hence ${\mc B}$ is a $\gl_{m+n}$-crystal with respect to $\td{\mathsf e}_i$ and $\td{\mathsf f}_i$ for $i\in I_{m+n}\setminus\{\ov{m}\}$, where $\varepsilon_i$ and $\varphi_i$ are defined in a usual way. So it suffices to define $\td{\mathsf e}_{\ov{m}}$ and $\td{\mathsf f}_{\ov{m}}$.\vskip 2mm

\textsc{Case 1}. Suppose that $\mf{g}=\mf{c}$ and $\mc{B}={\bf T}_{m+n}(a)$ for $0\leq a\leq m+n$.
For $T\in {\bf T}_{m+n}(a)$, let $t^{\tt L}$ and $t^{\tt R}$ be the top entries in $T^{\tt L}$ and $T^{\tt R}$, respectively. If $t^{\tt L}=t^{\tt R}=\ov{m}$, then we define $\td{\mathsf e}_{\ov{m}}T$ to be the tableau obtained  by removing
$\resizebox{.07\hsize}{!}{${\def\lr#1{\multicolumn{1}{|@{\hspace{.6ex}}c@{\hspace{.6ex}}|}{\raisebox{-.2ex}{$#1$}}}\raisebox{-.8ex}
{$\begin{array}[b]{cc}
\cline{1-1}\cline{2-2}
\lr{\ov{m}}&\lr{\ov{m}}\\
\cline{1-1}\cline{2-2}
\end{array}$}}$}$ from $T$. Otherwise, we define $\td{\mathsf e}_{\ov{m}}T={\bf 0}$. Also, if $t^{\tt L}, t^{\tt R}>\ov{m}$, then we define $\td{\mathsf f}_{\ov{m}}T$ to be the tableau obtained by adding
$\resizebox{.07\hsize}{!}{${\def\lr#1{\multicolumn{1}{|@{\hspace{.6ex}}c@{\hspace{.6ex}}|}{\raisebox{-.2ex}{$#1$}}}\raisebox{-.8ex}
{$\begin{array}[b]{cc}
\cline{1-1}\cline{2-2}
\lr{\ov{m}}&\lr{\ov{m}}\\
\cline{1-1}\cline{2-2}
\end{array}$}}$}$  on top of $T$. Otherwise, we define $\td{\mathsf f}_{\ov{m}}T={\bf 0}$. Here ${\bf 0}$ denotes a formal symbol conventionally used in abstract crystals. 

\textsc{Case 2}. Suppose that $\mf{g}=\mf{b}$ and $\mc{B}={\bf T}^{\rm sp}_{m+n}$.
For $T\in {\bf T}^{\rm sp}_{m+n}$, let $t$ be the top entry of $T$. If $t=\ov{m}$, then we define $\td{\mathsf e}_{\ov{m}}T$ to be the tableau  obtained by removing $\resizebox{.035\hsize}{!}{$\boxed{\ov{m}}$}$ from $T$. Otherwise, we define $\td{\mathsf e}_{\ov{m}}T={\bf 0}$. We define $\td{\mathsf f}_{\ov{m}}T$ in a similar way by adding $\resizebox{.035\hsize}{!}{$\boxed{\ov{m}}$}$ on top of $T$. 

For $T\in \mc{B}$ in the above two cases, we put
$\varepsilon_{\ov{m}}(T)=\max\{\,r\in\mathbb{Z}_{\geq 0}\,|\,\td{\mathsf e}_{\ov{m}}^r T\neq 0\,\}$ and
$\varphi_{\ov{m}}(T)={\rm wt}(T)+\varepsilon_{\ov{m}}(T)$.
Then we can check that ${\mc B}\cup\{{\bf 0}\}$ is invariant under $\td{\mathsf e}_{\ov{m}}$ and $\td{\mathsf f}_{\ov{m}}$ (we assume $\td{\mathsf e}_{i}{\bf 0}=\td{\mathsf f}_{i}{\bf 0}={\bf 0}$), and hence ${\mc B}$ is a ${\mf g}_{m+n}$-crystal with respect to ${\rm wt}$, $\varepsilon_i$, $\varphi_i$ and $\td{\mathsf e}_{i}$, $\td{\mathsf f}_{i}$ for $i\in I_{m+n}$.

\textsc{Case 3}. Suppose that $\mf{g}=\mf{b}, \mf{b}^\bullet$ and $\mc{B}={\bf T}_{m+n}(a)$ for $0\leq a\leq m+n$.
We regard ${\bf T}_{m+n}(a)$ as a subset of $({\bf T}^{\rm sp}_{m+n})^{\otimes 2}$ by identifying $T\in {\bf T}_{m+n}(a)$ with $T^{\tt R}\otimes T^{\tt L}$. 
Suppose that $T\in {\bf T}_{m+n}(a)$ is given with ${\rm sh}(T)=\lambda(a,b,c)$ and $\td{\mathsf e}_{\ov{m}}T\neq {\bf 0}$. If $\td{\mathsf e}_{\ov{m}}T= T^{\tt R}\otimes (\td{\mathsf e}_{\ov{m}}T^{\tt L})$, then  ${\rm sh}(\td{\mathsf e}_{\ov{m}}T)=\lambda(a,b+1,c-1)$ and the signature of $\td{\mathsf e}_{\ov{m}}T$ is $(a,b+1)$ by Lemma \ref{signature condition}. Similarly, if $\td{\mathsf e}_{\ov{m}}T=(\td{\mathsf e}_{\ov{m}} T^{\tt R}) \otimes T^{\tt L}$, then  ${\rm sh}(\td{\mathsf e}_{\ov{m}}T)=\lambda(a,b-1,c)$ with $b>0$ and the signature of $\td{\mathsf e}_{\ov{m}}T$ is $(a,b-1)$. Hence ${\bf T}_{m+n}(a)\cup\{{\bf 0}\}$ is invariant under $\td{\mathsf e}_{\ov{m}}$.
By almost the same argument, we can check that ${\bf T}_{m+n}(a)\cup\{{\bf 0}\}$ is also invariant under $\td{\mathsf f}_{\ov{m}}$.
Therefore, ${\bf T}_{m+n}(a)$ is a subcrystal of $({\bf T}^{\rm sp}_{m+n})^{\otimes 2}$ with respect to ${\rm wt}$, $\varepsilon_i$, $\varphi_i$ and $\td{\mathsf e}_{i}$, $\td{\mathsf f}_{i}$ for $i\in I_{m+n}$.\vskip 2mm


\begin{thm}\label{regularity of level 1 osp crystal}\mbox{}
\begin{itemize}
\item[(1)] ${\bf T}_{m+n}^{\rm sp}$ is isomorphic to the crystal of $L_q({\mf b}_{m+n},\Lambda_{\ov{m}})$.

\item[(2)] ${\bf T}_{m+n}(a)$ is isomorphic to the crystal of $L_q({\mf g}_{m+n},\Lambda_{m+n}((1^a),r))$ for $0\leq a\leq m+n$, where $r=1$ for $\mf{g}=\mf{c}$ and $r=2$ otherwise.
\end{itemize}
\end{thm}
\pf (1) Let $T\in {\bf T}^{\rm sp}_{m+n}$ be given. Let $(\sigma_a)_{a\in \J_{m+n}}$ be the sequence of $\pm$ such that $\sigma_a = -$ if and only if $a$ occurs as an entry of $T$. Then the map sending $T$ to $(\sigma_a)$ is isomorphism of $\mf{b}_{m+n}$-crystals from ${\bf T}^{\rm sp}_{m+n}$ to the crystal of the spin representation $L_q({\mf b}_{m+n},\Lambda_{\ov{m}})$ (cf. \cite[Section 5.4]{KashNaka}).

(2) Suppose that $\mf{g}=\mf{c}$. We first claim  that ${\bf T}_{m+n}(a)$ is connected. We use induction on the number of boxes in $T\in {\bf T}_{m+n}(a)$, say $|T|$, to show that $T$ is connected to $H_{(1^a)}\in SST_{\J_{m+n}}(1^a)\subset {\bf T}_{m+n}(a)$. Suppose that $T\in {\bf T}_{m+n}(a)$ is given with ${\rm sh}(T)=\mu$. Since ${\bf T}_{m+n}(a)$ is a $\gl_{m+n}$-crystal, $T$ is connected to a highest weight element in $SST_{\J_{m+n}}(\mu)$ whose columns have $\ov{m}$ as top entries. Hence $\td{\mathsf e}_{\ov{m}}T\neq {\bf 0}$ and $|\td{\mathsf e}_{\ov{m}}T|=|T|-2$, and by induction hypothesis, $T$ is connected to $H_{(1^a)}$.

Next we claim that ${\bf T}_{m+n}(a)$ is a regular ${\mf c}_{m+n}$-crystal, that is, for $J\subset I_{m+n}$ with $|J|\leq 2$ such that $(\langle \alpha_i^\vee,\alpha_j \rangle)_{i,j\in J}$ is of positive definite, it is isomorphic to a crystal of an integrable representation of type $(\langle \alpha_i^\vee,\alpha_j \rangle)_{i,j\in J}$. Since  ${\bf T}_{m+n}(a)$ is a regular $\gl_{m+n}$-crystal, it remains to consider the case when $J=\{\,\ov{m},\ov{m-1}\,\}$. Let $C(T)$ be the connected component of $T\in {\bf T}_{m+n}(a)$ with respect to $\td{\mathsf{e}}_{i}, \td{\mathsf{f}}_{i}$ for $i=\ov{m},\ov{m-1}$. It is straightforward to see that $C(T)$ is isomorphic to the crystal of a fundamental representation of type $C_2$ \cite{KashNaka}.  Hence ${\bf T}_{m+n}(a)$ is a regular ${\mf c}_{m+n}$-crystal and it is isomorphic to the crystal of an integrable $U_q({\mf c}_{m+n})$-module \cite{KKMMNN92}.
Since ${\bf T}_{m+n}(a)$ is a connected  crystal with highest weight ${\rm wt}(H_{(1^a)})=\Lambda_{m+n}((1^a),1)$, it is isomorphic to the crystal of $L_q({\mf c}_{m+n},\Lambda_{m+n}((1^a),1))$.

Suppose that $\mf{g}=\mf{b}$. By similar arguments as in $\mf{g}=\mf{c}$, we can show that ${\bf T}_{m+n}(a)$ is connected for $0\leq a\leq m+n$ with the highest weight element $H_{(1^a)}$. Since ${\bf T}_{m+n}(a)$ is a subcrystal of $({\bf T}^{\rm sp}_{m+n})^{\otimes 2}$, which is a regular ${\mf b}_{m+n}$-crystal by (1), ${\bf T}_{m+n}(a)$ is also regular with highest weight $\Lambda_{m+n}((1^a),2)$. Therefore, ${\bf T}_{m+n}(a)$ is isomorphic to the crystal of $L_q({\mf b}_{m+n},\Lambda_{m+n}((1^a),2))$.

As an $I_{m+n}$-colored oriented graph, the crystal of $L_q({\mf b}^\bullet_{m+n},\Lambda_{m+n}((1^a),2))$ is isomorphic to that  of $L_q({\mf b}_{m+n},\Lambda_{m+n}((1^a),2))$ (see \cite{BK}). This implies that ${\bf T}_{m+n}(a)$ is also  isomorphic to the crystal of $L_q({\mf b}^\bullet_{m+n},\Lambda_{m+n}((1^a),2))$. \qed\vskip 2mm

Let $(\lambda,\ell)\in \cP({\mf g})_{m+n}$ be given with $L$ as in \eqref{Length of tuples}. We consider ${\bf T}_{m+n}(\lambda,\ell)$ as a subset of
\begin{equation}\label{tensor product of fundamental crystals}
\begin{cases}
{\bf T}_{m+n}(\lambda'_1)\otimes\cdots \otimes {\bf T}_{m+n}(\lambda'_{L-1})\otimes {\bf T}^{\rm sp}_{m+n}, & \text{if $\mf{g}=\mf{b}$ with $\ell-2\lambda_1$ odd},\\
{\bf T}_{m+n}(\lambda'_1)\otimes\cdots \otimes {\bf T}_{m+n}(\lambda'_L), & \text{otherwise},\\
\end{cases}\\
\end{equation}
by identifying ${\bf T}=(T_L,\ldots, T_1) \in {\bf T}_{m+n}(\lambda,\ell)$ with $T_1\otimes \cdots \otimes T_L$, and apply $\td{\mathsf{e}}_{i}, \td{\mathsf{f}}_{i}$ on ${\bf T}_{m+n}(\lambda,\ell)$ for $i\in I_{m+n}$.

\begin{lem}\label{crystal invariance of osp tableaux}
${\bf T}_{m+n}(\lambda,\ell)\cup\{{\bf 0}\}$ is invariant under $\td{\mathsf{e}}_{i}, \td{\mathsf{f}}_{i}$ for $i\in I_{m+n}$. 
\end{lem}
\pf (1) Suppose that $\mf{g}=\mf{c}$.  Let ${\bf T}=T_1\otimes \cdots \otimes T_\ell \in {\bf T}_{m+n}(\lambda,\ell)$ be given. Recall that each $T_k$ can be viewed as $T_k^{\tt R}\otimes T_k^{\tt L}$ as an element of a $\gl_{m+n}$-crystal. For $1\leq k\leq \ell$, let ${\bf T}^{(k)}$ be given by replacing $T_k^{\tt R}\otimes T_k^{\tt L}$ with $^{\tt R}T_k\otimes {}^{\tt L}T_k$ in ${\bf T}$. Since the column insertion is compatible with the $\gl_{m+n}$-crystal structure on $(\J_{m+n})^{\otimes r}$, the map sending $T_k^{\tt R}\otimes T_k^{\tt L}$ with $^{\tt R}T_k\otimes {}^{\tt L}T_k$ commutes with $\td{\mathsf{e}}_{i}, \td{\mathsf{f}}_{i}$  for $i\in I_{m+n}\setminus\{\ov{m}\}$, and hence so does the map sending ${\bf T}$ to ${\bf T}^{(k)}$.

Suppose that $\td{\mathsf f}_i {\bf T}\neq {\bf 0}$ for some $i\in I_{m+n}$ and
\begin{equation*}
\td{\mathsf f}_i {\bf T}=T_1\otimes \cdots \otimes \td{\mathsf f}_i T_k \otimes \cdots \otimes T_\ell,
\end{equation*}
for some $1\leq k\leq \ell$.
Let ${\bf S}=\td{\mathsf f}_i {\bf T}=S_1\otimes\cdots\otimes S_\ell$. It is clear that $S_j\in {\bf T}_{m+n}(\lambda'_j)$ for $1\leq j\leq \ell$, and $(S_{j+1},S_j)=(T_{j+1},T_j)$ is admissible for $1\leq j\leq \ell-1$ with $j\neq k-1, k$.\vskip 1mm

\textsc{Case 1}. Suppose that $i\in I_{m+n}\setminus\{\ov{m}\}$.
Since ${\bf S}^{(k)}=\td{\mathsf f}_i ({\bf T}^{(k)})$,  we have
\begin{itemize}
\item[(i)] $^{\tt R}S_k(l)\leq {}S^{\tt L}_{k-1}(l)$ for $l\geq 1$, when $k\geq2 $,

\item[(ii)] $S^{\tt R}_{k+1}(l+\lambda'_{k}-\lambda'_{k+1})\leq{}^{\tt L}S_k(l)$ for $l\geq 1$, when $k\leq \ell-1$.
\end{itemize}
Similarly, since ${\bf S}^{(j)}=\td{\mathsf f}_i ({\bf T}^{(j)})$ for $j=k-1,k+1$, we have
\begin{itemize}
\item[(i)] $^{\tt R}S_{k+1}(l)\leq{}S_k^{\tt L}(l)$ for $l\geq 1$, when $k\leq \ell-1$,

\item[(ii)] $S^{\tt R}_{k}(l+\lambda'_{k-1}-\lambda'_{k})\leq{}^{\tt L}S_{k-1}(l)$ for $l\geq 1$, when $k\geq 2$.
\end{itemize}
This implies that $S_{k+1}\prec S_k$ and $S_k \prec S_{k-1}$. Hence $\td{\mathsf f}_i {\bf T}\in {\bf T}_{m+n}(\lambda,\ell)$.\vskip 1mm

\textsc{Case 2}.  Suppose that $i=\ov{m}$. First, we claim that ${\rm ht}(S_k^{\tt L})\leq {\rm ht}(T_{k-1}^{\tt L})$. If ${\rm ht}(S_k^{\tt L})>{\rm ht}(T_{k-1}^{\tt L})$, then we have ${\rm ht}(T_k^{\tt L})= {\rm ht}(T_{k-1}^{\tt L})$ with $t_k^{\tt L},t_k^{\tt R}>\ov{m}$ and $t_{k-1}^{\tt L}=\ov{m}$ by tensor product rule, where $t_{j}^{\tt L}$ (resp. $t_j^{\tt R}$) is the top entry of $T_j^{\tt L}$ (resp. $T_j^{\tt R}$). Since $t_k^{\tt L}>\ov{m}$, we have $^{\tt R}t_k>\ov{m}=t_{k-1}^{\tt L}$ by definition of ${}^{\tt R}T_k$, where $^{\tt R}t_k$ is the top entry of $^{\tt R}T_k$. This contradicts the admissibility of $(T_{k},T_{k-1})$ and hence proves our claim.
Now, it is clear that $S_{j+1}\prec S_j$ for $1\leq j\leq \ell-1$ since $\ov{m}$ is the smallest one in $\J_{m+n}$. Hence $\td{\mathsf f}_{\ov{m}} {\bf T}\in {\bf T}_{m+n}(\lambda,\ell)$.

By \textsc{Case 1} and \textsc{Case 2}, we have $\td{\mathsf f}_i {\bf T}\in {\bf T}_{m+n}(\lambda,\ell)$. The proof of $\td{\mathsf e}_i {\bf T}\in {\bf T}_{m+n}(\lambda,\ell)\cup \{{\bf 0}\}$ is similar.

(2) Suppose that $\mf{g}=\mf{b}, \mf{b}^\bullet$. 
Let ${\bf T}=T_1\otimes \cdots \otimes T_{L} \in {\bf T}_{m+n}(\lambda,\ell)$ be given. Suppose that $\td{\mathsf f}_i {\bf T}\neq {\bf 0}$ for some $i\in I_{m+n}$ and
\begin{equation*}
\td{\mathsf f}_i {\bf T}=T_1\otimes \cdots \otimes \td{\mathsf f}_i T_k \otimes \cdots \otimes T_L,
\end{equation*}
for some $1\leq k\leq L$.
Let ${\bf S}=\td{\mathsf f}_i {\bf T}=S_1\otimes\cdots\otimes S_L$. Then  $(S_{j+1},S_j)=(T_{j+1},T_j)$ is admissible for $1\leq j\leq L-1$ with $j\neq k-1, k$. 

\textsc{Case 1}. Suppose that $i\in I_{m+n}\setminus\{\ov{m}\}$. By the same argument as in (1), we have $\td{\mathsf f}_i {\bf T}\in {\bf T}_{m+n}(\lambda,\ell)$.

\textsc{Case 2}. Suppose that $i=\ov{m}$. First, assume that 
$\resizebox{.035\hsize}{!}{$\boxed{\ov{m}}$}$ 
has been added on top of $T^{\tt L}_k$. In this case, $ {}^{\tt R}S_k={}^{\tt R}T_k $ and $^{\tt L}S_k$ is obtained by adding 
$\resizebox{.035\hsize}{!}{$\boxed{\ov{m}}$}$ 
on top of ${}^{\tt L}T_k$ by construction of $^{\tt R}T_k$ and $^{\tt L}T_k$ (see the proof of Lemma \ref{R matrix and signature-2}). Hence $S_{k+1}\prec S_k$ and $S_k \prec S_{k-1}$. Next, assume that 
$\resizebox{.035\hsize}{!}{$\boxed{\ov{m}}$}$
has been added on top of $T^{\tt R}_k$. In this case, ${}^{\tt L}S_k={}^{\tt L}T_k$ and $^{\tt R}S_k$ is obtained by adding 
$\resizebox{.035\hsize}{!}{$\boxed{\ov{m}}$}$
on top of ${}^{\tt R}T_k$. If ${\rm ht}(S_k^{\tt R})>{\rm ht}(T_{k-1}^{\tt L})$, then ${\rm ht}(T_k^{\tt R})={\rm ht}(T_{k-1}^{\tt L})$ and $t^{\tt L}_{k-1}=\ov{m}$. Since $T_k\prec T_{k-1}$, we have ${}^{\tt R}t_{k}= \ov{m}$ and hence either ${t}^{\tt L}_{k}= \ov{m}$ or ${t}^{\tt R}_{k}= \ov{m}$. But  this contradicts the fact that 
$\resizebox{.035\hsize}{!}{$\boxed{\ov{m}}$}$
can be added  on top of $T^{\tt R}_k$. So we have ${\rm ht}(S_k^{\tt R})\leq {\rm ht}(T_{k-1}^{\tt L})$. Hence $S_{k+1}\prec S_k$ and $S_k \prec S_{k-1}$.

By \textsc{Case 1} and \textsc{Case 2}, we have $\td{\mathsf f}_i {\bf T}\in {\bf T}_{m+n}(\lambda,\ell)$. The proof of $\td{\mathsf e}_i {\bf T}\in {\bf T}_{m+n}(\lambda,\ell)\cup \{{\bf 0}\}$ is similar.
\qed\vskip 2mm

\begin{lem}\label{connectedness of osp tableaux}
${\bf T}_{m+n}(\lambda,\ell)$ is a connected crystal with highest weight $\Lambda_{m+n}(\lambda,\ell)$.
\end{lem}
\pf (1) Suppose that $\mf{g}=\mf{c}$. Put ${\bf H}_{(\lambda,\ell)}=H_{(1^{\lambda'_1})}\otimes\cdots\otimes H_{(1^{\lambda'_\ell})}$. Then ${\rm wt}({\bf H}_{(\lambda,\ell)})=\Lambda_{m+n}(\lambda,\ell)$.
We will prove that any ${\bf T}=T_1\otimes\cdots\otimes T_\ell\in {\bf T}_{m+n}(\lambda,\ell)$ is connected to ${\bf H}_{(\lambda,\ell)}$ under  $\td{\mathsf e}_i$ for $i\in I_{m+n}$ by using  induction on $|{\bf T}|=\sum_{i=1}^\ell|T_i|$, the sum of the boxes in ${\bf T}$. We may assume that ${\bf T}$ is a $\gl_{m+n}$-highest weight element.

Choose the smallest $k\geq 1$ such that $T_k^{\tt R}$ is non-empty. Suppose that there is no such $k$. Since $T^{\tt R}_i=\emptyset$, $T_{i+1}\prec T_i$ and $(T_{i+1}\rightarrow(\cdots(T_2\rightarrow T_1)))$ is a $\mf{gl}_{m+n}$-highest weight element for $1\leq i\leq k-1$, we have $T_i=H_{(1^{\lambda'_i})}$ for $1\leq i\leq k$, that is, ${\bf T}={\bf H}_{(\lambda,\ell)}$.
If $k=1$, then $t^{\tt L}_1=t^{\tt R}_1=\ov{m}$. Otherwise, $\te_i {\bf T}\neq {\bf 0}$ for some $i\in I_{m+n}\setminus\{\ov{m}\}$. This implies that $\td{\mathsf e}_{\ov{m}}{\bf T}\neq {\bf 0}$. By induction hypothesis, $\td{\mathsf e}_{\ov{m}}{\bf T}$ is connected to ${\bf H}_{(\lambda,\ell)}$ and so is ${\bf T}$.

Suppose that $k\geq 2$. We have $T_i=H_{(1^{\lambda'_i})}$ for $1\leq i\leq k-1$.
Note that ${\rm ht}({}^{\tt R}T_k)=d$ for some $\lambda'_k< d\leq \lambda'_{k-1}$ and ${}^{\tt R}T_k(r)\leq T^{\tt L}_{k-1}(r)$ for $1\leq r\leq d$. Since $({}^{\tt R}T_k\rightarrow (T_{k-1}\rightarrow( \cdots  \rightarrow T_1)))$ is also a $\mf{gl}_{m+n}$-highest weight element, it is equal to $H_\mu$, where $\mu=(\lambda'_1,\ldots,\lambda'_{k-1},d)'$. In particular, we have ${}^{\tt R}T_k=H_{(1^d)}$.

We claim that $t^{\tt R}_k=\ov{m}$. Suppose that $t^{\tt R}_k>\ov{m}$. Since  ${}^{\tt R}T_k=H_{(1^d)}$, we have $d\geq 2$ and $\ov{m}<{}^{\tt L}t_k\leq  {}^{\tt R}t'_k$, where ${}^{\tt R}t'_k$ is the largest entry in $^{\tt R}T_k$.
Then $({}^{\tt L}t_k \rightarrow H_\mu)$ is of shape $\nu=(\lambda'_1,\ldots,\lambda'_{k-1},d,1)'$ but not a highest weight element $H_\nu$ since $\ov{m}<{}^{\tt L}t_k\leq  {}^{\tt R}t'_k$. This contradicts the fact that ${\bf T}$ is a $\gl_{m+n}$-highest weight element, and proves our claim, that is, $t^{\tt L}_k=t^{\tt R}_k=\ov{m}$. Hence $\td{\mathsf e}_{\ov{m}}{\bf T}\neq {\bf 0}$ and by induction hypothesis, ${\bf T}$ is connected to ${\bf H}_{(\lambda,
\ell)}$.

(2) Suppose that $\mf{g}=\mf{b}, \mf{b}^\bullet$. Put ${\bf H}_{(\lambda,\ell)}=H_{(1^{\lambda'_1})}\otimes\cdots\otimes H_{(1^{\lambda'_{L}})}$, where we assume that $\lambda'_L=0$ or the last tensor factor is empty tableau in ${\bf T}_{m+n}^{\rm sp}$ when $\ell-2\lambda_1$ is odd. Then ${\rm wt}({\bf H}_{(\lambda,\ell)})=\Lambda_{m+n}(\lambda,\ell)$, and we can show that  ${\bf T}\in {\bf T}_{m+n}(\lambda,\ell)$ is connected to ${\bf H}_{(\lambda,\ell)}$ under  $\td{\mathsf e}_i$ for $i\in I_{m+n}$ in almost the same way as in (1). \qed\vskip 2mm

\begin{thm}\label{character formula for m+n}
For $(\lambda,\ell)\in \cP({\mf g})_{m+n}$, ${\bf T}_{m+n}(\lambda,\ell)$ is isomorphic to the crystal of $L_q({\mf g}_{m+n},\Lambda_{m+n}(\lambda,\ell))$.
\end{thm}
\pf Suppose that $\mf{g}=\mf{b}, \mf{c}$. By Theorem \ref{regularity of level 1 osp crystal}, ${\bf T}_{m+n}(\lambda'_i)$ and ${\bf T}_{m+n}^{\rm sp}$ are  regular crystals and so is the crystal \eqref{tensor product of fundamental crystals}. By Lemma \ref{crystal invariance of osp tableaux}, ${\bf T}_{m+n}(\lambda,\ell)$ is a regular crystal. Hence by Lemma \ref{connectedness of osp tableaux}, it is isomorphic to the crystal of $L_q({\mf g}_{m+n},\Lambda_{m+n}(\lambda,\ell))$. Finally, the crystal of $L_q({\mf b}^\bullet_{m+n},\Lambda_{m+n}(\lambda,\ell))$ is isomorphic to that of $L_q({\mf b}_{m+n},\Lambda_{m+n}(\lambda,\ell))$ as an $I_{m+n}$-colored oriented graph \cite{BK}, and hence isomorphic to ${\bf T}_{m+n}(\lambda,\ell)$. \qed\vskip 2mm

Let  $\Z[P_{m+n}]$ be a group ring of $P_{m+n}$ with a $\Z$-basis  $\{\,e^\mu\,|\,\mu\in P_{m+n}\,\}$, and  ${\rm ch}L_q({\mf g}_{m+n},\Lambda)=\sum_{\mu\in P_{m+n}}\dim L_q({\mf g}_{m+n},\Lambda)_\mu e^\mu$ for $\Lambda\in P_{m+n}$. The character of a ${\mf g}_{m+n}$-crystal is defined in the same way.  Put $z=e^{\Lambda_{\ov{m}}}$ and $x_a=e^{\delta_a}$ for $a\in \J_{m+n}$.

\begin{cor}\label{character formula for m+n-2}
For $(\lambda,\ell)\in \cP({\mf g})_{m+n}$, we have $${\rm ch}L_q({\mf g}_{m+n},\Lambda_{m+n}(\lambda,\ell))=S^{\mf g}_{(\lambda,\ell)}(x_{_{\J_{m+n}}}).$$
\end{cor}

Now,  we have the following, which is the main result in this section.

\begin{thm}\label{character formula for m|n}
For $(\lambda,\ell)\in \cP({\mf g})_{m|n}$, we have $${\rm ch}L({\mf g}_{m|n},\Lambda_{m|n}(\lambda,\ell))=S^{\mf g}_{(\lambda,\ell)}(x_{_{\J_{m|n}}}).$$
That is, the weight generating function of ortho-symplectic tableaux of type $\mf{g}$ and shape $(\lambda,\ell)$ is the character of $L({\mf g}_{m|n},\Lambda_{m+n}(\lambda,\ell))$.
\end{thm}
\pf Considering the classical limit of $L_q({\mf g}_{m+\infty},\Lambda_{m+\infty}(\lambda,\ell))$, we have
\begin{equation*}
\begin{split}
{\rm ch}L({\mf g}_{m+\infty},\Lambda_{m+\infty}(\lambda,\ell))&=S^{\mf g}_{(\lambda,\ell)}(x_{_{\J_{m+\infty}}})=\sum_{\mu\in\cP}K^{\mf g}_{\mu\,(\lambda,\ell)}s_\mu(x_{_{\J_{m+\infty}}})
\end{split}
\end{equation*}
by Corollaries \ref{character formula for m+n-2} and \ref{Schur positivity-2}.
Hence by Theorem \ref{super duality} (see \eqref{T of K and L}), we have
\begin{equation*}
\begin{split}
{\rm ch}L({\mf g}_{m|\infty},\Lambda_{m|\infty}(\lambda,\ell))
&=\sum_{\mu\in\cP}K^{
\mf g}_{\mu\,(\lambda,\ell)}s_\mu(x_{_{\J_{m|\infty}}})=S^{\mf g}_{(\lambda,\ell)}(x_{_{\J_{m|\infty}}}).
\end{split}
\end{equation*}
Now ${\rm ch}L({\mf g}_{m|n},\Lambda_{m|n}(\lambda,\ell))$ is obtained by specializing $x_a=0$ for $a>n$, which is equal to $S^{\mf g}_{(\lambda,\ell)}(x_{_{\J_{m|n}}})$. \qed

\begin{rem}{\rm \mbox{}

(1) A Weyl-Kac type character formula for $L({\mf g}_{m|\infty},\Lambda_{m|\infty}(\lambda,\ell))$ can be obtained by super duality (see also \cite{CKW} for its detailed expression). 

(2) The coefficient $K^{\mf g}_{\mu\,(\lambda,\ell)}$ gives the branching multiplicity with respect to $U_q(\mf{gl}_{m|n})$-submodules. Even when $n=0$, our formula for branching multiplicity with respect to $\mf{gl}_m\subset \mf{b}_m, \mf{c}_m$ seems to be new.

(3) We may regard $S^{\mf g}_{(\lambda,\ell)}(x_{_{\J_{m|n}}})$ as a natural super-analogue of the irreducible characters over the classical Lie algebras or as an ortho-symplectic analogue of super Schur functions since it is obtained by superizing symmetric functions in $x_{_{\J_{m+\infty}}}$ with respect to $x_{_{\J_{m|\infty}}}$ and then specializing certain variables to 0. 

(4) We have more irreducible characters by $S^{\mf g}_{(\lambda,\ell)}(x_{\mc{A}})$ with other choices of $\mc{A}$. For example, when $\mc{A}=\td{\I}_m^+$ with $m\geq 0$, we have irreducible characters in $\td{\mc{O}}$ (see \eqref{duality functors}) corresponding to those in $\mc{O}^{int}(m+\infty)$ via the functor $T$.}
\end{rem}

\subsection{Connection with Kashiwara-Nakashima tableaux}\label{Connection with KN}
Let us briefly discuss how ortho-symplectic tableaux are related with Kashiwara-Nakashima tableaux (simply, KN tableaux) \cite{KashNaka}
when $\mf{g}=\mf{b}, \mf{c}$ and $n=0$.

Suppose that $\mf{g}=\mf{c}$. Let $T=(T^{\tt L}, T^{\tt R})\in {\bf T}_{m+0}(a)$ be given for $0\leq a\leq m$. Consider $({}^{\tt L}T,{}^{\tt R}T)$. Let $\sigma({}^{\tt R}T)$ be the single-columned tableau of height $m-{\rm ht}({}^{\tt R}T)$ with entries $\{1<\cdots < m\}$ such that $k$ appears in $\sigma({}^{\tt R}T)$ if and only if $\ov{k}$ does not appear in ${}^{\tt R}T$. Let $\td{T}$ be the tableau of height $m-a$ obtained by gluing ${}^{\tt L}T$ at the bottom of $\sigma({}^{\tt R}T)$. 
Then $\td{T}$ is a KN tableau of type $C_{m}$. For example, when $m=5$ and $a=1$, we have
$$T=(T^{\tt L}, T^{\tt R})\ =\ {\def\lr#1{\multicolumn{1}{|@{\hspace{.75ex}}c@{\hspace{.75ex}}|}{\raisebox{-.04ex}{$#1$}}}\raisebox{-.6ex}
{$\begin{array}{cc}
\cline{1-1}\cline{2-2}
\lr{\overline{5}}&\lr{\overline{4}}\\ 
\cline{1-1}\cline{2-2}
\lr{\overline{3}}&\lr{\overline{1}}\\ 
\cline{1-1}\cline{2-2}
\lr{\overline{2}}\\ 
\cline{1-1}  
\end{array}$}}\ \ \ \ \ \  
({}^{\tt L}T,{}^{\tt R}T)\ = \
{\def\lr#1{\multicolumn{1}{|@{\hspace{.75ex}}c@{\hspace{.75ex}}|}{\raisebox{-.04ex}{$#1$}}}\raisebox{-.6ex}
{$\begin{array}{cc}
\cline{1-1}\cline{2-2}
\lr{\overline{5}}&\lr{\overline{4}}\\ 
\cline{1-1}\cline{2-2}
\lr{\overline{2}}&\lr{\overline{3}}\\ 
\cline{1-1}\cline{2-2}
&\lr{\overline{1}}\\ 
\cline{2-2}  
\end{array}$}}\ \ \ \ \ \  \td{T} \ = \
{\def\lr#1{\multicolumn{1}{|@{\hspace{.75ex}}c@{\hspace{.75ex}}|}{\raisebox{-.04ex}{$#1$}}}\raisebox{-.6ex}
{$\begin{array}{c}
\cline{1-1} 
\lr{2} \\ 
\cline{1-1} 
\lr{5} \\ 
\cline{1-1}   
\lr{\overline{5}}\\ 
\cline{1-1}
\lr{\overline{2}}\\ 
\cline{1-1}  
\end{array}$}}
$$ \vskip 2mm
Now, for ${\bf T}=(T_\ell,\ldots, T_1)\in {\bf T}_{m+0}(\lambda,\ell)$, we have $\td{\bf T}=(\td{T}_\ell,\ldots,\td{T}_1)$, which forms a tableau of shape $\mu$, where $\mu=(\ell-\lambda_\ell,\ldots, \ell-\lambda_1)$ and $\td{T}_k$ is the $k$-th column from the right. Then from the admissibility of $(T_{k+1},T_k)$ for $1\leq k\leq \ell-1$, it follows that $\td{\bf T}$ is a KN tableau of type $C_m$.

Suppose that $\mf{g}=\mf{b}$. Let $T=(T^{\tt L}, T^{\tt R})\in {\bf T}_{m+0}(a)$ be given for $0\leq a\leq m$. 
In this case, we define $\td{T}$ as follows: Define $\sigma({}^{\tt R}T)$ in the same way as in $\mf{g}=\mf{c}$.  We place a  single-columned tableau of height $a+{\rm ht}(T^{\tt R})-{\rm ht}(T^{\tt L})$ with entries $0$, at the bottom of $\sigma({}^{\tt R}T)$ and then glue it with ${}^{\tt L}T$. 
Then $\td{T}$ is a KN tableau of type $B_{m}$ with ${\rm ht}(\td{T})=m-a$. For example, when $m=5$ and $a=1$, we have
$$T=(T^{\tt L}, T^{\tt R})\ =\ {\def\lr#1{\multicolumn{1}{|@{\hspace{.75ex}}c@{\hspace{.75ex}}|}{\raisebox{-.04ex}{$#1$}}}\raisebox{-.6ex}
{$\begin{array}{cc}
\cline{2-2}
&\lr{\overline{5}}\\ 
\cline{1-1}\cline{2-2}
\lr{\overline{5}}&\lr{\overline{4}}\\ 
\cline{1-1}\cline{2-2}
\lr{\overline{3}}&\lr{\overline{1}}\\ 
\cline{1-1}\cline{2-2}
\lr{\overline{1}}\\ 
\cline{1-1}  
\end{array}$}}\ \ \ \ \ \  
({}^{\tt L}T,{}^{\tt R}T)\ = \
{\def\lr#1{\multicolumn{1}{|@{\hspace{.75ex}}c@{\hspace{.75ex}}|}{\raisebox{-.04ex}{$#1$}}}\raisebox{-.6ex}
{$\begin{array}{cc}
\cline{2-2}
&\lr{\overline{5}}\\ 
\cline{1-1}\cline{2-2}
\lr{\overline{5}}&\lr{\overline{4}}\\ 
\cline{1-1}\cline{2-2}
\lr{\overline{1}}&\lr{\overline{3}}\\ 
\cline{1-1}\cline{2-2}
&\lr{\overline{1}}\\ 
\cline{2-2}  
\end{array}$}}\ \ \ \ \ \  \td{T} \ = \
{\def\lr#1{\multicolumn{1}{|@{\hspace{.75ex}}c@{\hspace{.75ex}}|}{\raisebox{-.04ex}{$#1$}}}\raisebox{-.6ex}
{$\begin{array}{c}
\cline{1-1} 
\lr{2} \\ 
\cline{1-1} 
\lr{0} \\ 
\cline{1-1}   
\lr{\overline{5}}\\ 
\cline{1-1}
\lr{\overline{1}}\\ 
\cline{1-1}  
\end{array}$}}
$$ \vskip 2mm
If $T\in {\bf T}^{\rm sp}_{m+0}$, then let $\sigma=(\sigma_1,\ldots,\sigma_m)$ be a sequence of $\pm$'s such that $\sigma_k=-$ if and only $\ov{k}$ appears in $T$. Then $\sigma$ determines a unique KN tableau of spin shape.
Now, we can recover KN tableaux of type $B_m$ from ${\bf T}_{m+0}(\lambda,\ell)$ in the same way as in $\mf{g}=\mf{c}$.

\section{Crystal base of a highest weight module in $\mc{O}_q^{int}(m|n)$}\label{Crystal base for m|n}
We assume that $\mf{g}=\mf{b}, \mf{b}^\bullet, \mf{c}$. In this section, we show that $L_q(\mf{g}_{m|n},\Lambda_{m|n}(\lambda,\ell))$ is an irreducible $U_q(\g_{m|n})$-module in $\mc{O}^{int}_q(m|n)$ and it has a unique crystal base for $(\lambda,\ell)\in\cP(\mf{g})_{m|n}$.

\subsection{Crystal structure on $\mathscr{V}_q$}\label{Crystal structure on V_q}
Let us consider the crystal structure on $\mathscr{V}_q$ in \eqref{Fock space}.
We identify $\J_{m|n}$ with the crystal of the natural representation of $U_q(\gl_{m|n})$, where
\begin{equation*}
\begin{split}
&\ov{m}\ \stackrel{^{\ov{m-1}}}{\longrightarrow}\ \ov{m-1}\ \stackrel{^{\ov{m-2}}}{\longrightarrow}
\cdots\stackrel{^{\ov{1}}}{\longrightarrow}\ \ov{1}\ \stackrel{^0}{\longrightarrow}\ \tfrac{1}{2} \
\stackrel{^{\hf}}{\longrightarrow}\ \tfrac{3}{2} \stackrel{{}^{\frac{3}{2}}}{\longrightarrow}\cdots \ \ \ \ \ \ \ \ \ \ \ \ \ \ \ \ \ (n=\infty),\\
&\ov{m}\ \stackrel{^{\ov{m-1}}}{\longrightarrow}\ \ov{m-1}\ \stackrel{^{\ov{m-2}}}{\longrightarrow}
\cdots\stackrel{^{\ov{1}}}{\longrightarrow}\ \ov{1}\ \stackrel{^0}{\longrightarrow}\ \tfrac{1}{2} \
\stackrel{^{\hf}}{\longrightarrow}\ \tfrac{3}{2} \stackrel{{}^{\frac{3}{2}}}{\longrightarrow}\cdots \stackrel{{}^{n-\frac{3}{2}}}{\longrightarrow} n-\tfrac{1}{2}\ \ \ \,(n<\infty),
\end{split}
\end{equation*}
where ${\rm wt}(a)=\delta_a$ for $a\in \J_{m|n}$ \cite{BKK}. The set of finite words with letters in $\J_{m|n}$
is the crystal of the tensor algebra generated by the natural representation of $U_q(\gl_{m|n})$, where
we identify each non-empty word $w=w_1\cdots w_r$ with $w_1\otimes\cdots\otimes w_r\in (\J_{m|n})^{\otimes r}$.

\begin{rem}\label{crystal for gl type-2}{\rm
Our convention for a crystal base of a $U_q(\gl_{m|n})$-module is different from the one in \cite{BKK}, where it is a lower crystal base as a $U_q(\mf{gl}_{m|0})$-module and a upper crystal base as a $U_q(\mf{gl}_{0|n})$-module (cf. Remark \ref{crystal for gl type}). But we have the same results as in \cite{BKK}.
Hence in our setting, for $w=w_1\cdots w_r\in  (\J_{m|n})^{\otimes r}$, $\te_i w$ and $\tf_i w$ ($i\in I_{0|n}$) are obtained by applying the usual signature rule (see Section \ref{signature rule of a tableau}) to $w$ for crystals for symmetrizable Kac-Moody algebras, while for $\te_i w$ and $\tf_i w$ ($i\in I_{m|0}\setminus\{\ov{m}\}$) we apply the signature rule to the reverse word $w^{\rm rev}$. Also,  for $i=0$, choose the largest $k$ ($1\leq k\leq
r$) such that $(\beta_0|{\rm wt}(w_k))\neq 0$. Then $\te_0 w$ (resp.
$\tf_0w$) is obtained by applying $\te_0$ (resp. $\tf_0$) to $w_k$. If
there is no such $k$, then $\te_0 w= {\bf 0}$ (resp. $\tf_0w={\bf
0}$).}
\end{rem}

Recall that for a skew Young diagram $\lambda/\mu$, $SST_{\J_{m|n}}(\lambda/\mu)$ is equipped with an $(I_{m|n}\setminus\{\ov{m}\})$-colored oriented graph structure (\cite[Theorem 4.4]{BKK}). Here $\te_i$ and $\tf_i$ are defined via the embedding $SST_{\J_{m|n}}(\lambda/\mu) \longrightarrow \bigsqcup_{r\geq 0}(\J_{m|n})^{\otimes r}$, which maps $T$ to $w^{\rm rev}(T)$ (see Remark \ref{crystal for gl type-2}).
Then for $\lambda\in\cP$, $SST_{\J_{m|n}}(\lambda)$ is isomorphic to the crystal of an irreducible polynomial representation of $U_q(\gl_{m|n})$ with highest weight $\Lambda_{m|n}(\lambda,0)\in P_{m|n}$ \cite[Theorem 5.1]{BKK}. We denote by $H^\natural_\lambda$ the highest weight element in $SST_{\J_{m|n}}(\lambda)$ with highest weight $\Lambda_{m|n}(\lambda,0)$, which is also called a genuine highest weight element \cite[Section 4.2]{BKK}.
\vskip 2mm

Let
\begin{equation*}
{\bf T}_{m|n} =\bigsqcup_{p,q\geq 0} SST_{\J_{m|n}}(1^{p})\times  SST_{\J_{m|n}}(1^{q}).
\end{equation*}
By identifying $(T^-,T^+)\in {\bf T}_{m|n}$ with $T^-\otimes T^+$ and applying the tensor product rule in Section \ref{crystal base}, we have an $(I_{m|n}\setminus\{\ov{m}\})$-colored oriented graph structure on ${\bf T}_{m|n}$.
Let $t^\pm$ be the top entries in $T^\pm$. If $t^-=t^+=\ov{m}$, then we define $\te_{\ov{m}}(T^-,T^+)=(S^-,S^+)$, where $S^\pm$ is obtained  by removing 
$\resizebox{.035\hsize}{!}{$\boxed{\ov{m}}$}$
from $T^\pm$. Otherwise, we define $\te_{\ov{m}}(T^-,T^+)={\bf 0}$. Also, if $t^-, t^+>\ov{m}$, then we define $\tf_{\ov{m}}(T^-,T^+)=(S^-,S^+)$, where $S^\pm$ is obtained from $T^\pm$ by adding 
$\resizebox{.035\hsize}{!}{$\boxed{\ov{m}}$}$
on top of $T^\pm$. Otherwise, we define $\tf_{\ov{m}}(T^-,T^+)={\bf 0}$. Hence ${\bf T}_{m|n}$ has an $I_{m|n}$-colored oriented graph structure.
We also have an $I_{m|n}$-colored oriented graph structure on
${\bf T}^{\rm sp}_{m|n} :={\bf T}^{\rm sp}_{\J_{m|n}}$,
where $\td{e}_{\ov{m}}$ and $\td{f}_{\ov{m}}$ are defined by removing and adding 
$\resizebox{.035\hsize}{!}{$\boxed{\ov{m}}$}$, respectively.

Let ${\bf m}=(m_a)\in \B$ be given. Put $d^\pm=\sum_{a\in\pm \J_{m|n}}m_a$. Let  $T^\pm({\bf m})\in SST_{\J_{m|n}}(1^{d^{\pm}})$ be the unique tableaux such that the number of occurrence of $a$ in $T^{\pm}({\bf m})$ is equal to $m_{\pm a}$ for $a\in \J_{m|n}$. 

If we regard $\B$ as a crystal of a $U_q(\mf{c}_{m|n})$-module $\F_q$, and $\B^+$ as a crystal of  a $U_q(\mf{b}_{m|n})$-module $\mathscr{F}^+_{q^2}$, then  it is straightforward to check the following (see Theorem \ref{crystal base of a Fock space} and its proof).

\begin{prop}\label{B=T} The maps
\begin{equation*}
\begin{split}
\Psi : \B \longrightarrow {\bf T}_{m|n}, \ \
\Psi^+ : \B^+ \longrightarrow {\bf T}^{\rm sp}_{m|n}
\end{split}
\end{equation*}
given by $\Psi({\bf m})=(T^-({\bf m}),T^+({\bf m}))$ and $\Psi^+({\bf m})=T^+({\bf m})$ are  bijections which commute with $\te_i$ and $\tf_i$ for $i\in I_{m|n}$.
\end{prop}

From now on,  we  regard ${\bf T}_{m|n}$ (resp. ${\bf T}_{m|n}^{\rm sp}$) as a crystal of $\F_q$ (resp. $\F^+_{q^2}$), where ${\rm wt}$, $\varepsilon_i$ and $\varphi_i$ ($i\in I_{m|n}$) are induced from those on $\B$ (resp. $\B^+$) via $\Psi$ (resp. $\Psi^+$).

\subsection{Highest weight crystal ${\bf T}_{m|n}(\lambda,\ell)$}
Put
\begin{equation*}
\begin{split}
{\bf T}_{m|n}(a)={\bf T}^{\mf{g}}_{\J_{m|n}}(a), \ \
{\bf T}_{m|n}(\lambda,\ell)={\bf T}^{\mf{g}}_{\J_{m|n}}(\lambda,\ell),
\end{split}
\end{equation*}
for $a\in\Z_{\geq 0}$  and $(\lambda,\ell)\in \cP(\mf{g})_{m|n}$. We regard
\begin{equation*}
\begin{cases}
{\bf T}_{m|n}(a)\subset {\bf T}_{m|n}, & \text{if $\mf{g}=\mf{c}$},\\
{\bf T}_{m|n}(a) \subset \left({\bf T}_{m|n}^{\rm sp}\right)^{\otimes 2}, & \text{if $\mf{g}=\mf{b}, \mf{b}^\bullet$},
\end{cases}
\end{equation*}
by identifying $T$ with $T^{\tt L}\otimes T^{\tt R}$, and apply $\td{e}_i$ and $\td{f}_i$ on ${\bf T}_{m|n}(a)$ for $i\in I_{m|n}$. Recall that in case of  ${\bf T}_{m+n}(a)$ in Section \ref{crystal structure for ortho-symplectic tableaux m+n}, we identify $T$ with $T^{\tt R}\otimes T^{\tt L}$. This difference is due to the tensor product rule for $\td{e}_i$ and $\td{f}_i$ ($i\in I_{m|n}$) (see Remark \ref{crystal for gl type-2}). But when $n=0$, ${\bf T}_{m+0}(a)$ with respect to $\td{\mathsf e}_i$ and $\td{\mathsf f}_i$  is isomorphic to ${\bf T}_{m|0}(a)$ with respect to $\te_i$ and $\tf_i$ for $i\in I_{m+0}=I_{m|0}$.

\begin{lem}\label{Invariance of T m|n (a)}
${\bf T}_{m|n}(a)\cup\{{\bf 0}\}$ is invariant under $\td{e}_i$ and $\td{f}_i$ for $i\in I_{m|n}$.
\end{lem}
\pf The proof is almost the same as in the case of ${\bf T}_{m+n}(a)$. So we leave the details to the reader.
\qed\vskip 2mm

Let $(\lambda,\ell)\in \cP({\mf g})_{m|n}$ be given with $L$ as in \eqref{Length of tuples}. We consider ${\bf T}_{m|n}(\lambda,\ell)$ as a subset of
\begin{equation}\label{tensor product of fundamental crystals-2}
\begin{cases}
{\bf T}^{\rm sp}_{m|n}\otimes{\bf T}_{m|n}(\lambda'_{L-1})\otimes\cdots \otimes {\bf T}_{m|n}(\lambda'_{1}), & \text{if $\mf{g}=\mf{b}$ and $\ell-2\lambda_1$ is odd},\\
{\bf T}_{m|n}(\lambda'_L)\otimes\cdots \otimes {\bf T}_{m|n}(\lambda'_1), & \text{otherwise},\\
\end{cases}\\
\end{equation}
by identifying ${\bf T}=(T_L,\ldots, T_1) \in {\bf T}_{m|n}(\lambda,\ell)$ with $T_L\otimes \cdots \otimes T_1$, and apply $\te_{i}$ and $\tf_{i}$ on ${\bf T}_{m|n}(\lambda,\ell)$ for $i\in I_{m|n}$. We put
\begin{equation}\label{Highest weight element in T m|n lambda}
{\bf H}^{\natural}_{(\lambda,\ell)}=
 H_{L}\otimes\cdots\otimes H_{1},
\end{equation}
where $H_k\in SST_{\J_{m|n}}(1^{\lambda'_k})$ ($1\leq k\leq L$) are the unique tableaux such that
$(H_L \rightarrow (\cdots(H_2\rightarrow H_1)))=H^\natural_\lambda$. We should remark that $H_k$ is not necessarily a highest weight element $H^\natural_{(1^{\lambda'_k})}$ in $SST_{\J_{m|n}}(1^{\lambda'_k})$. 
Here we assume that $H_L$ is the empty tableau in ${\bf T}^{\rm sp}_{m|n}$ if $\mf{g}=\mf{b}$ and $\ell-2\lambda_1$ is odd.

\begin{thm}\label{connectedness of T m|n (lambda,ell)} For  $(\lambda,\ell)\in \cP({\mf g})_{m|n}$, ${\bf T}_{m|n}(\lambda,\ell)\cup\{{\bf 0}\}$  is invariant under $\td{e}_i$ and $\td{f}_i$ for $i\in I_{m|n}$. Moreover, ${\bf T}_{m|n}(\lambda,\ell)$ is a connected $I_{m|n}$-colored oriented graph with a highest weight element ${\bf H}^{\natural}_{(\lambda,\ell)}$ of weight $\Lambda_{m|n}(\lambda,\ell)$.
\end{thm}
\pf Let ${\bf T}=T_L\otimes\cdots\otimes T_1 \in {\bf T}_{m|n}(\lambda,\ell)$ be given.  For $1\leq k\leq L$, we define ${\bf T}^{(k)}$ in the same way (see the proof of Lemma \ref{crystal invariance of osp tableaux} (1)).
Since the column insertion of $\J_{m|n}$-semistandard tableaux is compatible with the $\gl_{m|n}$-crystal structure (see \cite{KK}), we have $\td{f}_i ({\bf T}^{(k)})=(\td{ f}_i {\bf T})^{(k)}$ for $i\in I_{m|n}\setminus\{\ov{m}\}$. Then by the same argument as in Lemma \ref{crystal invariance of osp tableaux}, we can show that  ${\bf T}_{m|n}(\lambda,\ell)\cup\{{\bf 0}\}$  is invariant under $\td{e}_i$ and $\td{f}_i$ for $i\in I_{m|n}$.

Next, let us show that  ${\bf T}_{m|n}(\lambda,\ell)$ is connected. We will prove this only in the case of $\mf{g}=\mf{c}$, since the proof for $\mf{g}=\mf{b}, \mf{b}^\bullet$ is similar. By \cite[Theorem 4.8]{BKK}, we may assume that ${\bf T}$ is a $\mf{gl}_{m|n}$-highest weight element, that is, $(T_{\ell}\rightarrow ( \cdots( T_2 \rightarrow T_1)))$ is a genuine highest weight element. Choose the smallest $k\geq 1$ such that $T_k^{\tt R}$ is non-empty. If there is no such $k$, then ${\bf T}={\bf H}^\natural_{(\lambda,\ell)}$.
If $k=1$, then $t^{\tt L}_1=t^{\tt R}_1=\ov{m}$, which implies $\td{\mathsf e}_{\ov{m}}{\bf T}\neq {\bf 0}$. By induction on the number of boxes in ${\bf T}$, $\td{\mathsf e}_{\ov{m}}{\bf T}$ is connected to ${\bf H}^\natural_{(\lambda,\ell)}$ and so is ${\bf T}$.

Suppose that $k\geq 2$. Since $T^{\tt R}_i=\emptyset$, $T_{i+1}\prec T_i$ for $1\leq i\leq k-1$ and $(T_{\ell}\rightarrow(\cdots(T_2\rightarrow T_1)))$ is a $\mf{gl}_{m|n}$-highest weight element, we have $T_i(\lambda'_i-r+1)=\ov{m-r+1}$ for $1\leq i\leq k-1$ and $1\leq r\leq \min\{m,\lambda'_i\}$.
Since ${\rm ht}({}^{\tt R}T_k)=d$ for some $\lambda'_k< d\leq \lambda'_{k-1}$ and ${}^{\tt R}T_k(r)\leq T^{\tt L}_{k-1}(r)$ for $1\leq r\leq d$, we also have ${}^{\tt R}T_k(d-r+1)=\ov{m-r+1}$ for $1\leq r\leq \min\{d,m\}$. (Otherwise, $(T_{\ell}\rightarrow ( \cdots( T_2 \rightarrow T_1)))$ can't be a genuine highest weight element.)

Now, if $t^{\tt R}_k>\ov{m}$, then we have $d\geq 2$, $\ov{m}<{}^{\tt L}t_k\leq  {}^{\tt R}t'_k$, where ${}^{\tt R}t'_k$ is the largest entry in $^{\tt R}T_k$, and $S=({}^{\tt L}t_k \rightarrow ({}^{\tt R}T_k\rightarrow (T_{k-1}\rightarrow  \cdots (T_2  \rightarrow T_1))))$ is of shape $\nu=(\lambda'_1,\ldots,\lambda'_{k-1},d,1)'$. Since $\ov{m}<{}^{\tt L}t_k$, the entry in the right-most column of $S$ is not $\ov{m}$, which implies that $(T_{\ell}\rightarrow(\cdots(T_2\rightarrow T_1)))$ is not of the form $H^\natural_\eta$ for some $\eta$. This is a contradiction.
Therefore, $t^{\tt L}_k=t^{\tt R}_k=\ov{m}$, and $\td{\mathsf e}_{\ov{m}}{\bf T}\neq {\bf 0}$. By induction hypothesis, ${\bf T}$ is connected to ${\bf H}^\natural_{(\lambda, \ell)}$.
\qed

\subsection{Crystal base of a highest weight module}\label{Main Theorem}

First, suppose that $\mf{g}=\mf{c}$, and consider a crystal base $(\mathscr{L},\mathscr{B})$ of a $U_q(\mf{c}_{m|n})$-module $\mathscr{F}_q$ in Theorem \ref{crystal base of a Fock space}. For $a\geq 0$, let ${\bf m}(a)\in \B$ be such that $\Psi({\bf m}(a))=(H^{\natural}_{(1^a)},\emptyset)$. Then ${\bf v}_a:=\psi_{{\bf m}(a)}|0\rangle$ is a highest weight vector with highest weight $\Lambda_{m|n}((1^a),1)$. By Theorem \ref{highest weight module for integrable quantum super} and Corollary \ref{tensor power of Fock space}, we have
\begin{equation}\label{eq:aux-1}
U_q(\mf{c}_{m|n}) {\bf v}_a  \cong L_q(\mf{c}_{m|n},\Lambda_{m|n}((1^a),1))\in \mc{O}^{int}_q(m|n).
\end{equation}

\begin{prop}\label{Crystal base of a fundamental weight of type C} 
For $a\geq 0$, let
\begin{equation*}
\begin{split}
\mathscr{L}(a)&=\sum\mathbb{A}\td{x}_{i_1}\cdots\td{x}_{i_r}{\bf v}_a,\\
\mathscr{B}(a)&=\left\{\,\pm\td{x}_{i_1}\cdots\td{x}_{i_r}{\bf v}_a \!\!\!\pmod{q\mathscr{L}(a)} \, \right\}\setminus\{0\},
\end{split}
\end{equation*}
where $r\geq 0$, $i_1,\ldots,i_r\in I_{m|n}$, and $x=e, f$ for each $i_k$.
Then $(\mathscr{L}(a),\mathscr{B}(a))$ is a crystal base of $L_q(\mf{c}_{m|n},\Lambda_{m|n}((1^a),1))$, and the crystal $\mathscr{B}(a)/\{\pm1\}$ is isomorphic to ${\bf T}_{m|n}(a)$.
\end{prop}
\pf Since ${\bf v}_a\in\mathscr{L}$,  $\mathscr{L}(a)\subset \mathscr{L}$ and it is invariant under $\td{e}_i$ and $\td{f}_i$ for $i\in I_{m|n}$. Also, $\mathscr{B}(a)\subset \mathscr{B}$ and hence it is a pseudo-basis of $\mathscr{L}(a)/q\mathscr{L}(a)$ over  $\mathbb{Q}$ since $\mathscr{B}/\{\pm1\}$ is linearly independent.
By Proposition \ref{B=T} and Lemma \ref{Invariance of T m|n (a)}, the map
\begin{equation}\label{Map from B(a)}
\td{x}_{i_1}\cdots\td{x}_{i_r}{\bf v}_a\longmapsto \td{x}_{i_1}\cdots\td{x}_{i_r} H^{\natural}_{(1^a)}
\end{equation}
with $r\geq 0$ and $i_1,\ldots,i_r\in I_{m|n}$ is a well-defined weight preserving injection from $\mathscr{B}(a)/\{\pm1\}\cup\{0\}$ to ${\bf T}_{m|n}(a)\cup\{{\bf 0}\}$, which commutes with $\td{e}_i$ and $\td{f}_i$ for $i\in I_{m|n}$. By Theorem \ref{connectedness of T m|n (lambda,ell)}, ${\bf T}_{m|n}(a)$ is connected, and hence the map \eqref{Map from B(a)} is a bijection.
Now it follows from Theorem \ref{character formula for m|n} that ${\rm rank}_{\mathbb{A}}\mathscr{L}(a)_{\mu}=\dim L_q(\mf{c}_{m|n},\Lambda_{m|n}((1^a),1))_\mu$ for all weight $\mu$. This implies that $\mathscr{L}(a)$ is an $\mathbb{A}$-lattice of $L_q(\mf{c}_{m|n},\Lambda_{m|n}((1^a),1))$, and hence $(\mathscr{L}(a),\mathscr{B}(a))$ is its crystal base.
\qed \vskip 2mm

Next, suppose that $\mf{g}=\mf{b}, \mf{b}^\bullet$. Consider $\mathscr{F}^+_{q^2}\otimes \mathscr{F}^+_{q^2}$ as a $U_q(\mf{g}_{m|n})$-module.
For $a\geq 0$, let ${\bf m}^+(a)\in \B^+$ be such that $\Psi^+({\bf m}^+(a))=H^{\natural}_{(1^a)}$.

\begin{lem}\label{Highest weight vector for fundamental weight of type B}
For $a\geq 0$, there exists  ${\bf v}_{a}\in \mathscr{F}^+_{q^2}\otimes \mathscr{F}^+_{q^2}$ such that
\begin{itemize}
\item[(1)] ${\bf v}_a$ is a highest weight vector with highest weight $\Lambda_{m|n}((1^a),2)$,

\item[(2)] ${\bf v}_a\in \mathscr{L}^+\otimes \mathscr{L}^+$ and
${\bf v}_a \equiv \psi_{{\bf m}^+(a)}|0\rangle \otimes |0\rangle \pmod{q\mathscr{L}^+\otimes \mathscr{L}^+}$.
\end{itemize}
\end{lem}
\pf Let ${\bf m}^{(1)}, {\bf m}^{(2)}\in \B^+$ be given with ${\bf m}^{(s)}=(m_{rs})_{r\in \J_{m|n}}$. For convenience, we identify the $(\J_{m|n}\times 2)$-matrix $[m_{rs}]=[{\bf m}^{(1)} : {\bf m}^{(2)}]$ with $\psi_{{\bf m}^{(1)}}|0\rangle\otimes \psi_{{\bf m}^{(2)}}|0\rangle \in \F^+_{q^2}\otimes \F^+_{q^2}$. Note that $m_{rs}\in\{0,1\}$ when $|r|=0$, and $m_{rs}\in\Z_{\geq 0}$ when $|r|=1$.

Let $a\geq 0$ be given. Put $b=\max\{0,a-m\}$. Let ${\bf M}(a)$ be the set of non-negative integral $(\J_{m|n}\times 2)$-matrices $M=[m_{rs}]$ satisfying the following conditions:
\begin{itemize}
\item[(1)] $m_{rs}=0$ for $r> \hf$ and $s=1,2$,

\item[(2)] $m_{r 1}+m_{r 2}=1$ for $\ov{m}\leq r\leq \ov{l+1}$ where $l=\max\{m-a,0\}$,

\item[(3)] $m_{\hf 1} + m_{\hf 2}=b$.
\end{itemize}

Let $M=[m_{rs}]\in {\bf M}(a)$ be given. We write $M \stackrel{\ov{m}}{\rightsquigarrow} M'$ if $[m_{\ov{m}}\  m_{\ov{m}2}]=[1\ 0]$ and $M'$ is obtained from $M$ by replacing $[m_{\ov{m}}\  m_{\ov{m}2}]=[1\ 0]$ with $[0\ 1]$. For $i\in \{\ov{m-1},\ldots,\ov{1}\}$, we write $M \stackrel{i}{\rightsquigarrow} M'$ if $m_{\ov{i+1}1}=0$, $m_{\ov{i}1}=1$ and $M'$ is obtained from $M$ by replacing
\begin{equation*}
\begin{bmatrix}
m_{\ov{i+1}1} &  m_{\ov{i}2} \\
m_{\ov{i}1}  &  m_{\ov{i}2}
\end{bmatrix}
=\begin{bmatrix}
0 &  1 \\
1  & 0
\end{bmatrix}\ \ \text{with}\ \
\begin{bmatrix}
1 &  0 \\
0  & 1
\end{bmatrix}.
\end{equation*}
Similarly, we write $M \stackrel{0}{\rightsquigarrow} M'$ if $m_{\ov{1}1}=0$, $m_{\hf1}\geq 1$ and $M'$ is obtained from $M$ by replacing
\begin{equation*}
\begin{bmatrix}
m_{\ov{1}1} &  m_{\ov{1}2} \\
m_{\hf 1}  &  m_{\hf 2}
\end{bmatrix}
=\begin{bmatrix}
0 &  1 \\
u  & v
\end{bmatrix}\ \ \text{with}\ \
\begin{bmatrix}
1 &  0 \\
u-1  & v+1
\end{bmatrix}.
\end{equation*}
Identifying $M$ with $\psi_{{\bf m}^{(1)}}|0\rangle\otimes \psi_{{\bf m}^{(2)}}|0\rangle$, 
we have 
\begin{equation}\label{M and M'}
e_{i}M=Q_{M,M'}(q) e_{i}M',
\end{equation}
for $M \stackrel{i}{\rightsquigarrow} M'$, where $Q_{M,M'}(q)$ is a monomial in $q$ of positive degree given by
\begin{equation}
Q_{M,M'}(q)=
\begin{cases}
q, & \text{if $i=\ov{m}$ and $\mf{g}=\mf{b}$},\\
(-1)^{|{\rm wt}({\bf m}^{(1)})|+1}q, & \text{if $i=\ov{m}$ and $\mf{g}=\mf{b}^\bullet$},\\
q^2, & \text{if $i=\ov{m-1},\ldots,\ov{1}$},\\
(-1)^{|{\rm wt}({\bf m}^{(1)})|+1}q^{2\langle \beta^\vee_0,{\rm wt}({\bf m}^{(2)})\rangle}, & \text{if $i=0$}.\\
\end{cases}
\end{equation}

Let $M(a)\in {\bf M}(a)$ be such that $m_{r1}=1$ for $\ov{m}\leq r\leq \ov{l+1}$, and $m_{\hf 1}=b$.
Then for $M\in {\bf M}(a)$, we have
\begin{equation}\label{path from M(a) to M}
M(a)=M_0\stackrel{i_1}{\rightsquigarrow}M_1\stackrel{i_2}{\rightsquigarrow}\cdots \stackrel{i_{t}}{\rightsquigarrow}M_t=M,
\end{equation}
for some $t\geq 0$, $i_1,\ldots,i_t\in \{\ov{m},\ldots,\ov{1},0\}$ and $M_1,\ldots, M_{t-1}\in{\bf M}(a)$.
Put
\begin{equation*}
h(M)=t, \ \ \ Q_M(q)=\prod_{k=0}^{t-1}Q_{M_k,M_{k+1}}(q).
\end{equation*} 

Note that $M\in {\bf M}(a)$ is completely determined by its second column ${\bf m}^{(2)}$, and with this identification the $\{\ov{m},\ldots,\ov{1},0\}$-colored graph structure on ${\bf M}(a)$ with respect to $\ \stackrel{i}{\rightsquigarrow}\ $ coincides with the ${\mf b}_{m|1}$-crystal structure on ${\bf T}^{\rm sp}_{m|1}$ (see Section \ref{Crystal structure on V_q}). Then we can check without difficulty that $h(M)$ and $Q_M(q)$ are well defined, that is,  independent of a path  \eqref{path from M(a) to M} from $M(a)$ to $M$.

Now, we define
\begin{equation*}
{\bf v}_a=\sum_{M\in {\bf M}(a)}(-1)^{h(M)}Q_M(q)M.
\end{equation*}
Then ${\bf v}_a\in \mathscr{L}^+\otimes \mathscr{L}^+$ and ${\bf v}_a \equiv \psi_{{\bf m}^+(a)}|0\rangle \otimes |0\rangle \pmod{q\mathscr{L}^+\otimes \mathscr{L}^+}$. It remains to show that ${\bf v}_a$ is a highest weight vector, that is, $e_i {\bf v}_a =0$ for $i\in \{\ov{m},\ldots,\ov{1},0\}$.

Consider the pairs $(M,M')$ for $M, M'\in{\bf M}(a)$ such that $M\stackrel{i}{\rightsquigarrow} M'$. We see that any $M\in {\bf M}(a)$ with $e_{i}M\neq 0$ belongs to one of these pairs.
Since $h(M')=h(M)+1$ and $Q_{M'}(q)=Q_{M}(q)Q_{M,M'}(q)$, we have by \eqref{M and M'}
\begin{equation*}
\begin{split}
e_{i}&\left\{(-1)^{h(M)}Q_{M}(q)M+(-1)^{h(M')}Q_{M'}(q)M'\right\} \\ &= (-1)^{h(M)}Q_{M}(q)\left\{e_{i}M- Q_{M,M'}(q) e_{i}M' \right\}=0.
\end{split}
\end{equation*}
This implies that $e_{i}{\bf v}_a =0$. \qed\vskip 2mm

By Lemma \ref{Highest weight vector for fundamental weight of type B}, Theorem \ref{highest weight module for integrable quantum super} and Corollary \ref{tensor power of Fock space}, we have
\begin{equation}\label{eq:aux-2}
U_q(\mf{g}_{m|n}) {\bf v}_a  \cong L_q(\mf{g}_{m|n},\Lambda_{m|n}((1^a),2))\in \mc{O}^{int}_q(m|n).
\end{equation} 

\begin{prop}\label{Crystal base of a fundamental weight of type B} 
Suppose that $\mf{g}=\mf{b}, \mf{b}^\bullet$. For $a\geq 0$, let
\begin{equation*}
\begin{split}
\mathscr{L}(a)&=\sum\mathbb{A}\td{x}_{i_1}\cdots\td{x}_{i_r}{\bf v}_a,\\
\mathscr{B}(a)&=\{\,\pm\td{x}_{i_1}\cdots\td{x}_{i_r}{\bf v}_a \!\!\!\pmod{q\mathscr{L}(a)} \, \}\setminus\{0\},\\
\end{split}
\end{equation*}
where $r\geq 0$, $i_1,\ldots,i_r\in I_{m|n}$, and $x=e, f$ for each $i_k$.
Then $(\mathscr{L}(a),\mathscr{B}(a))$ is a crystal base of $L_q(\mf{g}_{m|n},\Lambda_{m|n}((1^a),2))$, and the crystal $\mathscr{B}(a)/\{\pm1\}$ is isomorphic to ${\bf T}_{m|n}(a)$ for $a\geq 0$.
\end{prop}
\pf  We note that $\mathscr{L}(a)\subset \mathscr{L}^+\otimes \mathscr{L}^+$ is invariant under $\td{e}_i$ and $\td{f}_i$ for $i\in I_{m|n}$, and $\mathscr{B}(a)\subset \mathscr{B}\otimes \mathscr{B}$ is a pseudo-basis of $\mathscr{L}^+(a)/q\mathscr{L}^+(a)$ over $\mathbb{Q}$. Then it follows from the same argument as in Proposition \ref{Crystal base of a fundamental weight of type C} that $(\mathscr{L}(a),\mathscr{B}(a))$ is a crystal base of $L_q(\mf{g}_{m|n},\Lambda_{m|n}((1^a),2))$.\qed\vskip 2mm

Now we are ready to state and prove our main theorem in this paper.
\begin{thm}\label{Existence of crystal base}
For $(\lambda,\ell)\in \cP({\mf g})_{m|n}$, $L_q(\mf{g}_{m|n},\Lambda_{m|n}(\lambda,\ell))$ is an irreducible $U_q(\g_{m|n})$-module in $\mc{O}^{int}_q(m|n)$, and it has a unique crystal base up to scalar multiplication, whose crystal is isomorphic to ${\bf T}_{m|n}(\lambda,\ell)$.
\end{thm}
\pf Let $(\lambda,\ell)\in \cP({\mf g})_{m|n}$ be given  with $L$ as in \eqref{Length of tuples}. 
For $1\leq i\leq L$, let $V_i$ be the $U_q(\mf{gl}_{m|n})$-submodule of  $\mathscr{V}_q$ or $\mathscr{V}_q^{\otimes 2}$ generated by ${\bf v}_{\lambda'_i}$ in \eqref{eq:aux-1} and \eqref{eq:aux-2} (we assume that ${\bf v}_{\lambda'_L}=|0\rangle$ in $\mathscr{F}^+_{q^2}$ if  $\mf{g}=\mf{b}$ with $\ell-2\lambda'_1$ odd). Then $V_i$ is isomorphic to the irreducible $U_q(\mf{gl}_{m|n})$-module with highest weight $\Lambda_{m|n}((1^{\lambda'_i}),r)$, where $r$ is either 1 or 2. Consider a $U_q(\mf{gl}_{m|n})$-module
$V_{(\lambda,\ell)}=V_{L}\otimes\cdots\otimes V_1$.
 Then it is completely reducible and has a crystal base \cite{BKK}. By Propositions \ref{Crystal base of a fundamental weight of type C}  and \ref{Crystal base of a fundamental weight of type B}, we may assume that the crystal lattice of $V_{(\lambda,\ell)}$ is contained in $\mathscr{L}_{(\lambda,\ell)}$, where $\mathscr{L}_{(\lambda,\ell)}$ is $\mathscr{L}^+ \otimes\mathscr{L}(\lambda'_{L-1})\otimes\cdots\otimes\mathscr{L}(\lambda'_1)$ if ${\mf g}={\mf b}$ with $\ell-2\lambda'_1$ odd, and  $\mathscr{L}(\lambda'_L)\otimes\cdots\otimes\mathscr{L}(\lambda'_1)$ otherwise. 
 
 By the decomposition of $V_{(\lambda,\ell)}$ into irreducible $U_q(\mf{gl}_{m|n})$-modules (see for example \cite[Example 5.8]{KK}), there exists a unique $U_q(\mf{gl}_{m|n})$-highest weight vector ${\bf v}_{(\lambda,\ell)}$ in $V_{(\lambda,\ell)}$ (up to scalar multiplication) such that
$U_q(\mf{gl}_{m|n}){\bf v}_{(\lambda,\ell)}$ is isomorphic to the irreducible $U_q(\mf{gl}_{m|n})$-module with highest weight $\Lambda_{m|n}(\lambda,\ell)$ and ${\bf v}_{(\lambda,\ell)}\not\equiv 0 \pmod{q\mathscr{L}_{(\lambda,\ell)}}$. Since ${\bf v}_{(\lambda,\ell)}\in V_{(\lambda,\ell)}=V_{L}\otimes \cdots\otimes V_1$ and $e_{\ov{m}}V_i=0$ for $1\leq i\leq L$, we have $e_{\ov{m}}{\bf v}_{(\lambda,\ell)}=0$. Hence ${\bf v}_{(\lambda,\ell)}$ is a $U_q(\mf{g}_{m|n})$-highest weight vector and
\begin{equation*}
U_q(\mf{g}_{m|n}){\bf v}_{(\lambda,\ell)} \cong L_q(\mf{g}_{m|n},\Lambda_{m|n}(\lambda,\ell))\in \mc{O}^{int}_q(m|n),
\end{equation*}
by Theorem \ref{highest weight module for integrable quantum super}, which also implies $L_q(\mf{g}_{m|n},\Lambda_{m|n}(\lambda,\ell))$ is a direct summand of $\mathscr{V}_q^{\otimes M}$ for some $M\geq 1$.
We also have
\begin{equation*}
{\bf v}_{(\lambda,\ell)}\equiv \pm {\bf H}^\natural_{(\lambda,\ell)} \pmod{q\mathscr{L}_{(\lambda,\ell)}}
\end{equation*}
(see \eqref{Highest weight element in T m|n lambda}) since a crystal base of $V_{(\lambda,\ell)}$ can be embedded into that of a tensor power of the natural representation of $U_q(\mf{gl}_{m|n})$, which is also a direct sum of crystal bases of irreducible polynomial $U_q(\mf{gl}_{m|n})$-modules (see \cite[Section 5]{BKK}).

Now, we  let
\begin{equation*}
\begin{split}
\mathscr{L}(\lambda,\ell)&=\sum\mathbb{A}\td{x}_{i_1}\cdots\td{x}_{i_r}{\bf v}_{(\lambda,\ell)},\\
\mathscr{B}(\lambda,\ell)&=\{\,\pm\td{x}_{i_1}\cdots\td{x}_{i_r}{\bf v}_{(\lambda,\ell)} \!\!\!\pmod{q\mathscr{L}(\lambda,\ell)}\, \}\setminus\{0\},
\end{split}
\end{equation*}
where $r\geq 0$, $i_1,\ldots,i_r\in I_{m|n}$, and $x=e, f$ for each $i_k$.

Since $\mathscr{L}(\lambda,\ell)\subset \mathscr{L}_{(\lambda,\ell)}$, $\mathscr{L}(\lambda,\ell)$ is invariant under $\td{e}_i$ and $\td{f}_i$ for $i\in I_{m|n}$, and $\mathscr{B}(\lambda,\ell)$ is a pseudo-basis of $\mathscr{L}(\lambda,\ell)/q\mathscr{L}(\lambda,\ell)$ over $\mathbb{Q}$.
By Theorem \ref{connectedness of T m|n (lambda,ell)}, Propositions \ref{Crystal base of a fundamental weight of type C} and \ref{Crystal base of a fundamental weight of type B}, the map
\begin{equation}\label{Map from B(lambda,ell)}
\Psi_{(\lambda,\ell)} : \mathscr{B}(\lambda,\ell)/\{\pm1\}\cup\{0\} \longrightarrow {\bf T}_{m|n}(\lambda,\ell)\cup\{{\bf 0}\},
\end{equation}
given by $\td{x}_{i_1}\cdots\td{x}_{i_r}{\bf v}_{(\lambda,\ell)}\longmapsto \td{x}_{i_1}\cdots\td{x}_{i_r} {\bf H}^\natural_{(\lambda,\ell)}$
with $r\geq 0$ and $i_1,\ldots,i_r\in I_{m|n}$ is a well defined bijection  (see \eqref{tensor product of fundamental crystals-2} and \eqref{Highest weight element in T m|n lambda}), which commutes with $\td{e}_i$ and $\td{f}_i$ for $i\in I_{m|n}$.

Finally, by Theorem \ref{character formula for m|n},  $\mathscr{L}(\lambda,\ell)$ is an $\mathbb{A}$-lattice of $L_q(\mf{g}_{m|n},\Lambda_{m|n}(\lambda,\ell))$, and hence $(\mathscr{L}(\lambda,\ell),\mathscr{B}(\lambda,\ell))$ is its crystal base. 

Finally, the uniqueness of a crystal base of $L_q(\mf{g}_{m|n},\Lambda_{m|n}(\lambda,\ell))$ follows from the connectedness of ${\bf T}_{m|n}(\lambda,\ell)$ and \cite[Lemma 2.7 (iii) and (iv)]{BKK}.\qed

\begin{cor}
The category $\mc{O}^{int}_q(m|n)$ is equivalent to $\mc{O}^{int}(m|n)$, and each $U_q(\mf{g}_{m|n})$-module in $\mc{O}^{int}_q(m|n)$ has a crystal base.
\end{cor}

\begin{cor}
A highest weight  $U_q(\mf{g}_{m|n})$-module in $\mc{O}^{int}_q(m|n)$ is a direct summand of $\mathscr{V}_q^{\otimes M}$ for some $M\geq 1$.
\end{cor}

\begin{rem}{\rm \mbox{}
There is another combinatorial character formula for the irreducible module $L_q(\mf{g}_{m|n},\Lambda_{m|n}(\lambda,\ell))$ ($\mf{g}=\mf{b}, \mf{b}^\bullet, \mf{c}$) given in terms of Young bitableaux \cite{K12}, where there is no crystal theory for superalgebras is used. It would be interesting to find a more explicit connection with ${\bf T}_{m|n}(\lambda,\ell)$. }
\end{rem}

\end{document}